 \magnification=\magstep1 
\newcount\sectno
\newcount\subsectno
\newcount\parno
\newcount\equationno
\newif\ifsubsections
\subsectionsfalse

\def\sectnum{\the\sectno} 
\def\subsectnum{\sectnum\ifsubsections .\the\subsectno\fi} 
\def\parnum{\subsectnum .{\the\parno}}
\def\eqnum{\subsectnum .\the\equationno}

\def\abstract#1\endabstract
{
{\abstractfont
    \baselineskip=9pt
    \leftskip=4pc  \rightskip=4pc
    \bigskip
    \noindent
     ABSTRACT.\ #1
\medskip} 
}
 
\def\thanks[#1]#2\endthanks{\footnote{$^#1$}{\footnotefont\kern-6pt #2}}

\newcount\minutes
\newcount\scratch

\def\timestamp{%
\scratch=\time
\divide\scratch by 60
\edef\hours{\the\scratch}
\multiply\scratch by 60
\minutes=\time
\advance\minutes by -\scratch
\the \month/\the\day$\,$---$\,$\hours:\null
\ifnum\minutes< 10 0\fi
\the\minutes}

\def\today{\ifcase\month\or
January\or February\or March\or April\or May\or June\or
July\or August\or September\or October\or November\or December\fi
\space\number\day,\number\year}

\outer\def\newsection #1.\par{\vskip1.5pc plus.75pc \penalty-250
     \subsectno=0
     \parno=0
     \equationno=0
     \advance\sectno by1
     \leftline{\smalltitlefont \sectnum.\hskip 1pc  #1}
                \nobreak \vskip.75pc\noindent}

\outer\def\newsectiontwoline #1/#2/.\par{\vskip1.5pc plus.75pc \penalty-250
     \subsectno=0
     \parno=0
     \equationno=0
     \advance\sectno by1
     \leftline{\smalltitlefont \sectnum.\hskip 1pc  #1}
     \leftline{\smalltitlefont \hskip 22pt #2}
                \nobreak \vskip.75pc\noindent}

\outer\def\newsubsection #1.\par{\vskip1pc plus.5pc\penalty-250
     \parno=0
     \equationno=0
     \advance\subsectno by1
     \leftline{{\bf \subsectnum}\hskip 1pc  #1.}
                \nobreak \vskip.5pc\noindent}

\def\newpar #1.{\advance \parno by1
     \par
 \medbreak \noindent 
      {\bf \parnum. #1.} \hskip 6pt}

\long\def \newclaim #1. #2\par {\advance \parno by1
    \medbreak \noindent 
     {\bf \parnum \hskip 6 pt #1.\hskip 6pt} {\sl #2} \par \medbreak}

\def\eq $$#1$${\global \advance \equationno by1 $$#1\eqno(\eqnum)$$}

\def\rmarginsert[#1]{\hglue 0pt\vadjust
{\null\vskip -\baselineskip\rightline{\abstractfont\rlap{\hfil\  #1}}}}

\def\lmarginsert[#1]{\hglue 0pt\vadjust
{\null\vskip -\baselineskip\leftline{\abstractfont\llap{#1\ \hfill}}}}

\newif\ifproofmode
\proofmodefalse

\def\refpar[#1]#2.{\advance \parno by1
     \par
 \medbreak \noindent 
      {\bf \parnum \hskip 6 pt #2.\hskip 6pt}%
\expandafter\edef\csname ref#1\endcsname
{\parnum}\ifproofmode\rmarginsert[\string\ref#1]\fi}

\long\def \refclaim[#1]#2. #3\par {\advance \parno by1
    \medbreak \noindent 
{\bf \parnum \hskip 6 pt #2.\hskip 6pt}%
\expandafter\edef\csname ref#1\endcsname
{\parnum}\ifproofmode\rmarginsert[\string\ref#1]\fi
{\sl #3} \par \medbreak}

\def\refer[#1]{%
\expandafter\xdef\csname ref#1\endcsname
{\parnum}\ifproofmode\rmarginsert[\string\ref#1]\fi}

\def\refereq[#1]$$#2$$ {%
\eq$$#2$$%
\expandafter\xdef\csname ref#1\endcsname{(\eqnum)}%
\ifproofmode\rmarginsert[\string\ref#1]\fi
}

\def\refeq{\refereq}

\def \Definition #1\\ {\vskip 1pc \noindent 
      {\bf #1. Definition. \hskip 6pt}\vskip 1pc}

\def\proof{{PROOF.} \enspace}

\def\qedmark{\hbox{\vrule height 4pt width 3pt}}
\def\qedskip{\vrule height 4pt width 0pt depth 1pc}
\def\qed{\penalty 1000\quad\penalty 1000{\qedmark\qedskip}}

\def \a {\alpha}
\def \b {\beta}

\def \D {\triangle}
\def \e {\epsilon}

\def \G {\Gamma}

\def \l {\lambda}
\def \L {\Lambda}
\def \n {\,\vert\,}
\def \N {\,\Vert\,}
\def \o {\theta}
\mathchardef\p="011E    

\def \s {\sigma}


\def\Gtwo{{\mathop{{{\mbi G\/}}\kern-.5pt_{{}_2}}}}
\def\Ffour{{\mathop{{{\mbi F\/}}\kern-2.5pt_{{}_4}}}}
\def\Esix{{\mathop{{{\mbi E\/}}\kern-.5pt_{{}_6}}}}
\def\Eseven{{\mathop{{{\mbi E\/}}\kern-.5pt_{{}_7}}}}
\def\Eeight{{\mathop{{{\mbi E\/}}\kern-.5pt_{{}_8}}}}


\def\ca{{\liefont A}}

\def\cb{{\liefont B}}

\def\cg{{\liefont G}}

\def\ck{{\liefont K}}

\def\cm{{\liefont M}}

\def\cn{{\liefont N}}

\def\co{{\liefont O}}

\def\cp{{\liefont P}}

\def\cs{{\liefont S}}

\def\ct{{\liefont T}}

\def\cu{{\liefont U}}

\def\cy{{\liefont Y}}

\def\sdp{
\mathop{\hbox{$\raise 1pt\hbox{$\scriptscriptstyle |$}\kern-2.5pt\times$}}
        }

\def\dag{{\raise 1 pt\hbox{{$\scriptscriptstyle \dagger$}}}}

\def\*{\raise 1.5pt \hbox{*}}

\def\csch{\mathop{\rm csch}\nolimits}
\def\sech{\mathop{\rm sech}\nolimits}
\def\Ker{\mathop{\rm Ker}\nolimits}

\def\Im{\mathop{\rm Im}\nolimits}

\def\exp{\mathop{\rm exp}\nolimits}

\def\tr{\mathop{\rm tr}\nolimits}

\def\diag{\mathop{\rm diag}\nolimits}

\def\Ad{\mathop{\rm Ad}\nolimits}
\def\ad{\mathop{\rm ad}\nolimits}

\def\diag {\mathop{\rm diag}\nolimits}

\def \li{\langle}
\def \ri{\rangle}

                                                          
\def\ifundefined#1{\expandafter\ifx\csname#1\endcsname\relax}

\def\cross{\times}
\def\longerrightarrow{-\kern-5pt\longrightarrow}

\def\star{\lower 1pt\hbox{*}}
\def \nulset {
\raise 1pt\hbox{
\hskip -3pt$\not$\kern -0.2pt \raise .7pt\hbox{${\scriptstyle\bigcirc}$}}}
\def \norm|#1|{\Vert#1\Vert}

\def \interior(#1){#1\kern -6pt \raise 7.5pt 
      \hbox{$\scriptstyle \circ$}{}\hskip 2pt}

\def\twist_#1{\kern -.15em\cross\kern -.30 em{}_{{}_{#1}}\kern .07 em}

\font\cmr=cmr10 at 10pt
\font\cmrviii=cmr8
\font\cmrvi=cmr6

\font\cmrXIV=cmr12 at 14 pt
\font\cmrXX=cmr17 at 20 pt
\font\cmrXXIV=cmr17 at 24 pt
\font\cmbxXII=cmbx12
\font\cmbxsl=cmbxsl10
\font\cmbxslviii=cmbxsl10 at 8pt
\font\cmbxslv=cmbxsl10 at 5pt

      \font\tenrm=cmr10 at 10.3 pt
      \font\sevenrm=cmr7 at 7.21 pt
      \font\fiverm=cmr5 at 5.15 pt
      \font\teni=cmmi10 at 10.3 pt
      \font\seveni=cmmi7 at 7.21 pt
      \font\fivei=cmmi5 at 5.15 pt     
      \font\tensy=cmsy10 at 10.3 pt
      \font\sevensy=cmsy7 at 7.21 pt   
      \font\fivesy=cmsy5 at 5.15 pt    
      \font\tenex=cmex10 at 10.3 pt
      \font\tenbf=cmbx10 at 10.3 pt
      \font\sevenbf=cmbx7 at 7.21 pt    
      \font\fivebf=cmbx5 at 5.15 pt  

\def\UseComputerModern   
{
\textfont0=\tenrm \scriptfont0=\sevenrm \scriptscriptfont0=\fiverm
\def\rm{\fam0\tenrm}
\textfont1=\teni \scriptfont1=\seveni \scriptscriptfont1=\fivei
\def\mit{\fam1} \def\oldstyle{\fam1\teni}
\textfont2=\tensy \scriptfont2=\sevensy \scriptscriptfont2=\fivesy
\def\cal{\fam2}
\textfont3=\tenex \scriptfont3=\tenex \scriptscriptfont3=\tenex
\def\it{\fam\itfam\tenit} 
\textfont\itfam=\tenit
\def\sl{\fam\slfam\tensl} 
\textfont\slfam=\tensl
\def\bf{\fam\bffam\tenbf} 
\textfont\bffam=\tenbf \scriptfont\bffam=\sevenbf
\scriptscriptfont\bffam=\fivebf
\def\tt{\fam\ttfam\tentt} 
\textfont\ttfam=\tentt
\def\abstractfont{\cmrviii}
\def\footnotefont{\cmrviii}
\def\tinyfont{\cmrvi}
\def\smalltitlefont{\cmbxXII}  
\def\titlefont{\cmrXIV}
\def\bigtitlefont{\cmrXX}
\def\verybigtitlefont{\cmrXXIV} 
\textfont9=\cmbxsl \scriptfont9=\cmbxslviii \scriptscriptfont9=\cmbxslv
\def\mbi{\fam9}
\cmr
}  

\def\liefont{\cal}

\def \bs {\bigskip}
\def \ms {\medskip}
\def \ss {\smallskip}

\def \ni {\noindent}

\def\enditem{\item{}\par\vskip-\baselineskip\noindent}
\def\ei{\enditem}

   \baselineskip=14 true pt 
   \hsize 37 true pc \hoffset= 22 true pt
   \voffset= 0 true pt
   \vsize  54 true pc

\def\ni{\noindent}
\def\bs{\bigskip}
\def\ms{\medskip}

\def\ss{\smallskip}
\def\ct{{\cal T}}

   \baselineskip=14 true pt 
   \hsize 35 true pc \hoffset= 25 true pt
   \vsize  52 true pc
\UseComputerModern
\subsectionsfalse
\font\smalltitlefont=cmbx10 at 11 pt
\proofmodetrue
%


\def\Bibliography
{

\font\TRten=cmr10 at 10 true pt
\font\TIten=cmti10 at 10 true pt
\font\TBten=cmbx10 at 10 true pt

\def\ourindent{\hfil\vskip-\baselineskip}

 \frenchspacing
 \parindent=0pt

 \def \keyfnt{\TRten}
 \def \authornamefnt{\TRten}
 \def \booktitlefnt{\TIten}
 \def \articletitlefnt{\TRten}
 \def \journalnamefnt{\TIten}
 \def \volumefnt{\TBten}
 \def \publishernamefnt{\TRten}
 \def \pagesfnt{\TRten}
 \def \yearfnt{\TRten}
 \def \commentfnt{\TRten}

 \def \bookitem //##1//##2//##3//##4//##5//##6//##7//##8//
      { \goodbreak{\par\hskip-40pt{\keyfnt [##1]}\ourindent{\authornamefnt ##2,}}
              {\booktitlefnt ##3.\/}\thinspace
              {\publishernamefnt ##4,}
              {\yearfnt ##6.}
              {\commentfnt ##8}
       }

\def \b{\bookitem}

\def \articleitem //##1//##2//##3//##4//##5//##6//##7//##8//
      { \goodbreak{\par\hskip-40pt{\keyfnt [##1]}\ourindent{\authornamefnt ##2,}}
              {\articletitlefnt ##3},
              {\journalnamefnt ##4\/}
              {\volumefnt ##5}
              {\hbox{\yearfnt(\hskip -1pt ##6)}},
              {\pagesfnt ##7.}
              {\commentfnt ##8}
       }
\def \a{\articleitem}

\def \preprintitem //##1//##2//##3//##4//##5//##6//##7//##8//
      { \goodbreak{\par\hskip-40pt{\keyfnt [##1]}\ourindent{\authornamefnt ##2,}}
              {\articletitlefnt ##3},
              {\commentfnt ##8}
       }
\def \p{\preprintitem}

   \vskip 1in
   \centerline{References}
   \vskip .5in
}


\def\p{\partial}
\font\bfs=cmb10 at 9 pt

\proofmodefalse

\font\cmbxXIV=cmbx12 at 14pt

\centerline {\cmbxXIV B\"acklund Transformations and Loop Group Actions\/}
\bs\bs
\centerline { Chuu-Lian Terng\footnote{$^1$} {Research supported in part by 
NSF Grant DMS 9626130 and Humboldt Senior Scientist Award}
and  Karen Uhlenbeck\footnote{$^2$}{Research supported in part by Sid Richardson 
Regents' Chair Funds, University of Texas system}}

\font\bfs=cmb10 at 10 pt
\font\cmrviii=cmr8

\bs\bs\bs
\centerline {\bfs Abstract\/}
\bs\bs
{\cmrviii  We construct a local action of the group of rational maps from $S^2$
to $GL(n,C)$ on local solutions of flows of the ZS-AKNS $sl(n,C)$-hierarchy. We show
that the actions of simple elements (linear fractional transformations) give local
B\"acklund transformations, and we derive a permutability formula from different
factorizations of a quadratic element.  We prove that the action of simple
elements on the vacuum may give either global smooth solutions or solutions
with singularities. However, the action of the subgroup of the rational maps
that satisfy the $U(n)$-reality condition $g(\bar\l)^*g(\l)=I$ on the space of
global rapidly decaying  solutions of the flows in the 
$u(n)$-hierarchy is global, and the action of a simple element gives a global
B\"acklund transformation. 
 The actions of certain 
elements in the rational loop group on the vacuum give rise to explicit time
periodic multi-solitons  (multi-breathers).  We show that this theory
generalizes the classical B\"acklund theory of the sine-Gordon equation. 
The group structures of B\"acklund transformations for various hierarchies are
determined by their reality conditions. We
identify the reality conditions (the group structures) for the $sl(n,R)$,
$u(k,n-k)$,  KdV, Kupershmidt-Wilson, and Gel'fand-Dikii hierarchies.   The
actions of linear fractional transformations that satisfies a reality condition,
modulo the center of the group of rational maps, gives  B\"acklund and Darboux
transformations for the hierarchy defined by the reality condition.    Since the
factorization  cannot always be carried out under these reality condition, the
action is again local, and B\"acklund transformations only generate local solutions
for these hierarchies unless singular solutions are allowed.}

\vfil\eject

\centerline {\bf Table of Contents\/}
\ms

\+\kern .6in &\kern .3in & \cr
\+ & 1. & Introduction\cr
\+ & 2. & The ZS-AKNS $n\times n$ flows \cr
\+ & 3. & Reality conditions \cr
\+ & 4. & B\"acklund transformations for the $sl(n,C)$-hierarchy  \cr
\+ & 5.  &  B\"acklund transformations for the $u(n)$-hierarchy\cr
\+ & 6.  &  Permutability formula\cr
\+ & 7.  & $N$-soliton formula\cr
\+ & 8. & Scaling transformations \cr
\+ & 9. & B\"acklund transformations for $n$-dimensional systems \cr
\+ & 10. & Time periodic solutions \cr
\+ & 11. &  The $u(k,n-k)$-hierarchies\cr
\+ &12.  & The KdV hierarchy \cr
\+ &13.  & The Kupershmidt-Wilson hierarchy \cr
\+ &14.  & The Gel'fand-Dikii hierarchy \cr

\bs

\newsection Introduction.\par

The classical B\"acklund transformations are local geometric transformations, which
construct from a given surface of constant Gaussian curvature $-1$ a 
$2$-parameter family of such surfaces.  To find such transformations, one needs to solve a
system of compatible ordinary differential equations.  Since surfaces of Gaussian
curvature
$-1$ are classically known to be equivalent to local solutions of the sine-Gordon equation
$$q_{xt}=\sin q, \eqno(SGE)$$ this provides a method of deriving new solutions of a partial
differential equation from a given solution via the solution of ordinary differential equations
(cf. [Da], [Ei]).  Most of the known ``integrable systems'' possess transformations
of this type.

Applying B\"acklund transformations $n$ times to a
 solution of the  sine-Gordon equation produces a hierarchy of $2n$-dimensional
families of solutions. Moreover, the Bianchi permutability theorem
states that the second and higher families can be obtained
from the first family through algebraic formulas. This allowed the classical
geometers to write down explicit solutions for the sine-Gordon equation  and
explicit surfaces of curvature 
$-1$. For example, they applied one B\"acklund transformation to the vacuum solution of the
sine-Gordon equation to get the pseudosphere (stationary $1$-soliton) and  Dini
surfaces ($1$-soliton), applied B\"acklund transformations twice to the vacuum
solution to get the  K\"uen surface ($2$-soliton), and applied the Bianchi
permutability formula with two suitable complex conjugate parameters to get
breathers (time periodic solutions). 

The sine-Gordon equation  is clearly invariant under the Lorentz group, i.e., if $q$ is a
solution of the sine-Gordon equation and $r$ is a non-zero real number,  then $\hat
q(x,t)=q(r^{-1}x, rt)$ is again a solution of the sine-Gordon equation (note that we
are working in characteristic coordinates). This is called the Lie transformation for
$-1$ curvature surfaces in
$R^3$ in the classical surface theory.  

What are now called Darboux transformations  were discovered by
Darboux during his investigation of Liouville metrics. A metric
$ds^2= A(x,y)(dx^2+dy^2)$ is Liouville if there is a coordinate system
$(u,v)$ such that $ds^2$ is of the form
$$ds^2= (f(u)- g(v))(du^2+dv^2)$$ for some $f$ and $g$ of one variable.  The
classical geometers were interested in such metrics at least in part
because Liouville
had shown that all geodesics on such surfaces can be obtained by quadratures.
The question of deciding whether a metric  $ds^2$ is Liouville
led to the  study of the following special second order 
linear partial differential equation
$$w_{xy}= (f(x+y)-g(x-y))w.$$ Darboux was led to look for transformations 
of Hill's operators in the process of separating variables in this
equation.     The original 
analytic version of Darboux transformation ([Da] v. 2 Chap. 9) is the
following: Let
$q$ be a smooth function of one variable, $k_0$ a constant, and
suppose that $f$ satisfies $f'' =
(q+k_0)f$.  Set $$q^\sharp = f(f^{-1})'' - k_0.$$ If $y(x,k)$ is the
general
solution of the Hills operator with potential $q$:
$$y''=(q+k)y,$$ then $z=y'-(f/f') y$ is the general solution of the Hills
operator
with potential $q^\sharp$:
$$z''= (q^\sharp + k) z.$$
This Darboux theorem gives an algebraic algorithm (without quadrature) to
transform general solutions of $D^2-q-k$ to those of $D^2-q^\sharp-k$.
Next, suppose that we factor
$$D^2 -q-\l_0 = (D+v)(D-v),$$ In other words, suppose that $v$ 
satisfies $v_x+v^2 =q+\l_0$.
(Here $D={d\over dx}$.) Choose $f$ so that $f'/f= v$.  Then
$$(D-v)(D+v) = D^2-q^\sharp -\l_0.$$
Since, if $q(x,t)$ is a solution of KdV then the Hills operators with potential
$q(\cdot, t)$ are isospectral,   it follows that the Darboux transformations of 
the Hills operators
induce transformations on the space of solutions of KdV.  This is a
critical observation due to Adler and Moser [AM] and Deift [De].

We give another interpretation of the Darboux transformation. Write the Hills operator
${d^2\over dx^2}-q-\l^2$ as a first order system
$L_{q,\l}={d\over dx} - \pmatrix{\l& q\cr 1 & -\l\cr}$. Given an eigenfunction of the Hills
operator with potential $q$ and eigenvalue
$k_0=\a^2$ is the same as given a trivialization of
$L_{q,\a}$ (thought as a connection on the line).  Then Darboux' theorem can be
reformulated as follows: Given a trivialization of $L_{q,\a}$ for some $\a\in C$,
the map $q\mapsto q^\sharp$ transforms the trivialization of $L_{q,\l}$ to
that of  $L_{q^\sharp,\l}$ by an algebraic formula.  Most of the known integrable systems 
also possess transformations of this type.

The second author constructed an action of the rational loop group on the space of
solutions of harmonic maps from $R^{1,1}$ to $SU(n)$ in [U1], and showed that the action
of a simple element (i.e., a linear fractional transformation) can be obtained by
solving two compatible ordinary differential equations. The starting point of this
paper is the realization that this rational loop group action  in [U1]
can be generalized to solutions of other  partial differential equations having Lax pairs and
that satisfy the ``reality condition'' of a compact group. We give an analogous
construction of the action of the rational loop group on the space of global solutions of
the flows in the AKNS-ZS $u(n)$-hierarchy.  We will see:
\item {(1)} The action of a simple element (a linear fractional transformation) corresponds to
a global B\"acklund transformation.
\item {(2)} The Bianchi permutability formula arises from various ways of
factoring quadratic elements in the rational loop group into simple elements.
\item {(3)} The B\"acklund transformations can be computed from solutions of
ordinary differential equations given a known solution of the partial differential
equation. 
\item {(4)} Once given the trivialization of the Lax pair corresponding to a given
solution, the action of a simple element corresponds to a global Darboux transformation and is 
algebraic. 
\item {(5)} Lie transformations arise as the scaling transformations, which extend the action
of the rational loop group to the semi-direct product of the multiplicative group $R^*$ of
non-zero real numbers and the rational loop group. 

\ms

Since the sine-Gordon equation arises as part of the hierarchy (the $-1$-flow for
$su(2)$ with an involution constraint), we can check that we are generalizing the classical
theory. The choice of group structure depends on the choice of base point.
Hence the group structure is not canonical and was not apparent to the classical
geometers. 

An interesting observations is that appropriate choices of poles for
the rational loop yield time periodic solutions. This gives an insight
into the construction of the classical breathers of the
sine-Gordon equation ([Da]).  There are no
simple factors in the rational loop group corresponding to the placement of poles for time
periodic solutions. However, there are quadratic elements (product of two simple
elements), whose simple factors do not satisfy the algebraic constraints to
preserve sine-Gordon, but which nevertheless generate the well-known breathers
(one way to think of them is as the product of two complex conjugate B\"acklund
transformations). The product of these quadratic factors generate arbitrarily
complicated time periodic solutions. 

The sine-Gordon equation also arises as the equation for wave (or harmonic) maps from
the Lorentz space $R^{1,1}$ to $S^2$ (for example, see [P]).  Shatah and Strauss proved
in [SS] that the classical breather solutions for the sine-Gordon equation produce
homoclinic wave maps from
$S^1\times R$ to $S^2$.  Using a simple change of gauge for the Lax pair of the $-1$-flow,
the first author proved in [Te] that solutions of the $-1$-flow give rises to wave maps
from
$R^{1,1}$ to symmetric spaces.  In a forthcoming paper [TU2], we prove that the 
time periodic
$m$-solitons for the  $-1$-flow constructed in this paper also give rise
to homoclinic wave maps from $S^1\times R$ into compact symmetric spaces.   

The permutability formula has several useful applications. For example, 
one of the key ingredients of the study of discrete $-1$ curvature surfaces in $R^3$  by
Bobenko and Pinkall [BP] is the permutability formula for the sine-Gordon equation. Since
constant sectional curvature $n$-dimensional submanifolds in Euclidean spaces are given by
solutions of the $n$ commuting first flows (cf. [Te], [TU1]),  the generalized
permutability formula  should be useful in the study of discretization of constant
curvature submanifolds and soliton equations.  We also use the Bianchi
permutability formula to write down an explicit formula for  $m$-soliton solutions
of the $j$-th flow. 

 Local B\"acklund transformations for the
$j$-th flow were constructed by Zakharov and Shabat in [ZS 2], Sattinger and Zurkowski  in
[SZ 1, 2], by Beals, Deift and Tomei in [BDT], by Gu and Zhou in [GZ] and by Cherndik
in [Ch].  Our construction gives a group structure of these transformations,  and
provide a systematic method of finding such transformations for equations having
a Lax pair. The algebraic structure of these transformations also makes many of
the mysterious classical results for the sine-Gordon equation apparent.   

We give an outline of the method we use to construct B\"acklund transformations and explain
how the group structure for these transformations is obtained. Most of the evolution we
considered in this paper has a Lax pair with a parameter, i.e., it is given as the condition that a
one-parameter family of connections is flat:
$$\left[{\p\over \p x} + A(x,t,\l), {\p\over \p t} + B(x,t,\l)\right]=0,$$ where $A$ and $B$
are differential operators in $u$ and its derivatives in $x$ and $A,B$ are holomorphic for $\l\in
C$.  The trivialization $E$ of a solution $u$ normalized at $(0,0)$  is
defined to be the solution for the following linear system
$$E_x = E A, \quad E_t= EB, \quad E(0,0,\l)=I.$$ Then $E(x,t,\l)$ is holomorphic for $\l\in
C$. The general view of the specific construction is the Birkhoff factorization
theorem.  Let
$\co_\infty$ denote a neighborhood of $\infty$ in $C\cup\{\infty\}=S^2$, $L_+(GL(n,C))$
the group of holomorphic maps from $C$ to $GL(n,C)$, (under pointwise
multiplication), $L_-(GL(n,C))$ the group of holomorphic maps $h_-$ from $\co_\infty$ to
$GL(n,C)$ such that $h_-(\infty)=I$, and $L(GL(n,C))$ the group of holomorphic maps from
$\co_\infty\cap C$ to
$GL(n,C)$.  The Birkhoff factorization theorem states that the multiplication map
$$\mu: L_+(GL(n,C))\times L_-(GL(n,C))\to L(GL(n,C)), \quad (h_+, h_-)\mapsto h_+h_-$$
is one to one and the image is an open dense subset of $L(GL(n,C))$.   Hence, formally, there
is a ``dressing action'' of $L_-(GL(n,C))$ on $L_+(GL(n,C))$ defined as follows: given
$h_\pm\in L_\pm(GL(n,C))$, if $h_-h_+$ lies in the image of $\mu$ then there exists
unique
$f_\pm\in L_\pm(GL(n,C))$ such that $h_-h_+= f_+ f_-$.  Then the dressing action is
defined by
$h_-\sharp h_+= f_+$.  The Birkhoff factorization is not explicit.  Moreover,
singularities arise, typically on a codimension two set in the parameter spaces. 
Because the singular set is closed and its complement is dense, the action is local. 
 
To construct B\"acklund transformation, we choose a
linear fractional transformation $h_-\in L_-(GL(n,C))$. Since the trivialization
$E(x,t, \l)$ of a solution $u$ of the evolution equation is holomorphic for $\l\in C$, the map
$E(x,t)\in L_+(GL(n,C))$, where $E(x,t)(\l)= E(x,t,\l)$.  For each $(x,t)$, let
$h_-$ acts on  $E(x,t)$.  If 
$\tilde E(x,t)= h_-\sharp E(x,t)$ exists, then a new solution can be obtained from $\tilde
E^{-1}\tilde E_x$.  However, the factorization only involves poles and zeros, and is
explicit and algebraic, not abstract.  The group structure of these transformations
is clearly the one inherited from
$L_-(GL(n,C))$.  However, even these factorizations can not be always carried out. 
Hence, in general, it gives rise to a local theory. The only case we obtain a
good global theory is when the Lax pair satisfies the $u(n)$-reality condition, i.e., $A$ and
$B$ satisfies:
$$A(x,t, \bar \l)^* + A(x,t,\l)=0, \quad B(x,t,\bar\l)^* + B(x,t,\l)=0.$$    

Let $\cg$, $\cg_\pm$ denote the Lie algebra of $L(GL(n,C))$ and
$L_\pm(GL(n,C))$ respectively.  It was explained in an earlier paper [TU1] that $(\cg, \cg_+,
\cg_-)$ is a Manin-triple, and the $sl(n,C)$-hierarchy is a natural hierarchy of flows on the the
function space $C(R,\cg_+)$. Moreover,  the algebraic and symplectic properties of the
$sl(n,C)$-hierarchy are determined by this Manin-triple triple. All the hierarchies we considered
in this paper are obtained as the restriction of the $sl(n,C)$-hierarchy to $C(R,(\cg_+)_\o)$,
where $(\cg_+)_\o$ is the fixed point set of certain finite order Lie algebra automorphism
$\o$ of $\cg_+$.  Symplectic structures and the group structure of B\"acklund
transformations of these hierarchies only depend on the algebra
$(\cg_+)_\o$. Finite order automorphisms of $\cg$ are not difficult to find.  For
example, if
$\cu$ is a real form of
$sl(n,C)$ defined by a conjugate linear involution $\s$, then $\s$ induces an
 involution $\hat\s$ on $\cg$:
$$\hat \s(A)(\l)= \s(A(\bar\l)).$$ The restriction of the
$sl(n,C)$-hierarchy to the subspace
$C(R,(\cg_+)_{\hat\s})$ gives the $\cu$-hierarchy.
For example, the second flow in the $su(2)$-hierarchy is the focusing non-linear
Schr\"odinger equation, and the second flow in the $u(1,1)$-hierarchy is the defocusing
non-linear Schr\"odinger equation. In general, if $\s$ is an order $k$ automorphism of
$sl(n,C)$, then it induces naturally an order $k$ automorphism $\tilde \s$ on $\cg_+$:
$$\tilde \s(A)(\l)= \s(A(\a^{-1}\l)), \qquad {\rm where\,\,\,} \a=e^{2\pi i\over k}.$$
For example, the Kupershmidt-Wilson hierarchy is of this type.  However, the automorphisms
that give the KdV and Gel'fand-Dikii hierarchies are more difficult to find. We construct these
in section 12 and 14 respectively.

This paper is organized as follows: In section 2, we review the construction of the
ZS-AKNS $sl(n)$-hierarchy of flows. In section 3, we explain various restriction defined by Lie
algebra involutions (reality conditions). We will also see in this
section that the second flow in the
$su(2)$-hierarchy is the focusing non-linear Schr\"odinger equation, and the second
flow in the $u(1,1)$-hierarchy is the defocusing non-linear Schr\"odinger equation (i.e.,
with an opposite sign of the cubic term).  The modified KdV equation is the third flow of
the $su(2)$-hierarchy twisted by the involution
$\s(x)=-x^t$, and the modified KdV with an opposite sign of the cubic term is the third
flow of the $u(1,1)$-hierarchy twisted by the involution $\tau(x)= -x^t$.  So 
this group theoretic approach allows several different generalizations of the  two
non-linear Schr\"odinger equations and the two modified KdV equations depending on
the choices of the real form of $sl(n,C)$ and of various involutions.   Although the
hierarchies associated to two different real forms may look similar  algebraically, they
have quite dissimilar global analytic behavior.   In section 4, we review the construction of local
B\"acklund transformations for the flows in the $sl(n,C)$-hierarchy and give example to show
that the new solutions obtained from applying these transformations to a smooth solution can
have singularities. We deal with only the $u(n)$-hierarchy in the next six sections.  In section 5,
we construct an  action of the rational loop group on the space of solutions of the
$j$-th flow in the $u(n)$-hierarchy and show that the action of a simple element gives a
global B\"acklund transformation. In section 6, we prove a relation among simple
elements of the rational loop group and use this relation to prove an analogue of
the Bianchi permutability formula for the $j$-th flow. In section 7, we derive
$m$-soliton formula in closed form.  In section 8, we show that the scaling
transformations of the
$j$-th flow extend the action of the rational loop group to a one-dimensional extension of the
rational loop group. Since the $n$-dimensional system associated to $u(n)$ is given
by $n$ commuting first flows, B\"acklund theory for this system works the same
way as for the first flow. This is explained in section 9.  In section 10, we
use B\"acklund transformations and the permutability formula to  construct time
periodic solutions for the $j$-th flow and the $-1$-flow.  In section 11, we
construct B\"acklund transformations for local solutions of the 
$u(k,n-k)$-hierarchies of flows. The solutions obtained this way may be
singular.  This construction is evidence that the domain of the inverse scattering
transformation for equations whose Lax pairs satisfying the reality conditions of
a non-compact group is in general quite complicated.   In section 12, 13 and 14, we
find the reality conditions,  algebraic structures and B\"acklund
transformations for the KdV,  Kupershmidt-Wilson and Gel'fand-Dikii
hierarchies respectively. 

Some parts of the research in this paper were carried out while the first author
was  visiting MPI at Bonn in the spring of 1997 and a member of
IAS in 1997-98, and while the second author is the Distinguished Visiting Professor
of IAS in 1997-98.   We would like to thank both institutes for their generous
support.  

\bs

\newsection The ZS-AKNS $n\times n$ flows.\par

Study of the hierarchy of commuting Hamiltonian flows associated to the
non-linear Schr\"odinger equation began in 1972 with a paper by Zakharov and
Shabat ([ZS1]). Ablowitz, Kaup, Newell and Segur generalized these ideas in 1974
 [AKNS] to
any $2\times 2$ system, including sine-Gordon and modified KdV. The $n\times n$ case
was treated by Zakharov and Shabat in 1979 [ZS2]. Beals and Coifman made these
construction analytically rigorous ([BC1,2]). To understand the B\"acklund
transformations, only the algebraic description is needed. However, the scattering
theory is used in the proofs.  We review the construction of these flows and give a
simple method to construct local solutions whose reduced wave functions are
analytic at infinity. 

Let  $a=\diag(a_1, \cdots, a_n)$ be  a fixed non-zero diagonal matrix in $sl(n,C)$, and
$$\eqalign{sl(n)_a&=\{y\in sl(n,C)\n [a,y]=0\}, \cr 
sl(n)_a^\perp&=\{y\in sl(n,C)\n \tr(ay)=0\},\cr}$$
denote the centralizer of $a$ and its orthogonal complement (with respect to trace) in
$sl(n,C)$ respectively. Let $\cs(R, sl(n)_a^\perp)$ denote the space of maps in the
Schwartz class. To define the $sl(n,C)$-hierarchy of flows, we need part of the scattering
theory of Beals and Coifman ([BC1]).  Let 
$$\G_a=\{\l\in C\n {\rm Re\/}(\l(a_j-a_k))=0,  \, 1\leq j<k \leq n\}.$$ 

\refclaim[ip] Theorem ([BC1]).   If $u\in \cs(R, sl(n)_a^\perp)$,  then there exists 
$$m:R\times (C\setminus \G_a)\to GL(n,C)$$ such that 
\item {(i)} $m(x,\l)$ is meromorphic for $\l\in C\setminus \G_a$, has only poles in 
$C\setminus \G_a$, 
\item {(ii)} $m(x,\l)$ has an asymptotic expansion at $\l=\infty$:
$$m(x,\l) \sim I + m_1(x)\l^{-1} + m_2(x)\l^{-2} + \cdots,$$
\item {(iii)} $E(x,\l)= m(0,\l)^{-1}e^{a\l x} m(x,\l)$ is holomorphic for $\l\in C$. 
\item {(iv)} $E^{-1}E_x = a\l + u$ and $u=[a,m_1]$,
\item {(v)} $\lim_{x\to -\infty} m(x,\l)= I$.

\refpar[pi] Remark.  A {\it wave function\/} of $u(x)$ is a solution $\psi(x,\l)$ of 
$$\psi^{-1}\psi_x= a\l + u(x).$$ This is what in general called a trivialization of the
connection
$${d\over dx} + a\l + u.$$  For $m$ as in Theorem \refip{}, $e^{a\l x}m(x,\l)$ is a wave
function. We call $m$ the {\it global reduced wave function\/}. 

Note that if $m$ satisfies condition (ii) and (iii) of Theorem \refip{} then $m$ satisfies
condition (iv).  This follows from a direct computation because
$$E^{-1}E_x= m^{-1}am \l + m^{-1}m_x = a\l + [a, m_1] + \co(\l^{-1})$$ and $E^{-1}E_x$
is holomorphic in $\l\in C$ imply that $E^{-1} E_x= a\l + [a, m_1]$. 

\ms

\refpar[iq] Definition.  For $b\in sl(n,C)$ such that $[b,a]=0$, let $Q_{b,j}$ denote the
coefficient of $\l^{-j}$ in the asymptotic expansion of $m^{-1}bm$ at $\l=\infty$:
\refeq[ir]$$m^{-1}bm \sim Q_{b,0} + Q_{b,1}\l^{-1} + Q_{b,2}\l^{-2} + \cdots .$$

\ms 

\refpar[kx] Definition.  Let $I$ be an open interval of $R$, $\co_\infty$ an open
neighborhood of $\infty$ in $S^2=C\cup \{\infty\}$, $a\in sl(n)$ and $\G_a$ defined as
above by $a$.  A smooth map
$m:I\times (\co_\infty\setminus \G_a)\to GL(n,C)$ is called a {\it local reduced
wave function\/} of
$u:I\to sl(n)_a^\perp$ if $m$ satisfies conditions (i)-(iv) of Theorem \refip{}.

\ms

We will show that $Q_{b,j}$ is an operator in  $u=[a,m_1]$,  and will write it
as
$Q_{b,j}(u)$.  First note that $Q_{b,j}$ satisfies the following recursive formula:
\refeq[ar]$$(Q_{b,j}(u))_x + [u,Q_{b,j}(u)] = [Q_{b,j+1}(u), a].$$
 To see this, note that (iii) and (iv) of Theorem \refip{} implies that $\psi(x,\l)=e^{a\l x}
m(x,\l)$ satisfies
$$\psi^{-1}\psi_x= a\l+ u.$$ So 
 $$\left[{d\over d x} + a\l + u, \psi^{-1}b\psi\right] =0.$$ But $[b,a]=0$ implies
$\psi^{-1}b\psi= m^{-1}bm$.  Hence $Q_{b,j}$ satisfies \refar{}.  

Write 
$$Q_{b,j}= T_{b,j} + P_{b,j}\,\, \in\,\, sl(n)_a + sl(n)_a^\perp.$$ Then equation
\refar{} gives
\refeq[cf]$$\eqalign{P_{b,j} &= -\ad(a)^{-1}\left(\left(P_{b,j-1}\right)_x + \pi_1([u,
Q_{b,j-1}])\right),\cr (T_{b,j})_x &= -\pi_0([u,
P_{b,j-1}]),\cr}$$ where $\pi_0$ and $\pi_1$ denote the projection of
$sl(n,C)$ onto $sl(n)_a$ and $sl(n)_a^\perp$ with respect to $sl(n,C)=sl(n)_a +
sl(n)_a^\perp$ respectively. 

The following theorem is proved by Sattinger [Sa] if $a$ has distinct eigenvalues.  Since his
proof gives an explicit method to compute the $Q_{b,j}$'s, we will repeat
it here. 

\refclaim[iz] Theorem ([Sa]). Let $a\in sl(n)$, and $m$ a local reduced wave function
for
$u:I\to sl(n)_a^\perp$. If $b$ is a polynomial of $a$, then the coefficient $Q_{b,j}$ of
$\l^{-j}$ in the asymptotic expansion of $m^{-1}bm$ is an order $(j-1)$ polynomial
differential operator in
$u$. 

\proof  Since $b=p(a)$ for some polynomial $p$ and $m^{-1}bm=p(m^{-1}am)$, it suffices
to prove the Theorem for $b=a$.  
 It is easy
to see that $Q_{a,1}=u$. We will prove $Q_{a,j}$ is a polynomial differential operator in $u$
by induction. Suppose
$Q_{a,i}$ is a polynomial differential operator in $u$ for $i\leq  j$. Write
$$Q_{a,i}=P_{a,i}+T_{a,i}\,\,\in\,\, sl(n)_a^\perp + sl(n)_a$$ as before.  Using formula
\refcf{}, we see that $P_{a,j+1}$ is a polynomial differential operator in $u$. But we can not
conclude from formula \refcf{} that $T_{a,j+1}$ is a polynomial differential operator in $u$.  
Suppose $a$ has $k$ distinct eigenvalues $c_1, \cdots, c_k$. Then 
$$f(t)=(t-c_1)(t-c_2)\cdots (t-c_k)$$ is the minimal polynomial of $a$. So
$f(m^{-1}am)=0$, which implies that the formal power series 
\refeq[bu]$$f(a+Q_{a,1}\l^{-1} + Q_{a,2} \l^{-2} +\cdots )=0.$$
Notice that $f'(a)$ is invertible and $T_{a,j+1}$ commutes with $a$.  Now
compare coefficient of $\l^{-(j+1)}$ in equation \refbu{} implies that $T_{a, j+1}$ can
written in terms of $a, Q_{a, 1}, \cdots , Q_{a,j}$. This proves that $Q_{a,j+1}$ is a
polynomial differential operator in $u$.
\qed

The following Proposition follows from formula \refcf{}. 

\refclaim[gm] Proposition.  Suppose $u(\cdot, t)\in \cs(R,sl(n)_a^\perp)$ for all $t$,
$$\left[{\p\over
\p x} + a\l + u, \,\,{\p\over
\p t} + b\l^j + v_1\l^{j-1} + \cdots + v_j\right]=0$$ for some $v_1, \cdots, v_j$,  and
$\lim_{x\to -\infty}$ $v_k(x,t) =0$ for all
$1\leq k\leq j$. Then $v_k=Q_{b,k}(u)$.

\refpar[cd] Definition. 
The $j$-th flow in the $sl(n,C)$-hierarchy  on
$\cs(R,sl(n)_a^\perp)$ defined by $b$ is the evolution equation
\refeq[ae]$$u_t=(Q_{b,j}(u))_x  + [u, Q_{b,j}(u)] = [Q_{b,j+1}(u),a].$$

\ms

\refpar[gl] Example. 
If $n=2$ and
$a=\pmatrix{1&0\cr 0&-1\cr}$, then $$sl(2)_a^\perp=\left\{\pmatrix{0&q\cr
r&0\cr}\bigg|\, \, q,r\in C\right\}.$$ For 
$u=\pmatrix{0&q\cr r& 0\cr}$, we have 
$$\eqalign{Q_{a,1}(u)&= u =\pmatrix{0&q\cr r&0\cr},\cr
Q_{a,2}(u) &= \pmatrix{-{qr\over 2}& -{q_x\over 2}\cr {r_x\over 2}& {qr\over
2}\cr  },\cr
Q_{a,3}(u) &=\pmatrix{-{1\over 4}(qr_x-rq_x)& {1\over 4}(q_{xx}- 2q^2 r)\cr 
{1\over 4}(r_{xx}- 2qr^2)& {1\over 4}(qr_x-rq_x)\cr}\cr },$$ and the first three
flows in the $sl(2,C)$-hierarchy are 
\refeq[dh]$$\eqalign{ & q_t=  q_x, \quad r_t=r_x,\cr
& q_t= -{1\over 2}(q_{xx}-2q^2 r), \quad r_t= {1\over 2}(r_{xx} - 2qr^2),\cr 
&  q_t= {1\over 4}(q_{xxx}- 6 qrq_x), \quad r_t= {1\over 4}(r_{xxx} - 6
qrr_x).\cr  }$$

\ms
 As a consequence of the recursive formula \refar{}, we have

\refclaim[ab] Proposition. Let $\co$ be an open subset of $R^2$.  Then $u:\co\to
sl(n)_a^\perp$ is a solution of the
$j$-th flow equation
\refae{} if and only if 
$$\left[{\p\over \p x} + (a\l + u),\,\, {\p\over \p t} + (b\l^j + Q_{b,1}(u)\l^{j-1} +
\cdots + Q_{b,j}(u))\right]=0.$$  
(This pair of operators is called a {\it Lax pair\/} for the $j$-th flow \refae{}).

 The following proposition is elementary. 

\refclaim[an] Proposition. 
Given two smooth maps $A, B:R^2\to gl(n)$, the following
three statements are equivalent:
\item {(i)} $[{\p\over \p x} + A,  {\p\over \p t} + B]=0$,
\item {(ii)} $B_x-A_t + [A,B]=0$,
\item {(iii)} $\cases{E_x=EA,&\cr E_t=EB, &\cr}$ is solvable.

 It follows from Propositions \refab{} and \refan{} that if $u$ is a solution of the $j$-th flow
\refae{} defined by
$b$ then there exists a unique solution $E(x,t,\l)$ for
\refeq[af]$$\cases{E_x=E \,(a\l + u), &\cr E_t= E\, (b\l^j + Q_{b,1}(u)\l^{j-1} + \cdots
+ Q_{b,j}(u)), &\cr  E(0,0,\l)=I.&\cr}$$ 
Such $E$ will be called the
{\it trivialization\/} of $u$ normalized at $(0,0)$. 
 This is a different normalization than
found in scattering theory  and much of the algebraic
literature.   Note that this choice of base point $(0,0)$ is not canonical. Different
choices of base point will not change the definition of B\"acklund transformations
which we will construct. However, the group structure of B\"acklund transformations
depends on this choice
$(0,0)$. For the rest of this paper we fixed the normalization at $(0,0)$.

B\"acklund transformations are described algebraically by starting with the
solution $u$, constructing the trivilization $E$ of $u$, operating on $E$ to produce
$\tilde E$, and producing a new solution $\tilde u$ from $\tilde E$. So it is
important to describe when $E$ is the trivialization belonging to a solution $u$.
Proposition \refgm{} gives such a condition for global solutions, and we give a
condition for local solutions next.  

\refclaim[it] Proposition. Let $\co$ be an open subset of $R^2$,
$\co_\infty$ an open subset of $S^2$ at $\infty$, 
$a\in sl(n)$, $b=p(a)$ for some polynomial $p$,  $m:\co\times(\co_\infty\setminus
\G_a) \to GL(n,C)$ a smooth map.  Suppose $m(x, t, \l)$ is meromorphic in
$\l\in C\setminus \G_a$, the asymptotic expansion of
$m$ at $\infty$ is 
$$m(x,t,\l) \sim I + m_1(x,t)\l^{-1} + m_2(x,t)\l^{-2} + \cdots , $$ and
\refeq[mi]$$E(x,t,\l)= m(0,0,\l)^{-1} e^{a\l x + b\l^j t} m(x,t,\l)$$ is holomorphic in $\l\in
C$. Then $u=[a,m_1]:\co\to sl(n)_a^\perp$ is a solution of the $j$-th flow equation
\refae{} and $E$ is the trivialization of $u$ normalized at $(0,0)$.  (We call $m$ a 
reduced wave function of the local solution $u$). 

\proof  Since $E$ is holomorphic in $\l\in C$, $E^{-1}E_x$ and $E^{-1}E_t$ are holomorphic
in $\l\in C$. But 
$$\eqalign{E^{-1}E_x&= m^{-1}am\l +m^{-1}m_x,\cr
&= (a + [a,m_1]\l^{-1})\l + O(\l^{-1}) = a\l + u + \co(\l^{-1}).\cr}$$
So $E^{-1}E_x - (a\l+ u)$ is holomorphic, bounded in $\l\in C$ and tends to zero as $\l\to
\infty$.  By Liouville Theorem,
 $$E^{-1}E_x = a\l + u.$$ 

Note that $E^{-1}E_x= \Psi^{-1}\Psi_x$ and $m^{-1}bm=\Psi^{-1}b\Psi$, where
$\Psi(x,t,\l)=e^{a\l x+ b\l^j t} m(x,t,\l)$. So
$$[d_x+a\l + u, m^{-1}bm]=[d_x+a\l + u, \Psi^{-1}b\Psi]=0.$$ Then the proof of Theorem
\refiz{} implies that 
$$m^{-1}bm \sim b+ Q_{b,1}(u) \l^{-1} +Q_{b,2}(u)\l^{-2} + \cdots.$$
A direct computation gives $$E^{-1}E_t = m^{-1}bm\l^j + m^{-1} m_t = (m^{-1}bm\l^j)_+
+ \co(\l^{-1}).$$ Since $E^{-1}E_t$ is holomorphic in $\l\in C$, we get
$$E^{-1}E_t= b\l^j +
Q_{b,1}(u)\l^{j-1}+ \cdots + Q_{b,j}(u).$$ Propositions \refan{} and \refab{} imply that $u$
satisfies the $j$-th flow equation. \qed

\ms 
The above Proposition gives a simple method to construct a class of local solutions
of the
$j$-th flow, whose local reduced wave functions are analytic at $\l=\infty$.  To
explain this, we need the Birkhoff Factorization Theorem. Recall that $L_+(GL(n,C))$
is the group of holomorphic maps from $C$ to $GL(n,C)$ 
 (under pointwise
multiplication), $L_-(GL(n,C))$ is the group of holomorphic maps $h_-$ from $\co_\infty$ to
$GL(n,C)$ such that $h_-(\infty)=I$, and $L(GL(n,C))$ is the group of holomorphic maps from
$\co_\infty\cap C$ to
$GL(n,C)$, where $\co_\infty$ is an open subset near $\infty$ in $S^2=C\cup \{\infty\}$. 

\refclaim[ma] Birkhoff Factorization Theorem.  The multiplication map
$$\eqalign{&\mu:L_+(GL(n,C))\times L_-(GL(n,C))\to L(GL(n,C)), \cr
&\quad (f_+,f_-)\mapsto
f_+f_-,\cr}$$ is a diffeomorphism onto  an open dense subset of $L(GL(n,C))$.  

Given $f_\pm\in L_\pm(GL(n,C))$, if we can factor 
$$f_-f_+= \tilde f_+\tilde f_-\,\,\in L_+(GL(n,C))\times L_-(GL(n,C)),$$  then the left
dressing action of
$f_-$ on
$f_+$ (resp. the right dressing action of $f_+$  on $f_-$) is defined by 
$$f_\sharp f_+=\tilde f_+ \, \, \quad({\rm resp.\,\,} f_-\natural f_+= \tilde f_-).$$ 
These dressing actions are only defined locally. 
The image of $\mu$ is the top Bruhat cell, the singularities in the factorization
occur on lower dimensional cells, which have codimension at least two (cf. [PS]). 

\ms
Let $a\in sl(n,C)$, $b\in sl(n)_a$ (i.e., $[a,b]=0$), $j>0$ an
integer, and
$e_{a,b,j}(x,t)$ the two-parameter subgroup in $L_+(GL(n,C))$ defined by
$$e_{a,b,j}(x,t)(\l)= e^{a\l x+ b\l^jt}.$$  

\refclaim[mb] Proposition.  If $f_-\in L_-(GL(n,C))$, then there exists an open neighborhood
$\co$ of $(0,0)$ in
$R^2$ such that 
$f_-^{-1}e_{a,b,j}(x,t)$ can be factored uniquely as 
$$f_-^{-1} e_{a,b,j}(x,t)= E(x,t)m(x,t)^{-1}\,\, \in L_+(GL(n,C))\times L_-(GL(n,C))$$
for $(x,t)\in \co$.  Moreover,
\item {(i)} $m(x,t)(\l)$ and $E(x,t)(\l)$ are smooth in $(x,t)\in \co$,
\item {(ii)} $u_{f_-}(x,t)=[a,m_1(x,t)]:\co\to sl(n)_a^\perp$ is a solution of the $j$-th flow
\refae{}, where $m_1(x,t)$ is the coefficient
of  $\l^{-1}$ of the expansion of $m(x,t)(\l)$ at $\l=\infty$ and $E$ is the trivialization of
$u_{f_-}$, 
\item {(iii)} $m$ is a local reduced wave function for $u_{f_-}$ and $m(x,t)(\l)$ is
analytic at
$\l=\infty$. 

\proof 
 Since $f_-^{-1}e_{a,b,j}(0,0)= f_-^{-1}$ lies in the image of the multiplication map
$\mu$. By Birkhoff Theorem \refma{}, the image of $\mu$ is open.  So there  exists $\co$
such that
 if $(x,t)\in \co$ then $f_-^{-1}e_{a,b,j}(x,t)$ can be factored uniquely as 
$$f_-^{-1}e_{a,b,j}(x,t)=E(x,t)m(x,t)^{-1}\,\, \in L_+(GL(n,C))\times L_-(GL(n,C)).$$ Since
the map
$\mu$ is smooth, (i) follows.   The proof of Proposition \refit{} implies (ii) and (iii). \qed

\ms
\refpar[mk] Remark. The class of local solutions constructed by
Proposition \refmb{}  contains
multi-soliton solutions and algebraic geometry solutions.  We will prove later that
if $f\in L_-(GL(n,C))$ is rational then solution $u_f$ can be given explicitly.  The
reduced wave function described using scattering theory typically have only
asymptotic expansions at $\infty$ and are meromorphic off the scattering rays
$\G_a$ described in Theorem \refip{}. The relevant factorizations will all extend to
cover this case, which is described in detail in [TU1] for $\G_a= R$ and the
$SU(n)$-reality condition.

\bs

\newsection Reality conditions.\par

To get the focusing and non-focusing non-linear Schr\"odinger
equations we need impose reality conditions on the $sl(2,C)$-hierarchy.
We explain reality conditions given by involutions of $sl(n,C)$. This group theoretic approach
allows several different generalizations of the two non-linear Schr\"odinger equations and two
modified KdV equations depending on the choices of the involutions of $sl(n,C)$. 

\refpar[me] Definition.   Let
$\cu$ denote a real form of $sl(n,C)$, i.e., $\cu$ is the fix point set of
some complex conjugate linear, Lie algebra involution
$\s$ of $sl(n,C)$. 
\item {(a)} A map
$A$ from $C$ to $sl(n,C)$ is said to satisfies the $\cu$-{\it reality condition\/} if
$$\s(A(\bar\l))=A(\l), \quad \rm {for\, all\,\,} \l\in C.$$
\item {(b)} A Lax pair $[{\p\over \p x} + A(x,t,\l), {\p\over \p t} + B(x,t,\l)]= 0$
 is said to satisfy the $\cu$-reality condition if $\s(A(x,t,\bar \l))= A(x,t,\l)$ and
$\s(B(x,t,\bar\l))= B(x,t,\l)$.

\ms

 It is clear that
$A=\sum_{k\leq n_0} u_k\l^j$ satisfies the $\cu$- reality condition if
and only if $u_k\in \cu$ for all $k$.  For example,
$A$ satisfies
\item {(i)} $su(n)$-reality condition if
$A(\bar \l)^*+A(\l)=0$ for all $\l\in C$,
\item {(ii)} $su(1,n-1)$-reality condition if $A(\bar\l)^*J+JA(\l)=0$ for all
$\l\in C$, where $J=\diag(1, -1,\cdots, -1)$,
\item {(iii)} $sl(n,R)$-reality condition if $\overline {A(\bar\l)}=A(\l)$. 

\ms

For $a\in \cu$, let 
$$\eqalign{\cu_a &=\{y\in \cu\n [a,y]=0\},\cr
\cu_a^\perp &=\{ y\in \cu\n
\tr(ay)=0\}=sl(n)_a^\perp\cap \cu.\cr}$$

\refclaim[gn] Proposition. Let $\cu$ be a real form of $sl(n,C)$, $a, b\in \cu$ such
that $[a,b]=0$, and $u\in \cs(R, \cu_a^\perp)$.  Then
\item {(1)} $Q_{b,j}(u)\in \cu$ for all $j$,
\item {(2)} the Lax pair of the $j$-th flow satisfies the $\cu$-reality condition,
\item {(3)} the $j$-th flow in the $sl(n,C)$-hierarchy leaves
$\cs(R,\cu_a^\perp)$ invariant.

\proof Let $\s$ denote the involution defines $\cu$. Set 
$$\eqalign{&A(x,t,\l)=a\l + u(x,t),\cr
&B(x,t,\l) = b\l^j + Q_{b,1}(u)\l^{j-1} +\cdots + Q_{b,j}(u).\cr}$$ 
 It follows from $a, u \in \cu$ that $A$ satisfies the $\cu$-reality condition.  Since $\s$ is a
homomorphism of $sl(n,C)$,
$$\eqalign{&\s\left(\left[{\p\over \p x} + A(x,t,\bar\l),\,\, {\p\over \p t} +
B(x,t,\bar\l)\right]\right)\cr &=\left [{\p\over \p  x} +
\s(A(x,t,\bar\l)), \,\,{\p\over \p t} + \s(B(x,t,\bar\l))\right]\cr &= \left[{\p\over \p x} + 
A(x,t,\l),
\,\,{\p\over \p t} + \s(B(x,t,\bar\l))\right].\cr}$$
  Proposition \refgm{} implies $\s(B(x,t,\bar\l))=B(x,t,\l)$, which proves (1) and (2).  
Statement (3) follow from (1).
\qed

\refpar[gp] Definition. Let $\cu$ be a real form of $sl(n,C)$. The restriction of the
$sl(n,C)$-hierarchy of flows to
$\cs(R, \cu_a^\perp)$ is called the $\cu$-hierarchy.   

\refclaim[go] Corollary. Let $\cu$ be the real form of $sl(n,C)$ defined by $\s$, 
$\hat \s$ the induced involution on $SL(n,C)$, and $U$ fixed point set of $\hat\s$.  If 
$u$ is a solution of the
$j$-th flow in the $\cu$-hierarchy, then the trivialization of $u$
satisfies the $U$-reality condition:
$\hat\s(E(x,t,\bar\l))=E(x,t,\l)$. 

\refpar[ad] Examples. 

\ni (1) The $su(2)$-hierarchy.  Note that $su(2)$ is the fixed point set of the
involution $\s(y)=-y^*$ on $sl(2,C)$.  For
$a=\pmatrix{i&0\cr 0&i\cr}$,  
$\cu_a^\perp=\left\{\pmatrix{0& q\cr -\bar q&0\cr}\biggl| \, q\in C\right\}$.
So the space $\cs(R,\cu_a^\perp)$ can be identified as $\cs(R, C)$.  The  first three flows in
the $su(2)$-hierarchy are
$$\eqalign{&q_t=q_x,\cr & q_t= {i\over 2}(q_{xx} + 2\n q\n^2 q),\cr
&q_t = -{1\over 4} (q_{xxx} + 6\n q\n^2 q_x).\cr }$$ 
Note that the first flow just gives translation, 
the second flow is the focusing non-linear Schr\"odinger equation, and the
sequence of flows is the hierarchy of commuting flows associated to the 
non-linear Schr\"odinger equation. 
\ms

\ni (2)  The $su(n)$-hierarchy.  If $a=\diag(a_1, \cdots, a_n)\in su(n)$ has
distinct eigenvalues and
$b=\diag(b_1,\cdots, b_n)\in su(n)$, then 
$$\cu_a^\perp=\{ (u_{ij})\in su(n)\n u_{ii}=0 \,\, {\rm for\, all\,\,} 1\leq
i\leq n\}.$$
The first flow in the $su(n)$-hierarchy on $\cs(R,\cu_a^\perp)$ defined by $a$ is  the
translation
$$u_t=u_x.$$ The first flow defined by $b$ ($a, b$ linearly
independent) is the
$n$-{\it wave equation \/} ([ZMa1, 2]):
\refeq[am]$$(u_{ij})_t={b_i-b_j\over a_i-a_j} (u_{ij})_x + \sum_{k\not= i, j}
\left({b_k-b_j\over a_k-a_j} - {b_i-b_k\over a_i-a_k}\right) u_{ik} u_{kj}, \qquad
i\not=j.$$

\ms

\ni (3) The $u(n)$-hierarchy.  Let $a=\diag(i,\cdots, i, -i, \cdots, -i)$
be the diagonal matrix with eigenvalues $i, -i$ and multiplicities $k, n-k$ respectively.   Then 
$$\cs(R,\cu_a^\perp)=\left\{u=\pmatrix{0&q\cr -q^*&0\cr}\bigg|  \, q\in
\cs(R,\cm_{k\times (n-k)})\right\},$$ where $\cm_{k\times (n-k)}$ is the space of
$k\times (n-k)$ complex matrices and $q^\ast= \bar q^t$. So
$\cs(R,\cu_a^\perp)$ is naturally identified as
$\cs(R,\cm_{k\times (n-k)})$.    For 
$u=\pmatrix{0&q\cr -q^*&0\cr}$, we have 
$$\eqalign{Q_{a,0}(u)&=a, \cr Q_{a,1}(u) &= u,\cr Q_{a,2}(u)&=\pmatrix{{1\over 2i}
qq^*& {i\over 2} q_x\cr {i\over 2} q^*_x& -{1\over 2i} q^*q\cr}.\cr}$$ So  the first
three flows on $\cs(R,\cm_{k\times(n-k)})$ in the $su(n)$-hierarchy defined by $a$ are 
$$\eqalign{q_t &= q_x\cr
q_t &= {i\over 2} (q_{xx} + 2 qq^\ast q)\cr
q_t &= -{1\over 4} q_{xxx} - {3\over 4} (q_xq^\ast q + qq^\ast q_x).\cr}$$ 
Note that the second flow is the matrix non-linear Schr\"odinger equation studied
by Fordy and Kulish [FK].  

\ms

Next we recall the definition of $u(k, n-k)$.  Let 
$J=\diag(\e_1, \cdots, \e_n)$ with $\e_i= 1$ for $1\leq i\leq k$ and $\e_j= -1$ if
$k<j\leq n$, and let $$\li v_1, v_2\ri_J = v_1^*Jv_2$$ denote the Hermitian bilinear form
on
$C^n$ defined by $J$. Let $U(k,n-k)$ denote the group of linear maps of $C^n$ that
preserve $\li\,  , \,\ri_J$, and $u(k,n-k)$ its Lie algebra.  Then 
$$\eqalign{U(k,n)&=\{g\in GL(n,C)\n g^*Jg=J\},\cr
u(k,n) &=\{X\in gl(n,C)\n X^*J + JX = 0\}. \cr}$$ 
The involution that defines $u(k,n-k)$ is $\s(y)= -J^{-1}y^*J$, and the induced involution
on $U(k,n-k)$ is $\hat \s(g)=J^{-1}(g^*)^{-1}J$. 

\ss
\refpar[mt] Example.  The $u(1,1)$-hierarchy. Here 
$$u(1,1)=\left\{y\in sl(2,C)\n y^*J+Jy= 0\right\}=\left\{\pmatrix{ir& q\cr \bar q&
-ir\cr}\biggl|\,\, r\in R, q\in C\right\},$$ where $J=\diag(1,-1)$. Let $a=\diag(i,-i)\in u(1,1)$.
Then
 $$\cu_a^\perp=\{y\in sl(2,C)\n y^*J+ Jy =
0.\}=\left\{\pmatrix{0&q\cr
\bar q &0\cr}\bigg| \,q\in C\right\}.$$  The second flow in the $u(1,1)$-hierarchy is  the {\it
defocusing non-linear Schr\"odinger equation\/}:
\refeq[dk]$$q_t={i\over 2} (q_{xx}- 2\n q\n^2 q).$$ 

\ms

The classical B\"acklund transformation is a transformation of solutions,
not of the $j$-th flow, $j\geq 1$, but of the $-1$ flow. As an
evolution, the
$-1$ flow is poorly defined. The $(x,t)$ are characteristic coordinates.
However, the B\"acklund transformation operates algebraically on
solutions.  For the physical problem, we would not expect $u$ to
necessarily be in the Schwartz space along the characteristic coordinate
$x$ for all solutions. However, for the solutions we  construct in this
paper, this is the case.  

The  $-1$ flow in the $su(n)$-hierarchy defined by $b\in su(n)_a$ ([Te], [TU1])  is 
\refeq[bt]$$\cases{u_t=[a,g^{-1}bg],&\cr g^{-1}g_x= u, \quad \lim_{x\to -\infty}
g(x,t)= I.&\cr}$$ Its Lax pair is 
$$\left[ {\p\over \p x} + a\l + u,\quad \l^{-1} g^{-1}bg\right]=0.$$
Everything we said about the $j$-th flows, $j\geq 1$, applies to the $-1$
flow, except that the trivialization $E(x,t,\l)$ has a singularity at $0\in C$,
or $E(x,t,\l)$ is holomorphic in $\l\in C-\{ 0\}$.

We require a further step to connect with the classical theory, since
sine-Gordon is the $-1$ flow in the $su(2)$-hierarchy restricted (or twisted) by an
involution (this is related to the twisted affine Kac-Moody algebras). This is also
referred to in the literature as reduction, which is a unfortunate terminology as
reduction has a specific meaning in symplectic geometry.   Assume
$\tau$ is a complex conjugate linear, Lie algebra involution of $sl(n,C)$ and 
$\s$ is a complex linear, Lie algebra involution of $sl(n,C)$ such that $\s\tau=\tau\s$. 
Then $(\tau, \s)$ defines a symmetric space as follows:    Let
$\cu$ denote the real form defined by
$\tau$. Then
$\s(\cu)\subset
\cu$.    Let
$\ck,\cp$ denote the
$1, -1$ eigenspaces of $\s$ on $\cu$ respectively.  Then $\cu = \ck + \cp$, and 
$$[\ck, \ck]\subset \ck,  \quad [\ck,\cp]\subset \cp, \quad [\cp, \cp]\subset \ck.$$
Let $U$ and $K$ denote the subgroup corresponding to $\cu$ and $\ck$ respectively. Then
$U/K$ is a symmetric space, and $\cu=\ck+\cp$ is the {\it Cartan
decomposition\/} of $U/K$.  Fix $a\in \cp$. Let   
$$\eqalign{\ck_a &=\{x\in \ck\n [x,a]=0\}= \ck \cap \cu_a,\cr
 \cu_{a,\s}^\perp &=\ck\cap \cu_a^\perp \cr}.$$

\refpar[iv] Definition.  Let $\cu$ be the real form of $sl(n,C)$ defined by the complex
conjugate linear, Lie algebra involution $\tau$ of $sl(n,C)$,  $\s$ a complex linear, Lie
algebra involution of $sl(n,C)$ such that $\tau\s=\s\tau$, and $U/K$ the corresponding
symmetric space.   We say that
$A(\l)$ satisfies the
$\cu$-reality condition twisted by $\s$ or the $U/K$-reality condition if 
\refeq[iu]$$\tau(A(\bar\l))=A(\l), \quad \s(A(-\l))= A(\l).$$

\ms
A direct computation shows that  $A(\l)=\sum_j v_j\l^j$ satisfies the $\cu$-reality
condition twisted by $\s$ if 
$v_j\in \ck$ if $j$ is even, and $v_j\in \cp$ if $j$ is odd. 

If $u\in \cu_{a,\s}^\perp$, then
the recursive formula \refar{} implies that $Q_{b,j}(u)\in \ck$ if $j$ is odd and is in $\cp$ if
$j$ is even. Since $a\in \cp$, $[Q_{b,j}(u),a]$ is in $\cp$ if $j$ is even and is in $\ck$ if $j$
is odd. This proves

\refclaim[cu] Theorem ([Te]). Let $\tau, \s, \cu, \ck, \cp$ be as above, and 
$a,b\in \cp$ such that $[a,b]=0$.  Let 
$\cu_{a,\s}^\perp=\ck\cap \cu_a^\perp$.   If $u\in \cs(R,\cu_{a,\s}^\perp)$, then
\item {(i)} $Q_{b,j}(u)\in \ck$ if $j$ is odd and is in $\cp$ if
$j$ is even,
\item {(ii)} $[Q_{b,j}(u),a]\in \cp$ if $j$ is even, and is in $\ck$ if $j$ is odd,
\item {(iii)}  $\cs(R,\cu^\perp_{a,\s})$ is invariant under the odd
flows, and is stationary under the even flows. 

\refpar[cv] Definition. The $j$-th  ($j$ odd)  flow in the $\cu$-hierarchy {\it twisted by\/}
$\s$ defined by $a, b$  (or the hierarchy associated to the symmetric space $U/K$) is the
$j$-th flow in the $\cu$-hierarchy restricted to $\cs(R, \cu_{a,\s}^\perp)$:  
\refeq[cw]$$ u_t=(Q_{b,j}(u))_x + [u, Q_{b,j}(u)] = [Q_{b,j+1}(u),a], \quad  u:R^2\to
\cu_{a,\s}^\perp.$$

It follows easily from Theorem \refcu{} that we have

\refclaim[gz] Corollary. If
$u:R^2\to \cu_{a,\s}^\perp$ is a solution of the $j$-th flow ($j$ odd) in the $\cu$-hierarchy,
then
$$\eqalign{A(x,t,\l) &=a\l+u(x,t),\cr
B(x,t,\l) &=b\l^j + Q_{b,1}(u)\l^{j-1} + \cdots + Q_{b,j}(u)\cr}$$ satisfy the
$\cu$-reality condition and the $\s$-reality condition:
$$\s(A(x,t,-\l))= A(x,t,\l), \quad \s(B(x,t,-\l))= B(x,t,\l).$$
In particular, the trivialization $E$ of $u$ satisfies the following reality conditions
$$\hat\tau (E(x,t,\bar\l))=E(x,t,\l), \quad \hat \s(E(x,t,-\l)) = E(x,t,\l).$$

\refpar[by] Example. The hierarchy associated to $SU(n)/SO(n)$.   Let $\tau(y)=-y^*$,
$\s= -y^t$, and
$a=\diag(i,-i,
\cdots, -i)$. Then $\cu=su(n)$, and
$$\cs(R,\cu_{a,\s}^\perp)=\left\{\pmatrix{0&v\cr
-v^t&0}\biggl| \,\, v\in \cs(R,\cm_{1\times (n-1)})\right\},$$ where 
$\cm_{1\times (n-1)}$ is the space of real $1\times (n-1)$ matrices.   The even flows
vanishes on
$\cs(R,\cu_{a,\s}^\perp)$, and the odd flows are extensions of the usual hierarchy of
flows for the modified KdV. The third flow twisted by $\s$, written in terms of
$v:R\to
\cm_{1\times (n-1)}$,  is the {\it matrix modified KdV equation\/}:
\refeq[js]$$v_t=-{1\over 4}\left(v_{xxx} + 3(v_xv^tv+vv^tv_x)\right).$$
(When $n=2$, $v=iq$ is a scalar function and the above equation is the classic
modified KdV equation: $q_t=-{1\over 4}(q_{xxx} - 6q^2q_x)$.)  \ss

\refpar[jo] Example. The  third flow in the
$sl(n,R)$-hierarchy twisted by $\s(y)=-y^t$ defined by $a=\diag(-1, 1,\cdots, 1)$ (the
$SL(n, R)/SO(n)$ hierarchy) is the equation \refjs{}. When $n=2$, $SL(2,R)/SO(2)$ is 
the hyperbolic $2$-plane $H^2$.   The third flow in
the $H^2$-hierarchy is the other modified KdV:
\refeq[dr]$$q_t= {1\over 4}(q_{xxx} + 6q^2q_x).$$ 
(Here $u=\pmatrix{0&q\cr -q&0\cr}$ with $q\in R$.)

\ss

\refpar[bn] Example. The $-1$-flow associated to $SU(2)/SO(2)=S^2$.  Let $\tau$ and $\s$
be the involution in Example
\refby{}, and
$a=\diag(i,-i)$. Then  
$$\cu_{a,\s}^\perp=\left\{\pmatrix{0&v\cr -v&0\cr}\biggl| \,\,
v\in R\right\}.$$  The
$-1$ flow defined by
$b=-a/4$ twisted by $\s$ is the equation \refbt{} for $u=\pmatrix{0&{q_x\over 2}\cr
-{q_x\over 2}&0\cr}$ is
 the Sine-Gordon equation
\refeq[dt]$$q_{xt}=\sin q.$$

\refpar[ju] Example. The $u(1,1)$-hierarchy twisted by the involution
$\tau(y)=-y^t$ is the hierarchy associated to the Lorentzian symmetric space
$$U(1,1)/(U(1,1)\cap O(2,C)).$$ The third flow in this hierarchy is the modified KdV with an
opposite sign in the cubic term:
\refeq[ds]$$q_t= -{1\over 4}(q_{xxx} - 6q^2q_x).$$  
(Here $u=\pmatrix{0&iq\cr -iq&0\cr}$  for some real function $q$.)

\ss
   
\refpar[ky] Example.  The $u(1,n-1)$-hierarchy twisted by $\s(y)=-y^t$ is the hierarchy
associated to the Lorentzian symmetric space $$U(1, n-1)/(U(1,n-1)\cap O(n,C)).$$
The third flow in this hierarchy defined by $a=\diag(i,-i,\cdots, -i)$ is a generalization of
the other modified KdV \refds{}: 
$$v_t= -{1\over 4}(v_{xxx} -3(v_xv^tv + vv^t v_x)), \quad v:R^2\to \cm_{1\times
(n-1)}.$$ 

\bs

\newsection B\"acklund transformations for the $sl(n,C)$-hierarchy.\par

The scattering data for the integrable systems we are considering 
have two parts, discrete data in $C\setminus \G_a$ and continuous data along
$\G$. The group which generates the discrete data is intimately connected
with B\"acklund transformations.
 Local Darboux and B\"acklund transformations for the $j$-th flow in the $sl(n,C)$-hierarchy
were constructed by many authors (Zakharov and Shabat [ZS 2], Sattinger and Zurkowski 
[SZ 1, 2], Gu and Zhou [GZ] and Cherdnik [Ch]):

\refclaim[hmb] Theorem.  Suppose $u$ is a solution of the $j$-th flow in the
$sl(n,C)$-hierarchy that admits a local reduced wave function.  Let $E$ denote the
trivialization of $u$. Let $\a_1, \a_2\in C$, and
$V_1, V_2$ complex linear subspace of $C^n$ such that 
$C^n=V_1\oplus V_2$. Set $\tilde V_i(x,t) =  E(x,t,\a_i)^{-1}(V_i)$. Suppose $\tilde
V_1(x,t)\cap \tilde V_2(x,t)= 0$ for
$(x,t)$ in an open subset $\co$.   Then
 $\tilde u= u +(\a_1-\a_2)[a, \tilde \pi]$ is a solution of the $j$-th flow on
$\co$, where $\tilde \pi(x,t)$ is the projection onto
$\tilde V_1(x,t)$ with respect to $C^n= \tilde V_1(x,t)\oplus \tilde V_2(x,t)$. 

We will reformulate this theorem in terms of the dressing
action and give a proof, which will be used in the later sections.  
Let $L_\pm(GL(n,C))$ be as in section 2. 
By the Birkhoff Factorization Theorem \refma{}, there is a local dressing action
$\sharp$ of
$L_-(GL(n,C))$ on $L_+(GL(n,C))$: Given $g\in L_-(GL(n,C))$ and $f\in
L_+(GL(n,C))$, if we can factor 
$$gf=\tilde f\tilde g\in L_+(GL(n,C))\times L_-(GL(n,C)),$$ then $g\sharp f$ is defined to
be $\tilde f$. This action $\sharp$ is only defined
for $f$ in an open dense subset of $L_+(GL(n,C))$. In certain cases, this factorization can be
constructed explicitly.   First we choose simple elements (linear fractional transformations) in
$L_-(GL(n,C))$.  Given constants
$\a_1,\a_2\in C$ and a linear projection $\pi$ of $C^n$ (i.e., $\pi$ is complex linear
and $\pi^2=\pi$), let $\pi'=I-\pi$ and 
$$h_{\a_1,\a_2, \pi}(\l) = {\l-(\a_1\pi +
\a_2\pi')\over \l- \a_1}= I + {(\a_1-\a_2)\over \l-\a_1}\pi'.$$ Then $h_{\a_1,\a_2,\pi}\in
L_-(GL(n,C))$ and 
$$h_{\a_1,\a_2, \pi}(\l)^{-1}={\l-(\a_2\pi + \a_1\pi')\over \l- \a_2}.$$
We will call $h_{\a_1, \a_2, \pi}$ a {\it simple element\/}.

\refclaim[gu] Proposition. Let $\a_1,\a_2\in C$, $\pi$ a projection of $C^n$,  
$V_1$ and $V_2$ denote the image of $\pi$ and $\pi'=I-\pi$ respectively, and
$f\in  L_+(GL(n,C))$.  If 
\refeq[ig]$$(f(\a_1)^{-1}(V_1)) \cap (f(\a_2)^{-1}(V_2))= 0,$$
then $h_{\a_1,\a_2,\pi}f$ can be factored uniquely as
$$h_{\a_1,\a_2,\pi} f = \tilde f h_{\a_1,\a_2.\tilde \pi}\,\,\in L_+(GL(n,C))\times
L_-(GL(n,C)),$$ where 
$\tilde \pi$ is the projection onto $f(\a_1)^{-1}(V_1)$ with respect to
$$C^n=f(\a_1)^{-1}(V_1)\oplus
f(\a_2)^{-1}(V_2).$$

 \proof   It suffices to prove that
$$\tilde f  = h_{\a_1,\a_2,\pi}f
h_{\a_1,\a_2,\tilde\pi}^{-1}$$  lies in $L_+(GL(n,C))$.  Since the
right hand side of $\tilde f$ is holomorphic in $C\setminus\{\a_1,
\a_2\}$ and has simple poles at $\a_1$ and
$\a_2$, we only need to prove that the residues of
$\tilde f$ are zero at both $\a_1$ and $\a_2$.   But
 $$\eqalign{{\rm Res\,}(\tilde f, \a_1) &= (\a_1-\a_2) (I-\pi) f(\a_1)\tilde
\pi,\cr {\rm Res\,}(\tilde f,
\a_2) &=(\a_2-\a_1) \pi f(\a_2)(I-\tilde \pi).\cr}$$ Since Im$(\tilde \pi) =
f(\a_1)^{-1}(V_1)$ and Im$(I-\tilde \pi)= f(\a_2)^{-1}(V_2)$, both residues are zero. 
\qed

\ms

Local Darboux and B\"acklund transformations of the $j$-th flow in the
$sl(n,C)$-hierarchy are obtained by factoring the product of a
simple element and the trivialization of a solution of the
$j$-th flow in the $sl(n,C)$-hierarchy at each $(x,t)$:

\refclaim[hma] Theorem. Let $u:\co_1\to sl(n)_a^\perp$ be a local solution of the $j$-th
flow \refae{}, $E$ the trivialization of $u$,  $h_{\a_1,\a_2,
\pi}$ a simple element in $L_-(GL(n))$, and $V_1, V_2$ denote the image of $\pi$ and
$I-\pi$ respectively. Assume $m:\co\times (C\setminus \G)\to
GL(n,C)$ is a local reduced wave function for $u$. Then there exists an open
subset $\co \subset \co_1$ such that $\tilde V_1(x,t)\cap \tilde V_2(x,t)= 0$ for
$(x,t)\in\co$, where 
$\tilde V_i(x,t) =  E(x,t,\a_i)^{-1}(V_i)$. 
Moreover, let 
$\tilde \pi(x,t)$ denote the projection onto
$\tilde V_1(x,t)$ with respect to $C^n= \tilde V_1(x,t)\oplus \tilde V_2(x,t)$, then 
\item {(i)} $\tilde u:\co\to sl(n)_a^\perp$ defined by $\tilde u= u +(\a_1-\a_2)[a, \tilde
\pi]$ is a solution of the
$j$-th flow with $\tilde E$ is the trivialization and $\tilde m$ as a local reduced
wave function, where 
$$\eqalign{\tilde E(x,t,\l)&= h_{\a_1,\a_2,\pi} E(x,t) h_{\a_1,\a_2,\tilde \pi(x,t)}^{-1},\cr
\tilde m(x,t,\l) &=m(x,t,\l)h_{\a_1,\a_2,\tilde \pi(x,t)}(\l)^{-1},\cr}$$
\item {(ii)} $\tilde \pi$ is the solution of 
\refeq[hq]$$\cases{\tilde \pi_x = -[u,\tilde \pi] -[a,\tilde \pi]
(\a_2 + (\a_1-\a_2)\tilde\pi), &\cr \tilde \pi_t = -\sum_{k=0}^j
\left[Q_{b,j-k}(u),\tilde \pi \right](\a_2 + (\a_1-\a_2)\tilde\pi)^k,&\cr
\tilde \pi^2=\tilde \pi, \quad \tilde \pi(0,0)=\pi.&\cr}$$

\proof Since $\tilde V_1(0,0)\cap \tilde V_2(0,0)=V_1\cap V_2=0$ and two linear
subspaces are in general position is an open condition, there exists an open subset $\co$ of
$(0,0)$ in $\co_1$ such that $\tilde V_1(x,t)\cap \tilde V_2(x,t)=0$ for all $(x,t)\in\co$. 

To prove (i), we first note that
 $E(x,t)\in L_+(GL(n,C))$, where $E(x,t)(\l)=E(x,t,\l)$.  By Proposition
\refgu{}, $E$ is holomorphic in $\l\in C$. This proves that $\tilde m$ is a local
reduced wave function for $\tilde u$. 

 A simple computation shows that the coefficient of
$\l^{-1}$ in the asymptotic expansion of $\tilde m$ at $\l=\infty$ is 
$$\tilde m_1(x,t)= m_1(x,t) + (\a_2-\a_1)\tilde \pi'(x,t).$$  Then (i) follows from
Proposition \refit{}. 

To prove (ii), we note that 
$$\eqalign{a\l + \tilde u &= \tilde E^{-1}\tilde E_x \cr
&=(h_{\a_1, \a_2,\pi}Eh_{\a_1,\a_2,\tilde \pi}^{-1})^{-1}(h_{\a_1,
\a_2,\pi}Eh_{\a_1,\a_2,\tilde \pi}^{-1})_x\cr 
&=h_{\a_1,\a_2,\tilde \pi}(E^{-1}E_x)h_{\a_1,\a_2,\tilde \pi}^{-1} - (h_{\a_1,\a_2,\tilde
\pi})_xh_{\a_1,\a_2,\tilde \pi}^{-1}.\cr}$$ 
Multiply the above equation by $(\l-\a_1)h_{\a_1,\a_2,\tilde \pi}$ on the right to get
$$(a\l+\tilde u)(\l-(\a_1\tilde \pi + \a_2 \tilde \pi')) =(\l-(\a_1\tilde \pi + \a_2 \tilde \pi'))
(a\l + u) + (\a_1-\a_2)\tilde\pi_x.$$
Compare coefficient of $\l^i$ for $i=0, 1$ to get the ODE for $\tilde \pi$ in $x$ variable.
Similarly, we have 
$$\eqalign{&\left(\sum_{k=0}^jQ_{b,k}(\tilde u)\l^{j-k}\right)(\l-(\a_1\tilde \pi + \a_2
\tilde
\pi')) \cr &\,\,=(\l-(\a_1\tilde \pi + \a_2 \tilde \pi'))\left(\sum_{k=0}^jQ_{b,k}(
u)\l^{j-k}\right) + (\a_1-\a_2)\tilde\pi_t.\cr}$$
Compare coefficient of $\l^i$ for $0\leq i\leq j$ to get
$$\eqalign{&Q_{b,k}(\tilde u) = Q_{b,k}(u) + \sum_{i=1}^{k}[Q_{b,k-i}(u),
(\a_1-\a_2)\tilde
\pi](\a_1\tilde \pi +\a_2\tilde \pi')^{i-1},\cr & (\a_1-\a_2)\tilde
\pi_t=(\a_1\tilde\pi +\a_2\tilde \pi')Q_{b,j}(u) - Q_{b,j}(\tilde u)(\a_1\tilde \pi+\a_2\tilde
\pi').\cr}$$  Substitute the first equation for $k=j$ to the second equation to get the ODE for
$\tilde \pi$ in $t$ variable. \qed

\ms

\refpar[mu] Remark. The assumption the $u$ admits a local reduced wave function in
Theorem \refhma{} is necessary.  Without this assumption we can only
conclude that 
$$\cases{\tilde E^{-1}\tilde E_x = a\l + \tilde u,&\cr 
\tilde E^{-1}\tilde E_t= b\l^j + v_1 \l^{j-1} + \cdots + v_j,&\cr}$$
for some $\tilde u$ and $v_1, \cdots, v_j$. In general, it is not clear whether
the $v_i$ is equal to $Q_{b,i}(\tilde u)$. But if $\lim_{x\to
-\infty}v_i(x,t)=0$ for all $t$ then Proposition \refgm{} implies that $v_i=
Q_{b,i}(\tilde u)$ for all $1\leq i\leq j$. 

\ms

\refpar[mf] Definition.  Let $u$ be a local solution of the $j$-th flow \refae{} that admits
a local reduced wave function, and $h_{\a_1,\a_2,\pi}$ a simple element.  Define
$h_{\a_1,\a_2,\pi}\ast u=
\tilde u$, where $\tilde u$ is the new local solution obtained in Theorem \refhma{}. 
\ms

The following example explains why the above constructions only provide local solutions
from a global solution of the $j$-th flow in the $sl(n,C)$-hierarchy on
$\cs(R,sl(n)_a^\perp)$. 

\refpar[dj] Example. We apply Theorem \refhma{} to the vacuum solution
$u=0$ of the $j$-th flow \refdh{} in the $sl(2,C)$-hierarchy defined by $a=\diag(1,-1)$. 
Let $\a_1, \a_2\in C$, and $V_1$ (resp. $V_2$) the  subspace spanned by $v_1=(c_1,
c_2)^t$ (resp.
$v_2=(d_1, d_2)^t$). A direct computation implies that $\tilde u$ constructed in
Theorem \refhma{} is
\refeq[gd]$$\tilde u(x,t)= {2(\a_1-\a_2)\over c_2d_1e^{\eta(x,t)}-c_1d_2
e^{-\eta(x,t)}}\pmatrix{0&c_1d_1e^{-\xi(x,t)}\cr c_2d_2e^{\xi(x,t)}
&0\cr },$$  where  $\xi=(\a_1+\a_2)x+
(\a_1^j+\a_2^j)t$ and
$\eta=(\a_1-\a_2)x+(\a_1^j-\a_2^j)t$. 
Note that $\tilde u$ is not defined at $(x_0, t_0)$ when
$e^{2((\a_1-\a_2)x_0+(\a_1^j-\a_2^j)t_0)}= {c_1d_2\over  c_2d_1}$.  For example,
for $j=2$, 
\item {(i)} if $\a_1= 2$, $\a_2=1$, and $v_1=(1,1)^t, v_2=(1,2)^t$, then $\tilde u(x,t)$ is
singular along the line $x+3t = {1\over 2}\ln 2$,
\item {(ii)} if $\a_1= 2$, $\a_2=1$, and $v_1=(1,1)^t, v_2=(-1, 2)^t$, then $\tilde u$ is
smooth on $R^2$ but $\tilde u(x,t)$ goes to infinity when $x\to\pm \infty$,
\item {(iii)} if $\a_1, \a_2\in R$ satisfying  $\n \a_1+\a_2\n < \n\a_1-\a_2\n$ and
$v_1=(c_1,c_2)^t, v_2=(d_1,d_2)^t$  in $R^2$ satisfying $c_1c_2d_1d_2<0$, 
 then $\tilde u(x,t)$ is smooth on $R^2$ and is rapidly decay in $x$. 

\ms 

The following theorem explains B\"acklund transformations for the restricted flows. 

\refclaim[jk] Theorem.  Let $\cu$ be a real form of $sl(n,C)$ defined by $\tau$, and $U/K$
the symmetric space defined by $\s$. If $u$ is a solution of the $j$-th flow in the
$\cu$-hierarchy (resp. $U/K$-hierarchy) and $h_{\a_1, \a_2,\pi}$ is a simple element
satisfying the $U$-reality condition (resp. the $U/K$-reality condition), then $h_{\a_1, \a_2,
\pi}\ast u$ is again a solution of the $j$-th flow in the $\cu$- (resp. $U/K$-) hierarchy. 

\proof Let $E$ be the trivialization of $u$ normalized at $(0,0)$. Then $E$ satisfies the
$U$-reality condition.  Let $h=h_{\a_1,\a_2, \pi}$, and $hE= \tilde E\tilde h$ as
in Theorem \refhma{}.  If $h$ satisfies the $U$-reality condition, then so is $hE$.  But
$$h(\l)E(\l)= \tau(h(\l)E(\l))=\tau(\tilde E(\bar\l)\tilde h(\bar\l))= \tau(\tilde E(\bar\l))
\tau(h(\bar\l))=\tilde E(\l)
\tilde h(\l).$$ Uniqueness of the Birkhoff decomposition implies that $\tau(\tilde
E(\bar\l))=\tilde E(\l)$.  Hence $\tilde u$ is a solution of the $j$-th flow in the
$\cu$-hierarchy.  The same proof works for the $U/K$-hierarchy.
\qed

Next we give some relations among simple elements in $L_-(GL(n,C))$. 

\refclaim[gt] Proposition. Let $h_{\a_1,\a_2,\pi_1}$ and
$h_{\b_1,\b_2, \pi_2}$  be two simple elements.  If 
$$\phi= (\a_2+(\a_1-\a_2)\pi_1) -
(\b_2 + (\b_1-\b_2) \pi_2)$$ is invertible, then
\refeq[gq]$$h_{\a_1,\a_2,\tau_1}h_{\b_1,\b_2,\pi_2} = h_{\b_1,\b_2,\tau_2}
h_{\a_1,\a_2,\pi_1}$$ if and only if $\tau_i= \phi \pi_i \phi^{-1}$ for $i=1,2$.

\proof Set $$\eqalign{&Y_1=\a_2+(\a_1-\a_2)\tau_1, \quad Y_2= \b_2+
(\b_1-\b_2)\pi_2, \cr
&Z_1=\a_2+(\a_1-\a_2)\pi_1, \quad Z_2= \b_2+(\b_1-\b_2)\tau_2.\cr}$$
Then equality \refgq{} gives 
$$(\l-Y_1)(\l-Y_2) = (\l-Z_2)(\l-Z_1).$$  
This  holds if and only if 
$$\cases{Y_1-Z_2= Z_1-Y_2, &\cr Y_1 Y_2 - Z_2 Z_1= 0. &\cr}$$ Multiply the
first equation by $Y_2$ on the right and subtract the second equation to get
$$Z_2=(Z_1-Y_2)Y_2(Z_1-Y_2)^{-1},$$ where $Z_1-Y_2=\phi$ is invertible by assumption. 
Multiply the first equation by
$Z_1$ on the right and subtract the second equation to get
$$Y_1=(Z_1-Y_2)Z_1(Z_1-Y_2)^{-1}.$$ This finishes the proof.  
\qed

As a consequence of Proposition \refgt{}, we obtain an analogue of the Bianchi permutability
formula:

\refclaim[ix] Corollary. Let $h_{\a_1,\a_2,\pi_1}$ and
$h_{\b_1,\b_2, \pi_2}$  be two simple elements such that $\phi= (\a_2+(\a_1-\a_2)\pi_1)
- (\b_2 + (\b_1-\b_2) \pi_2)$ is invertible. Let $u$ be a local solution of the $j$-th flow
\refae{}, which admits a reduced wave function.  Let
$$\eqalign{ u_1&=h_{\a_1,\a_2,\pi_1}\ast u = u + (\a_1-\a_2)[a,\tilde \pi_1],\cr
 u_2 &= h_{\b_1,\b_2,\pi_2}\ast u = u + (\b_1-\b_2)[a,\tilde \pi_2],\cr}$$  as  in
Theorem \refhma{}.  Set   $\tilde \tau_i= \tilde\phi \tilde\pi_i \tilde\phi^{-1}$, where 
$$\tilde\phi= (\a_2+(\a_1-\a_2)\tilde\pi_1)
- (\b_2 + (\b_1-\b_2)\tilde \pi_2).$$
Then $$\eqalign{u_3&= h_{\a_1,\a_2, \tau_1}\ast u_2= u_2 + (\a_1-\a_2)[a,\tilde
\tau_1],\cr 
&=h_{\b_1,\b_2,\tau_2}\ast u_1 = u_1 + (\b_1-\b_2)[a,\tilde \tau_2].\cr}$$

\bs

\newsection B\"acklund transformations for the $u(n)$-hierarchy.\par

 In this section, we consider
B\"acklund transformations for the flows in the $u(n)$-hierarchy, in which the Birkhoff
factorization can always be carried out.  Let $G_-^m$ denote the subgroup of 
rational maps $g\in L_-(GL(n,C))$ such that $g$ satisfies the $U(n)$-reality condition.  We
obtain an action of $G_-^m$ on the space of solutions of the $j$-th flow in the
$u(n)$-hierarchy on
$\cs(R,\cu_a^\perp)$. Simple elements generate $G_-^m$ and the action of these simple
elements gives global B\"acklund transformations.  

Let $z\in C$, and  $\pi$
the Hermitian projection of $C^n$ onto  a complex linear subspace $V$, i.e., $\pi^*=\pi$ and
$\pi^2=\pi$. Let $\pi^\perp= I-\pi$ be the Hermitian projection of $C^n$ onto the
orthogonal complement $V^\perp$. Let
$$g_{z,\pi}(\l)=h_{\bar z,z,\pi}(\l)= \pi + {\l- z\over \l- \bar z} \pi^\perp. $$ It is easy to
check that $g_{z,\pi}(\bar\l)^*g_{z,\pi}(\l)=I$. So
$g_{z,\pi}\in G_-^m$, and will be called a {\it simple element\/} of $G_-^m$. 
The following theorem was proved by the
second author [U1]

\refclaim[ay] Theorem ([U1]). The set $\{g_{z,\pi}\n z\in C\setminus R,$ $\pi$ a Hermitian
projection of $C^n\}$ generates $G_-^m$.  

Let $G_\pm$ denote the subgroup of $f\in L_\pm(GL(n,C))$ such that $f$ satisfies the
$U(n)$-reality condition, i.e, $f(\bar\l)^*f(\l)=I$.  If for $f\pm\in G_\pm$ we can factor
$$f_-f_+ = g_+g_- \in L_+(GL(n,C))\times L_-(GL(n,C)),$$
then the uniqueness of the Birkhoff Theorem implies that $g_\pm\in G_\pm$. This implies
that the dressing action of $G_-$ leaves $G_+$ invariant.  Let
$\hat G_+$ (resp. $\hat G_-$) denote the group of  holomorphic map $f:C\setminus\{0\}\to
GL(n,C)$ (resp. $f:\co_0\cup \co_\infty\to GL(n,C)$) satisfying the $U(n)$-reality condition,
where $\co_0$ and $\co_\infty$ are open neighborhood of $0$ and $\infty$ in $S^2$
respectively. A similar argument implies that the dressing action of $\hat G_-$ leaves
$\hat G_+$ invariant. We have seen in section 4 that the
dressing action of
$L_-(GL(n,C))$ on
$L_+(GL(n,C))$ is only defined locally.  
 However, we will show that the $U(n)$-reality condition implies that the
simple elements act on $G_+$ (resp. $\hat G_+$) globally and explicitly. Since simple
elements generate $G_-^m$, the group $G_-^m$ acts globally on $G_+$ (resp. $\hat G_+$).
We  explain these in more detail below.

\refclaim[ic] Proposition.  Let $z\in C$, $\pi$ a Hermitian projection of $C^n$ onto $V$,
$g_{z,\pi}$ a simple element of
$G_-^m$, and
$f\in G_+$ (resp. $\hat G_+$).  Then $g_{z,\pi}f$ can be factored as 
$$g_{z,\pi}f = \tilde f g_{z,\tilde \pi}\in G_-^m\times G_+\,\, ({\rm resp.\,}
G_-^m\times\hat G_+),$$ where $\tilde \pi$ is the Hermitian projection of $f(\bar
z)^{-1}(V)$.

\proof  Since $g_{z,\pi}=h_{\bar z,z,\pi}$, it follows from Proposition \refgu{} that if $f(\bar
z)^{-1}\cap f(z)(V^\perp)=0$ then we can factor 
$$g_{z,\pi}f= \tilde f h_{\bar z, z, \tilde \pi},$$ where $\tilde \pi$ is the projection onto
$f(\bar z)^{-1}(V)$ with respect to  $C^n=f(\bar z)^{-1}(V)\oplus f(z)^{-1}(V^\perp)$.
Since $f$ satisfies the $U(n)$-reality condition, $f(\bar \l)^*f(\l)=I$.  So we have 
\refeq[iy]$$\eqalign{<f(\bar z)^{-1}(V), \, f(z)^{-1}(V^\perp)> &= <f(z)^*(V), \,
f(z)^{-1}(V^\perp)> \cr &= <V, \, V^\perp>  = 0,\cr}$$ where
$<v_1,v_2>=v_1^*v_2$.  Since $< , >$ is positive definite, $$(f(\bar z)^{-1}(V))\cap
(f(z)^{-1}(V^\perp)=0.$$ So the factorization can also be done.  Equation \refiy{} also implies
that $f(\bar z)^{-1}(V)$ is perpendicular to  $f(z)^{-1}(V^\perp)$. Hence $\tilde \pi$ is the
Hermitian projection of $C^n$ onto $f(\bar z)^{-1}(V)$. So $h_{\bar z,z,\tilde \pi}$
satisfies the
$U(n)$-reality condition and
$h_{\bar z, z,\tilde \pi}= g_{z,\pi}$. \qed

\refclaim[cx] Theorem. The action $\sharp: G_-^m\times G_+\to G_+$
(resp. $G_-^m\times \hat G_+\to \hat G_+$) is globally defined, where 
$g\sharp f= \tilde f$ such that $\tilde f^{-1}gf\in G_-^m$.

\proof $G_-^m$ is generated by the simple elements. Hence the algorithm
for the factorization of the simple elements extends to all of $G_-^m$.
\qed

Henceforth in this section, let $\cu=u(n)$, and $\cm_{a,b,j}$ denote the space of solutions
of the $j$-th flow defined by $b$ in the $u(n)$-hierarchy on $\cs(R,\cu_a^\perp)$. 

 The construction of an action of
$G_-^m$  on  $\cm_{a,b,j}$ uses this dressing action.  We describe the action for
$j\geq 1$. The main difference between $j\geq 1$ and $j=-1$ is the difference between the
group $G_+$ and $\hat G_+$. The action on $\cm_{a,b,j}$ is induced from the ``dressing
action''  $\sharp$ of $G_-^m$ on $G_+$ using trivializations of elements in $\cm_{a,b,j}$. 

 Since the trivialization $E(x,t,\l)$ of $u\in \cm_{a,b,j}$ satisfies the $U(n)$-reality
condition, $E(x,t)\in G_+$, where $E(x,t)(\l)=E(x,t,\l)$.  The following
theorem is a consequence of Theorem \refhma{} and Proposition \refic{}. 

\refclaim[hm] Theorem. Let $u\in \cm_{a,b,j}$, $E$ the trivialization of $u$,
$z\in C\setminus R$, and $\pi$ the projection of $C^n$ onto  a complex linear subspace $V$ of
$C^n$.  For each $(x,t)\in R^2$,  set
$$\eqalign{ \tilde V(x,t) &=  E(x,t,z)^*(V),\cr \tilde \pi(x,t) &=\,{\rm the\,
projection\, of\,\,} C^n \,{\rm onto\,} \tilde V(x,t),\cr \tilde E(x,t,\l) &=
g_{z,\pi}(\l)E(x,t,\l) g_{z,\tilde
\pi(x,t)}(\l)^{-1}\cr &= (\pi+{\l-z\over \l-\bar
z}\pi^\perp)E(x,t,\l)\left(\tilde\pi(x,t)+{\l-\bar z\over \l- z} \tilde
\pi(x,t)^\perp\right).\cr }$$  Then 
\item {(i)} the smooth map $\tilde u$ from $R^2$ to
$\cu_a^\perp$ defined by $\tilde u= u + (z-\bar z)[\tilde\pi, a]$ is a solution of the $j$-th
flow equation in the $u(n)$-hierarchy and
$\tilde E$ is the trivialization of $\tilde u$,
\item {(iv)} $\tilde \pi$ satisfies 
\refeq[hq]$$\cases{(\tilde \pi)_x + [az+u,\tilde \pi] = (\bar z-z)[\tilde\pi,
a]\tilde\pi, &\cr (\tilde \pi)_t =\sum_{k=0}^j [\tilde \pi, 
Q_{b,j-k}(u)](z+(\bar z- z)\tilde\pi)^k,&\cr
\tilde\pi^*=\tilde\pi, \quad \tilde\pi^2=\tilde\pi, \quad \tilde\pi(0,0)=\pi,&\cr}$$

Next we want to prove that $\tilde u$ in Theorem \refhm{} is a solution of the $j$-th flow
in the $u(n)$-hierarchy and $\tilde u(\cdot, t)$ lies in the Schwartz class for all $t$. To do
this,  we need a theorem proved in [TU 1] (Theorem 6.6 of [TU1]). We will not repeat the
somewhat technical proof in  this paper. 

\refclaim[ba] Theorem ([TU 1]). Given $u\in \cs(R,\cu_a^\perp)$, if $\tilde \pi$ is a
solution of  
\refeq[av]$$(\tilde \pi)_x + [az+u,\tilde \pi] = (\bar z-z)[\tilde\pi, a]\tilde\pi,$$ 
then $[\tilde \pi, a]$ is in the Schwartz class, and $\lim_{x\to\pm\infty}\tilde
\pi(x,t)$ exists and commutes with $a$. 

As a consequence of Theorems \refhm{} and \refba{}, we have

\refclaim[bi] Corollary.  The function $\tilde u$ given in Theorem \refhm{} lies in
$\cm_{a,b,j}$, i.e., $\tilde u$ is a solution of the $j$-th flow  on $\cs(R,\cu_a^\perp)$
defined by $b$.

\ms
 To summarize, we have

\refclaim[bj] Corollary. Let $u\in \cm_{a,b,j}$, $E$ the trivialization of $u$,
$z\in C\setminus R$, $V$ a complex linear subspace of $C^n$, and $\tilde \pi(x,t)$ the
Hermitian projection of $C^n$ onto $E(x,t,z)^*(V)$. Then $\tilde u= u + (z-\bar z)[\tilde\pi,
a]$ is in  $\cm_{a,b,j}$, $\tilde E$ defined in Theorem \refhm{} is the trivialization of
$\tilde u$, and $g_{z,\pi}E(x,t)= \tilde E(x,t)g_{z, \tilde \pi(x,t)}\in G_+\times
G_-^m$.

 \refclaim[ht] Corollary. If  $u\in \cm_{a,b,j}$, then system \refhq{} is solvable. 
Moreover, if $\tilde \pi$ is a solution of system \refhq{} then $\tilde u= u+(z-\bar
z)[\tilde\pi,a]$ is again in $\cm_{a,b,j}$. \ei

The above results give methods to construct solutions of the $j$-th flow from
a given solution. This is done either by an algebraic formula if the trivialization of the given
solution is known (Darboux transformation) or by solving two compatible systems of ordinary
differential equations (B\"acklund transformation).

Corollary \refbi{} is part of the construction of an action of $G_-^m$  on
$\cm_{a,b,j}$:

\refclaim[as] Theorem. Let $g\in G_-^m$, $u\in \cm_{a,b,j}$, and $E$ the
trivialization of $u$. Then:
\item {(i)} $gE(x,t)$ can be factored uniquely as
$$gE(x,t)= \tilde E(x,t) \tilde g(x,t)\in G_+\times G_-^m$$
with $E(0,0)=I$. 
\item {(ii)} $\tilde E^{-1}\tilde E_x =\tilde  A$, where $\tilde A(x,t,\l)= a\l + \tilde
u(x,t)$ for some $\tilde u(x,t)\in \cm_{a,b,j}$.
\item {(iii)} $g\ast u= \tilde u$ defines an action of $G_-^m$ on $\cm_{a,b,j}$. 
\item {(iv)} $\tilde E$ is the trivialization of $\tilde u$. 

\proof   To prove uniqueness, we suppose $gE$ has two factorizations:
$$gE(x,t)=E_1(x,t)g_1(x,t)=E_2(x,t)g_2(x,t)\in G_+\times G_-^m.$$
Then 
\refeq[cz]$$E_1^{-1}(x,t,\l)E_2(x,t,\l)= g_1(x,t)(\l)g_2(x,t)^{-1}(\l).$$ But the left hand
side of
\refcz{} is holomorphic for $\l\in C$ and the right hand side is holomorphic at
$\l=\infty$ for all $(x,t)$. Hence by Liouville Theorem, it must be constant.
The right hand side at $\l=\infty$ is equal to
$I$, which proves the uniqueness. 

 By Theorem \refay{} the $g_{z,\pi}$'s generate $G_-^m$. To prove
the existence of the factorization it suffices to prove that we can factor
$g_{z,\pi}E(x,t)$. This is done in Theorem \refhm{}, since 
$$\tilde E(x,t)=g_{z,\pi}E(x,t)g_{z,\tilde \pi(x,t)}^{-1}$$ can be rewritten as
$$g_{z,\pi}E(x,t)=\tilde E(x,t) g_{z,\tilde \pi(x,t)} \in G_+\times G_-^m.$$ This
completes the proof of (i). 
\ss

Statement (ii) follows from Corollary \refbi{} and the
fact that
$g_{z,\pi}$'s generate $G_-^m$.
\ss

To prove $\ast$ defines an action, we need to prove 
$$(gh)\ast u= g\ast (h\ast u)$$ for $g,h\in G_-^m$ and $u\in \cm_{a,b,j}$.    Let $E$ denote the
trivialization of $u$. We factor $$\eqalign{hE(x,t) &=E_1(x,t)h_1(x,t)\in G_+\times G_-^m,\cr
gE_1(x,t)&=E_2(x,t)g_2(x,t)\in G_+\times G_-^m.\cr}$$ Then by definition of  $\ast$, 
$$\cases{E_1^{-1}(E_1)_x=
A_1.&\cr E_2^{-1}(E_2)_x= A_2,&\cr}$$ where 
$$\cases{A_1(x,t,\l)= a\l + (h\ast u)(x,t),&\cr A_2(x,t,\l)= a\l + (g\ast(h\ast
u))(x,t).&\cr}$$ But $(gh)E$ can be factored as
$$(gh)E= g(hE)=g(E_1h_1) =
E_2g_2h_1= E_2(g_2h_1)\in G_+\times G_-^m.$$ 
So by definition of $\ast$,  we have $$E_2^{-1}(E_2)_x = a\l + (gh)\ast u.$$ This proves
that 
$(gh)\ast u= g\ast (h\ast u)$. 
\qed

\refpar[ao] Definition. The transformation on $\cm_{a,b,j}$ defined by $u\mapsto g_{z,\pi}\ast
u$ is called a B\"acklund transformation for the $j$-th flow on $\cs(R,\cu_a^\perp)$
defined by $b$ with parameter $z$ and initial condition $\pi$. 

\ss

The center of $G_-^m$ is the subgroup of all elements in $G_-^m$ of the form $gI$ for
some $g:C\to C$ (here $I$ is the identity matrix) satisfying $g(\bar\l)^*g(\l)=1$.  We show
below that the center acts trivially on $\cm_{a,b,j}$.

\refclaim[ie] Corollary.  Suppose $g:C\to C$ is a rational function
such that
$g(\bar\l)^*g(\l)$ $ =1$, i.e., $gI$ lies in the center of $G_-^m$. Then $(gI) \ast u=u$ for all
$u\in
\cm_{a,b,j}$.  In particular,  if $g_1, g_2\in G_-$ such that 
$g_1=gg_2$, then
$g_1\ast u=g_2\ast u$.

\proof Let $E$ be the trivialization of $u\in \cm_{a,b,j}$ normalized at $(0,0)$. Since 
$(gI)E(x,t)= E(x,t)(gI)$, by definition $\ast$ we have $(gI)\ast u=u$.  
\qed

The B\"acklund theory for the $-1$ flow is identical. However, the trivialization $E$ has
singularities at both $0$ and $\infty$.  Since $G_-^m$ is always holomorphic on the
real axis, we do not need to add an extra condition at $0$. In this case, the trivialization $E$
is holomorphic in $C\setminus \{0\}$.  Use the same proof for the positive flows to this
case to yield

\refclaim[ki] Theorem. Let $\cm_{a,b,-1}$ denote the space of solutions of
the $-1$-flow on $\cs(R,\cu_a^\perp)$ defined by $b$:
$$\cases{u_t=[a,g^{-1}bg],&\cr g^{-1}g_x=u, \quad \lim_{x\to -\infty} g(x,t)=
I.&\cr}$$ Then the group $G_-^m$ acts on the space $\cm_{a,b,-1}$.  Moreover, let
$u\in
\cm_{a, b,-1}$,
$E$ the trivialization of $u$,  and
$g_{z,\pi}$ a simple element of $G_-^m$. Then: 
\item {(i)} $g_{z,\pi}\ast u= u+ (z-\bar z)[\tilde \pi, a]$, where $\tilde \pi(x,t)$ is the
projection of
$C^n$ onto $$E(x,t,z)^*(\pi(C^n)).$$
\item {(ii)} $\tilde \pi$ is the solution to
\refeq[kk]$$\cases{(\tilde \pi)_x + [az+u,\tilde \pi] = (\bar z-z)[\tilde\pi, a]\tilde\pi, &\cr
(\tilde \pi)_t ={1\over \n z\n^2} \left( (z-\bar z)\tilde \pi g^{-1}bg\tilde \pi - zg^{-1}bg\tilde \pi
+ \bar z \tilde
\pi g^{-1}bg\right),&\cr
\tilde\pi^*=\tilde\pi, \quad \tilde\pi^2=\tilde\pi, \quad \tilde\pi(0,0)=\pi.&\cr}$$

Now we turn to the twisted case. 
Given an involution $\s$ of $SU(n)$,  let $\cm^\s_{a,b,j}$ denote the space of solutions of the
$j$-th  flow on the subspace $\cs(R,\cu_{a,\s}^\perp)$ (here we use the same
notations as in section 3).  Let $G_-^{m,\s}$ denote the subgroup of $G_-^m$ of
$g\in G_-^m$ such that $\s(g(-\l))=g(\l)$.  Since the trivialization $E$ of $u\in
\cm^\s_{a,b,j}$ satisfies the same reality condition $\s(E(x,t,-\l))= E(x,t,\l)$, we
obtain:

\refclaim[cy] Corollary. For $j=-1$ or $j$ a positive integer, then the action of
$G_-^{m,\s}$ leaves $\cm_{a,b,j}^\s$ invariant, where $\cm^\s_{a,b,j}$ is the space of
solutions of the $j$-th flow in the $su(n)$-hierarchy twisted by an involution $\s$ defined by
$b$ on $\cs(R,\cu_{a,\s}^\perp)$.

Use a direct computation to get

\refclaim[bm] Proposition. Let $\s$ denote the involution on $SU(n)$ defined by
$\s(y)=(y^t)^{-1}$. Then
\item {(i)} $g_{z,\pi}\in G_-^{m,\s}$ if and only if $z=-\bar z$ and $\bar \pi=\pi$,
\item {(ii)} if $z\in C$ and $\bar \pi= \pi$, then $g_{z,\pi}g_{-\bar z,\pi}\in G_-^{m,\s}$.

\refpar[bk] Example. The trivialization of the vacuum solution $u=0$ in $\cm_{a,b,j}$ is 
$$E(x,t,\l)= e^{a\l x + b\l^j t}.$$ Suppose $U$ is a $n\times k$ matrix such that the columns of
$U$ form a basis of the linear subspace $V$ of $C^n$. By elementary linear algebra,  the
Hermitian projection of $C^n$ onto $V$ is $\pi = U(U^*U)^{-1} U^*$. Then Corollary \refbj{}
implies that 
\refeq[cb]$$g_{z,\pi}\ast 0 = (z-\bar z)\left[e^{-a\bar z x - b \bar z^j t}U(U^*e^{a (z-\bar z) x+
b(z^j-\bar z^j)t}U)^{-1} U^*e^{azx+bz^j t}, \,\, a\right]$$ is in $\cm_{a,b,j}$. These are the
$1$-solitons for the $j$-th flow in the $su(n)$-hierarchy.  So the space of $1$-solitons for the
$j$-th flow on
$\cs(R,\cu_a^\perp)$ defined by $b\in \cu_a$ is parametrized by the set 
$$\bigcup_{k=1}^{\left[{n\over 2}\right]} \left((C\setminus R)\times {\rm Gr\/}(k,n)\right).$$
 Here we use only Gr$(k,n)$ with $k\leq [{n\over 2}]$. This is  because
$g_{\bar z,\pi^\perp}= {\l-\bar z\over \l- z} g_{z, \pi}$ and Corollary \refie{} implies that
$g_{z,\pi}\ast 0= g_{\bar z,
\pi^\perp}\ast 0$. 

\ms
\refpar[bl] Example. Let $a=\diag(-i, i, \cdots, i)$, $z \in C\setminus R$, and
$\pi$  the Hermitian projection on the subspace spanned by  $(1, v)^t=(1, v_2,
\cdots, v_n)^t$. Then the one-solitons,  generated by B\"acklund transformations from the
vacuum solution, for the $j$-th flow on $\cs(R,\cu_a^\perp)$ defined by $a$ is 
$\pmatrix{0&B(x,t)\cr-B^\ast(x,t)&0\cr}$, where
$$B(x,t)= {4\Im (z) e^{2i({\rm Re\/}(z)x+{\rm Re\/}(z^j)t)}\,\bar v\over e^{-2(\Im(z) x +
\Im(z^j) t)}+e^{2(\Im(z) x +
\Im(z^j) t)}\N v\N^2}.$$
\ss

We turn now to the classical description of B\"acklund transformation. 
The classical B\"acklund transformations for the sine-Gordon equation are based on
ordinary differential equations:

\refclaim[ja] Theorem ([Da], [Ei]).  Suppose $q$ is a solution of the sine-Gordon
equation
\refdt{}, and $s\not=0$ is a real number. Then the following first order system is solvable
for $q^*$:
\refeq[jb]$$\cases {(q^*-q)_x= 4s \sin ({q^*+q\over 2})\cr 
(q^*+q)_t = {1\over s} \sin ({q^*-q\over 2}).\cr}$$
Moreover,  $q^*$ is again a solution of the sine-Gordon equation.

\refpar[jc] Definition. If $q$ is a solution of the sine-Gordon equation, then given any $c_o\in
R$ there is a  unique solution $q^*$ for  equation \refjb{}
such that $q^*(0,0)= c_o$.  Then $B_{s,c_o}(q)=q^*$ is a transformation on the space of solutions
of the sine-Gordon equation, which is the classical  B\"acklund transformation
for the sine-Gordon equation.  

\ss

We now relate the classical B\"acklund transformations and
the action of $G^{m,\s}_-$ on the space of solutions of the sine-Gordon equation (i.e., the
space $\cm^\s_{a,b,-1}$ with $\s$, $a, b$ defined as in Example \refbn{}).  Note if
$s\in R$,
$\tilde\pi ^*=\tilde \pi =(\tilde
\pi)^t$, then by Proposition \refbm{}, $g_{is,\tilde \pi}\in
G_-^{m,\s}$. Hence $\tilde \pi$ is a projection of $C^2$ onto $\pmatrix{\cos{f\over 2}\cr
\sin{f\over 2}\cr}$ for some  function $f$. In other words, 
$$\tilde\pi= \pmatrix{\cos^2{ f\over 2} &\sin { f\over 2}\cos { f\over 2}\cr\sin { f\over 2}\cos {
f\over 2}& \sin^2 { f\over 2}\cr}.$$  
So the first order system \refkk{}  for $\tilde\pi$  becomes 
\refeq[hta]$$\cases{f_x= {q_x\over 2} + 2s \sin f, &\cr f_t= {1\over 2s} \sin
(f-q).&\cr}$$ Write
$$\tilde u= g_{is,\b}\ast u= \pmatrix{0& \tilde q_x/2\cr -\tilde q_x/2&0\cr}.$$ But
$\tilde u= u + 2is[\tilde\pi,a]$, hence we have
$\tilde q= 2f-q$.  Writing equation \refhta{} in terms of $\tilde q$, we get 
$$\cases {(\tilde q-q)_x= 4s \sin ({\tilde  q+q\over 2})\cr 
(\tilde  q+q)_t = {1\over s} \sin ({\tilde  q-q\over 2}),\cr}$$
which is the classical B\"acklund transformation \refjb{} for the sine-Gordon equation. We
summarize this computation in the following Proposition:

\refclaim[jg] Proposition. Let $q$ be a solution of the sine-Gordon equation \refdt{}, and
$0<c_0<\pi$. Set 
$$\eqalign{u &= \pmatrix{0&{q_x\over 2}\cr -{q_x\over 2} & 0},\cr  f_o &={1\over 2} ( q(0,0)+
c_o),\cr}$$  and $\pi=$ the Hermitian projection onto the complex linear subspace spanned
by
$(\cos {f_0\over 2}, \sin {f_0\over 2})^t$, i.e., 
$$\pi=\pmatrix{\cos^2{ f_0\over 2} &\sin { f_0\over 2}\cos { f_0\over 2}\cr\sin { f_0\over
2}\cos { f_0\over 2}& \sin^2 { f_0\over 2}.\cr}$$  Then 
$B_{s,c_o}(q)= g_{is,\pi}\ast u$.

Notice the different roles of the parameter $s$, and the parameter $c_0$. The
location of the pole comes from $s$, and the angle of the projection comes from
$c_0$. The location of the pole is independent of choice of the base point in the
trivialization and is canonical, but the angle of the projection changes with the
trivialization. 

\bs

\newsection  Permutability formula.\par

There is a Bianchi permutability theorem for surfaces with Gaussian curvature $-1$ in
$R^3$, which gives the following analytical formula for the sine-Gordon equation:

\refclaim[je] Theorem ([Da], [Ei]). Suppose $q_0$ is a solution of the sine-Gordon
equation, $c_1,c_2\in R$, and $s_1, s_2\in R$ such that $s_1^2\not= s_2^2$ and
$s_1s_2\not=0$. Let
$q_i=B_{s_i,c_i}(q_0)$ for $i=1,2$. Then there exist $d_1,d_2\in R$,
which can be constructed algebraically,  such that 
\item {(1)} $B_{s_1,d_1}B_{s_2,c_2}= B_{s_2, d_2}B_{s_1,c_1}$,
\item {(2)} let $q_3= B_{s_1,d_1}B_{s_2,c_2}(q_0)$, then 
\refeq[jf]$$\tan {q_3-q_0\over 4}={s_1+s_2\over s_1-s_2}\tan {q_1-q_2\over
4}.$$  This is called the {\it Bianchi permutability formula\/} for the sine-Gordon
equation.   

The Bianchi permutability formula for  the sine-Gordon equation
is a consequence of factoring quadratic elements in the rational loop group $G_-^m$ in two
ways as product of two simple elements. We can also derive an analogue of the 
Bianchi permutability formula for the $j$-th flow using these two different ways of
factoring a quadratic elements in $G_-^m$. Note that contrary to the name, 
the permutability theorem does not mean that the
B\"acklund transformations generated by two simple elements are commuting.  

Notice that the singularities (both poles and zeroes) of an element in
$G_-^m$ comes in pairs
$(z,\bar z)$ due to the $U(n)$-reality condition.  In the $n\times n$ case, $n>2$, if
$\pi$ is the projection to $V$, $V=V_1+V_2$ and $V_1\perp V_2$, then 
$$ g_{z,\pi}= {\l-\bar z\over \l-z} g_{z,\pi_1}g_{z,\pi_2},$$ where $\pi_i$
is the projection to $V_i$. Hence there are infinitely many ways to factor
$g_{z,\pi_1}g_{z,\pi_2}$ as product of simple elements. 
 But we will prove later that a quadratic element with two different
singularities factors in exactly two ways. 

First we derive some relations among generators of $G_-^m$:

\refclaim[jp] Theorem.  Let $z_1, z_2\in C\setminus R$, and $\pi_1,\pi_2$
Hermitian projections of $C^n$.  If 
\refeq[ax]$$\phi=(z_2-z_1)I + (z_1-\bar z_1) \pi_1 - (z_2-\bar z_2) \pi_2$$ is
non-singular, then 
\item {(i)} $\tau_i =\phi\pi_i\phi^{-1}$ is a Hermitian projection for $i=1, 2$,
\item {(ii)} $g_{z_2,\tau_2}g_{z_1,\pi_1} = g_{z_1,\tau_1}g_{z_2,\pi_2}$, 
\item {(iii)} if $g_{z_2,\tau_2}g_{z_1,\pi_1} = g_{z_1,\tau_1}g_{z_2,\pi_2}$, then $\pi_1, \pi_2$
and $\tau_1, \tau_2$ are related as in (i). 

\proof Since $g_{z,\pi}=h_{\bar z,z,\pi}$, this theorem follows from Proposition \refgt{} if
we can prove that the $\tau_i$  is a Hermitian projection.
Since
$\tau_i^2=\tau_i$, we only need to prove that $\tau_i^\ast=\tau_i$ for $i=1,2$. To prove
this, we first set up some notations.   Let
$z_i=r_i + \sqrt{-1} \, s_i$.  Set
$\b_i=\sqrt{-1}\,(\pi_i-\pi_i^\perp)$ for $i=1,2$. Then $\b_i\in u(n)$, $\b_i^2=-I$, and 
$$g_{z,\pi_i} = {\l -r_i +s_i\b_i\over \l -r_i + \sqrt{-1} s_i}.$$
 Set $y=s_1\b_1-s_2\b_2$.  We claim that
$y^2\b_i=
\b_i y^2$.  To see this, we note that
$$y^2= -(s_1^2+s_2^2)I - s_1s_2(\b_1\b_2+ \b_2\b_1).$$ A direct computation 
gives  
$$\eqalign{(\b_1\b_2+\b_2\b_1)\b_1 &= \b_1\b_2\b_1 -\b_2,\cr 
\b_1(\b_1\b_2+\b_2\b_1) &= -\b_2+ \b_1\b_2\b_1.\cr}$$ So $y^2$ and $\b_1$
commute.  Similarly, $y^2$ and $\b_2$ commute. This proves our claim. 

Set  $\xi_i=\sqrt{-1}(\tau_i-\tau_i')$, $S_i= -r_i+s_i\b_i$, and $r=r_2-r_1$, where
$\tau_i'=I-\tau_i$.  Then 
$$\eqalign{\phi &=S_1-S_2 = rI +y, \cr \xi_i &= (rI+y)^{-1}\b_i (rI+y).\cr}$$
To prove $\tau_i$ is a projection is equivalent to prove that $\xi_i^2=-I$ and $\xi_i\in u(n)$. 
Since $\b_i^2=-I$ and $\xi_i$ is conjugate to $\b_i$,  $\xi_i^2=-I$. To prove that
$\xi_i^*=-\xi_i$, we compute directly
\refeq[jm]$$\xi_i^* = -((rI+ y)^*)^{-1} \b_i(rI+y)^*  =-(rI-y)^{-1}\b_i(rI-y).$$  But
$(rI-y)(rI+y)= r^2I - y^2$, which commutes with $\b_i$ because both $I$ and $y^2$ commute
with $\b_i$.  So we have
$(r^2-y^2)\b_i= \b_i (r^2-y^2)$, which implies that
\refeq[jn]$$(rI-y)^{-1}\b_i(rI-y) = (rI+y)\b_i (rI+y)^{-1}.$$
So $\xi_i^*=-\xi_i$.  Hence $\tau_i^\ast=\tau_i$.  
\qed

The following Proposition gives a sufficient condition on $z_1, z_2$ so that
$\phi$ defined by formula \refax{} is non-singular. 

\refclaim[jt] Proposition. Let $z_1=r_1+is_1$, and $z_2=r_2+is_2$, $\pi_1, \pi_2$
Hermitian projections, and  $\phi$ as in formula \refax{}. If $z_1\not= z_2$ and
$z_1\not=\bar z_2$, then
$\phi$ is non-singular.  

\proof Set $\b_i=\sqrt{-1}(\pi_i-\pi_i^\perp)$ as in the proof of Theorem \refjp{}. Then $\phi$
given by formula \refax{} can be written as
$$\phi= -(z_1-z_2)+2i(s_1\pi_1-s_2\pi_2)= -(r_1-r_2)+(s_1\b_1-s_2\b_2).$$ Since
$(s_1\b_1-s_2\b_2)\in u(n)$, its eigenvalues are pure imaginary. So if $r_1\not= r_2$, then
all eigenvalues of $ -(r_1-r_2)+(s_1\b_1-s_2\b_2)$ are not zero. 
	
Since $\b_i\in u(n)$ and $\b_i^2=-I$, $\N \b_i(x)\N = \N x\N$ for all
$x\in C^n$. So
$$\N (s_1\b_1-s_2\b_2)(x)\N \geq \n \n s_1\n-\n s_2\n \n \N x\N$$
for all $x\in C^n$.  If $s_1^2-s_2^2\not= 0$, then $(s_1\b_1-s_2\b_2)$ is
non-singular. Hence all eigenvalues of
$(s_1\b_1-s_2\b_2)$ are non-zero and pure imaginary, which implies
that  $ -(r_1-r_2)+(s_1\b_1-s_2\b_2)$ is non-singular. \qed

Use $g_{z,\pi}^{-1}(\l)= g_{z, \pi}(\bar\l)^\ast = g_{\bar z, \pi}(\l)$ and Theorem \refjp{} to
get

\refclaim[bb] Corollary. Given $g_{z_1, \pi_1}, g_{z_2,\pi_2}$ in $G_-^m$ such that
$z_1\not=z_2$ and
$z_1\not= \bar z_2$, then there exist uniquely Hermitian projections
$\tau_1, \tau_2$ such that 
$$g_{z_1,\pi_1}g_{z_2,\pi_2} =g_{z_2, \tau_2} g_{z_1,\tau_1},$$ where $\tau_i =
\phi \pi_i\phi^{-1}$ and $\phi=(\bar z_2-z_1)I +(z_1-\bar z_1)\pi_1 +(z_2-\bar
z_2)\pi_2$.

The following theorem follows easily from Theorem \refjp{} and Proposition \refjt{}: 

\refclaim[jpa] Theorem.  Let $z_1, z_2\in C\setminus R$ such that
$z_1\not= z_2$ and $z_1\not=\bar z_2$, and
$\pi_1,\pi_2$ Hermitian projections of $C^n$. Let $u_0\in \cm_{a,b,j}$ ($j=-1$ or $j\geq
1$),  and 
$$u_i=g_{z_i,\pi_i}\ast u_0=u_0+ (z_i-\bar z_i)[\tilde \pi_i, a]$$ for $i=1,2$ as given in
Theorem \refhm{}.  Set   
$$\eqalign{\phi &=(z_2-z_1)I + (z_1-\bar z_1) \pi_1 - (z_2-\bar z_2) \pi_2,\cr
\tilde \phi &= (z_2-z_1)I + (z_1-\bar z_1) \tilde\pi_1 - (z_2-\bar z_2) \tilde\pi_2,\cr
\tau_i &=\phi \pi_i \phi^{-1},\cr
\tilde \tau_i &=\tilde \phi \tilde \pi_i \tilde\phi^{-1}.\cr}$$
Then  
\refeq[jl]$$\eqalign{u_3 &= (g_{z_2,\tau_2}g_{z_1,\pi_1})\ast u_0= u_0+ (z_1-\bar
z_1)[\tilde\pi_1, a]  + (z_2-\bar z_2)[\tilde\tau_2,a] \cr &=
(g_{z_1,\xi_1}g_{s_2,\pi_2})\ast u_0 = u_0+ (z_1-\bar z_1)[\tilde\tau_1,a] 
+(z_2-\bar z_2)[\tilde \pi_2, a].\cr}$$

As a consequence of Proposition \refjg{} and Theorem \refjpa{}, we have

\refclaim[kn] Proposition. Formula \refjl{} for the $-1$-flow in the $su(2)$-hierarchy twisted
by the involution $\s(y)=-y^t$ on $\cs(R,\cu_{a,\s}^\perp)$ defined by $b={a\over
4}={1\over 4}\diag(i,-i)$ is the
 permutability formula \refjf{} for the sine-Gordon
equation.

\bs

\newsection $N$-soliton formula.\par

B\"acklund transformations defined in section 6 give an algebraic algorithm to
compute the solution
$g\ast 0$ of the $j$-th flow in the $su(n)$-hierarchy for a rational loop $g$.  The procedure
is as follows:
\item {(i)} Factor $g=g_{z_N,\pi_N} \cdots g_{z_1,\pi_1}$ as product of simple
elements.
\item {(ii)} Apply Theorem \refhm{} repeatedly to obtain $g\ast 0$.  In other
words, we set $u_0=0$,
$$E_0(x,t,\l)=e^{a\l x+ b\l^jt},$$ and define $u_j, \tilde \pi_j, E_j$
for $1\leq j\leq N$ by induction as follows:
$$\eqalign{&u_k = u_{k-1} +(z_k-\bar z_k)[\tilde \pi_k, a],\cr
 &\tilde \pi_k(x,t) = \, {\rm projection\, onto\/}\,\,
E_{k-1}(x,t,z_k)^*(\pi_k(V)), \cr
&E_k(x,t,\l) =
g_{z_k,\pi_k}(\l)E_{k-1}(x,t,\l)g_{z_k, \tilde
\pi_k(x,t)}(\l)^{-1}.\cr}$$ 
 
\noindent
 Then $g\ast 0= u_N$.

\ms
 Although this algorithm  is explicit,  it is difficult to write down $g\ast 0$ as a
formula in closed form in terms of $(z_1,\cdots, z_N,\pi_1, \cdots, \pi_N)$. 
We are motivated by the $2\times 2$ case considered in the book of Faddeev and
Takhtajan [FT] to use the permutability formula to  give  a formula
for $g\ast 0$ in closed form. 

To obtain the formula for $g\ast 0$, we first construct local coordinates for
$G_-^m$.   Note that $$g_{z,\pi}(\l)=\pi+ {\l-z\over \l- \bar z}\pi^\perp$$ has a
simple pole at
$\l=\bar z$,  and is holomorphic but not invertible at  $z$. We will call $\bar z$ a
zero of $g_{z,\pi}$. For $g\in G_-^m$, the  zeros and poles occur in pairs $(z_j,\bar
z_j)$.  Since it is more convenient to denote the pole of a simple element as $z$ in
our computation below, we change our notation for simple elements slightly.  
 Set 
$$h_{z,\pi}(\l) = I + { z-\bar z\over \l-z} \pi = \pi^\perp + {\l-\bar z\over \l-z}
\pi.$$  Then we have $$h_{z,\pi}(\l) = g_{\bar z, \pi^\perp}(\l) = \left({\l-\bar z\over
\l-z}\right) g_{z,\pi}.$$ 
Since the center of $G_-^m$ acts trivially on $\cm_{a,b,j}$ (by Corollary \refie{}),
$$h_{z,\pi}\ast u = g_{z,\pi}\ast u.$$ 

\ms
\refpar[da] Definition. A rational map $g\in G_-^m$ is called {\it regular\/} if $g$
has only simple poles and all the poles and zeros of $g$ are distinct.

\ms

\refclaim[se] Proposition. If $g$ is regular and $z$ is a simple pole of $g$,  then
there exists a unique projection $\pi$ such that $gh_{z,\pi}^{-1}$ is 
holomorphic and non-singular at $\l=z$. 

\proof Existence follows from Theorem \refay{} and Corollary \refbb{}. So it remains to
prove uniqueness. Assume $g=h_1h_{z,\pi_1}=h_2 h_{z,\pi_2}$, where $h_1,
h_2$ have no  zeros and poles at $z$. So $h_1(\bar z)$ and $h_2(\bar z)$ are
non-singular. But 
$$g(\bar z)= h_1(\bar z)\pi_1^\perp = h_2(\bar z)\pi_2^\perp.$$ Let $V_i=$
Im$(\pi_i)$. Then the above equation implies that 
$$h_1(\bar z)\pi_1^\perp(V_1)=0 =h_2(\bar z)\pi_2^\perp(V_1).$$ Since
$h_2(\bar z)$ is non-singular, $\pi_2^\perp(V_1)=0$. Hence $V_1\subset
V_2$.  Similarly, $V_2\subset V_1$.  This proves that $\pi_1=\pi_2$.
\qed
 
\ms
Given regular $g\in G_-^m$ with $N$ simple poles, there exist uniquely
$\G_g= (z_1, \cdots, z_N,\pi_1, \cdots,  \pi_N)$  such that $gh_{z_k,\pi_k}^{-1}$
is holomorphic at $z_k$ for all $1\leq k\leq N$.  We call $\G_g$
the {\it singularity data\/} of $g$.  

\bs
In the rest of this section, we will derive a formula for $g\ast 0$ in $\cm_{a,b,j}$
in terms of the singularity data $\G_g=(z_1, \cdots, z_N, \pi_1, \cdots, \pi_N)$ of
$g$.  Let 
$$e_{a,b,j}(x,t)(\l)=e^{a\l x+b\l^j}$$ denote the trivialization of the vacuum solution
$u=0$ in $\cm_{a,b,j}$. 

\refclaim[db] Proposition. Let $g\in G_-^m$ be a regular element with singularity
data
$\G_g=(z_1, \cdots, z_N, \pi_1, \cdots, \pi_N)$, and  $v_k$ a $n\times r_k$
matrix of  rank $r_k$ such that $\pi_k$ is the projection  onto the space spanned
by columns of
$v_k$. Factor $$ge_{a,b,j}(x,t) = E(x,t) \tilde g(x,t)\in G_+\times G_-^m$$ as in
Theorem \refas{}.  Then  
\refeq[sa]$$\tilde g(x,t)= I + \sum_{k=1}^N {P_k(x,t)\over \l-z_k},$$ where
$P_k(x,t)$ is an $n\times n$ matrix of rank $r_k$ for all $(x,t)\in R^2$.  Moreover, 
there exists smooth maps 
$\xi_k:R^2\to \cm_{n\times r_k}$ such that 
\item {(i)} $\xi_k(x,t)$ has rank $r_k$ for all $(x,t)$,
\item {(ii)} 
\refeq[dc]$$\sum_{m=1}^N {1\over z_m-\bar z_k}\, \xi_m(x,t)v_m^* \, e^{a(\bar z_m-
z_k)x+ b(\bar z_m^j-z_k^j)t} v_k = e^{-(az_k x+ bz_k^j t)}\, v_k$$ for $1\leq k\leq N$,
\item {(iii)} $P_k(x,t)=\xi_k(x,t)\,v_k^*\, e^{a\bar z_k
x+b\bar z_k^j t}$.

\proof   It follows from
Corollary \refbj{} that $\tilde g(x,t)$ is regular,  has only simple poles at $\l=
z_1, \cdots, z_N$, and is equal to $I$ at $\l=\infty$. So we can write $\tilde g$ in terms of
partial fractions:
$$\tilde g(x,t)(\l)= I + \sum_{k=1}^N {P_k(x,t)\over \l-z_k}$$ for some $n\times n$ matrix
function $P_k(x,t)$.  

First we claim that the rank of $P_k(x,t)$ is equal to $r_k$. By definition of the
singularity data,  we have
$g=h_kh_{z_k,
\pi_k}^{-1}$ for some
$h_k\in  G_-^m$ such that
$h_k$ is holomorphic and non-degenerate at $\l= z_k$. So it follows from
Theorem \refas{} and Corollary \refbj{} that there exist 
$\tilde h_k(x,t)$ in $G_-^m$ and projections $\tilde \pi_k(x,t)$ such that $$\tilde
g(x,t)=
\tilde h_k(x,t) h_{z_k,\tilde \pi_k(x,t)}$$ and $\tilde h_k(x,t)$ is holomorphic and
non-singular at
$\l=z_k$ for all $(x,t)\in R^2$.  Hence the residue of $\tilde g(x,t)(\l)$ at $\l=z_k$
is
$$ P_k(x,t)=(z_k-\bar z_k)\tilde h_k(x,t)(z_k)(\tilde \pi_k(x,t)).$$ 
So the rank of $P_k(x,t)$ is equal to that of $\tilde\pi_k(x,t)$,
which is $r_k$. This proves our claim.

It follows from Corollary \refbj{} that  $\tilde \pi_k(x,t)$ is the projection onto
the  space spanned by columns of
\refeq[dd]$$\tilde v_k(x,t) = e_{a,b,j}(\bar z)^* (v_k) =e^{-az_k x -bz_k^j
t}\,v_k.$$
 Since $gh_{z_k,\pi_k}^{-1}$ is holomorphic at $\l=z_k$, 
$$\tilde g(x,t) h_{z,\tilde \pi_k(x,t)}^{-1}= \left(I +\sum_{j=1}^N {P_j(x,t)\over
\l-z_j}\right)(I + {\bar z_k-z_k\over \l-\bar z_k}\tilde \pi_k(x,t))$$ is also holomorphic
at $\l= z_k$. So its residue at $z_k$ is zero. This implies that
$$P_k(x,t)\tilde \pi_k(x,t)^\perp=0.$$ 
Therefore the kernel of $P_k(x,t)$ contains the orthogonal complement of the
image of  $\tilde \pi_k(x,t)$. But the rank of $P_k(x,t)$ is $r_k$, which is the
rank of $\tilde \pi_k(x,t)$. So there exists $n\times r_k$ matrix $\xi_k(x,t)$ of
rank
$r_k$  such that
$$P_k(x,t)=\xi_k(x,t)\tilde v_k^\ast(x,t) = \xi_k(x,t) v_k^* \, e^{a\bar z_kx +
b\bar z_k^j t}.$$  

To prove (i) and (iii), it remains to prove that $\xi_k$ satisfies the linear system \refdc{}.  
Write $\tilde g(x,t) = \tilde h_k(x,t) h_{z_k,\tilde \pi_k(x,t)}$. 
Since $$h_{z_k,\tilde\pi_k(x,t)}(\bar z_k)= \tilde\pi^\perp_k(x,t),$$
we have 
\refeq[sc]$$\tilde g(x,t)(\bar z_k)(\tilde v_k(x,t)) =
 0$$ for all $(x,t)$ and $1\leq k\leq N$. 
Use formulas \refsa{} and \refsc{} to get
\refeq[de]$$\tilde v_k(x,t) + \sum_{m=1}^N{\xi_m(x,t)\tilde v^\ast_m(x,t)\tilde
v_k(x,t)\over \bar z_k- z_m}=0.$$
Substitute formula \refdd{} in \refde{} to get \refdc{}.
\qed

\ms

Suppose $g\in G_-$ is regular and $\G_g=(z_1, \cdots, z_N, \pi_1,
\cdots, \pi_N)$. We compute $g\ast 0$ when the rank of
$\pi_k$ is equal to $1$ for all $1\leq k\leq N$. In this case, 
$v_k,\xi_k$ are $n\times 1$. Set
$$f_{mk}(x,t)= {v_m^\ast\, e^{a(\bar z_m -z_k)x + b(\bar z_m^j-z_k^j)t}\, v_k
\over z_m-\bar z_k}.$$  Then  equation \refdc{} becomes
$$\sum_{m=1}^N  \xi_m(x,t)  f_{mk}(x,t)= e^{-(az_kx+bz_k^jt)}v_k, \quad 1\leq k\leq
N.$$ If the matrix $F(x,t)=(f_{mk}(x,t))$ is non-singular, then we can solve 
$$\xi_k(x,t) = \sum_{m=1}^N e^{-(az_mx+bz_m^jt)}v_m f^{mk}(x,t),$$ where
$(f^{mk})$ is the inverse of $F=(f_{mk})$. 

\refclaim[dg] Theorem.  Suppose $g\in G_-^m$ is regular with singularity data
$$\G_g=(z_1, \cdots, z_N, \pi_1, \cdots, \pi_N),$$ and $\pi_k$ is the
projection onto  $v_k\in C^n$ for each $1\leq k\leq N$.  Let $F(x,t)=(f_{km}(x,t))$
be the $N\times N$ matrix defined by
$$f_{mk}(x,t)= {v_m^\ast e^{a(\bar z_m -z_k)x + b(\bar z_m^j-z_k^j)t} v_k
\over z_m-\bar z_k}.$$  Suppose $F(x,t)$ is invertible. Set $F(x,t)^{-1}=(f^{km}(x,t))$ and 
 $$P_k(x,t)  =\left(\sum_{m=1}^N e^{-(az_mx+bz_m^jt)}v_mf^{mk}(x,t)\right)
v_k^\ast e^{a\bar z_kx+b\bar z_k^jt}.$$
Then the $N$-soliton $\tilde u= g\ast 0$ of the $j$-th flow on $\cs(R,\cu_a^\perp)$ defined
by $b$ and its trivialization $\tilde E(x,t)$ are given below: 
 $$\eqalign{&\tilde u = (g\ast 0)(x,t) = \sum_{k=1}^N \left[P_k(x,t),\, a\right],\cr
&\tilde E(x,t,\l) = g(0) e_{a,b,j}(x,t) \left(I +\sum_{k=1}^N
{P_k^\ast(x,t)\over \l-\bar z_k}\right).\cr}$$

\proof The formula for $\tilde E$ follows from Theorem \refas{}. By Theorem
\refas{} (ii), 
$g\ast 0$ is equal to the constant coefficient of the power series expansion of 
$\tilde E^{-1}\tilde E_x$ at $\l=\infty$.  Note that the expansion of $\tilde
g(x,t,\l)$ at
$\l=\infty$ is
$$\eqalign{\tilde g(x,t,\l) &= I +
\sum_{k=1}^N {P_k(x,t)\over \l-z_k} \cr &= I + \left(\sum_{k=1}^N
P_k(x,t)\right)
\l^{-1} + \cdots\cr}.$$ We will omit the variables $x, t$ in the following computation. So we
have
\refeq[du]$$\eqalign{&\tilde E^{-1}\tilde E_x = \l \tilde g a\tilde g^{-1} - \tilde g_x
\tilde g^{-1}\cr 
&\,=(I+\sum_{k=1}^N P_k \l^{-1} + \cdots)a\l (I-\sum_{k=1}^N
P_k\l^{-1} +\cdots) \cr &\quad + (\sum_{k=1}^N (P_k)_x \l^{-1} +
\cdots)(I-\sum_{k=1}^N P_k \l^{-1} + \cdots), \cr}$$ which is equal to $a\l+ \tilde u$. But
the constant term in equation \refdu{} is 
$\sum_{k=1}^N [P_k,a]$.    Hence
$\tilde u(x,t) = (g\ast 0)(x,t)=\sum_{k=1}^N \left[P_k(x,t),a\right]$. \qed

\bs

\newsection Scaling transformations.\par

The sine-Gordon equation is clearly invariant under the Lorentz transformations, for example,

\refclaim[jd] Proposition ([Da], [Ei]). If $q$ is a solution of the sine-Gordon equation, then
$L_r(q)(x,t)= q(r^{-1}x,rt)$ is also solution of
the sine-Gordon equation.  ($L_r$ is called a {\it Lie transformation\/} in the classical
literature). 

It is clear that $L_{r_1r_2}=L_{r_1}L_{r_2}$. In other word, Lie transformations give an action
of the multiplicative group $R^\ast$ of non-zero real numbers on the space of solutions of
the sine-Gordon equation.  The following result relating B\"acklund and Lie transformations is
known in classical surface theory (cf. [Da, Ei]).  

\refclaim[ko] Proposition ([Da], [Ei]). B\"acklund transformations and Lie transformations of
the sine-Gordon equation are related by the following formula:
$$B_{s,c_o}= L_s^{-1}B_{1,c_o}L_s.$$

It is known that an analogue of Lie transformations exists for the $j$-th flow, which
 will be called {\it scaling transformations\/}.  We describe these scaling transformations
next.  If
$u\in \cm_{a,b,j}$ and $r\in R^\ast$, then
 $$\hat u(x,t)=r^{-1}u(r^{-1}x, r^{-j}t)$$ is again a solution of the
$j$-th flow.  So $r\ast u= \hat u$ defines an action of $R^\ast$ on $\cm_{a,b,j}$.

The main goal of this section is to explain the relation between the scaling
transformation and B\"acklund transformations.  In fact, the
scaling transformation extends the action of $G_-^m$ to the action of the
semi-direct product $R^*\sdp G_-^m$ (defined below) on $\cm_{a,b,j}$, and
Proposition
\refko{} follows from the multiplication law of the group $R^\ast \sdp G_-^m$.
 
First we outline a proof for $r\ast u\in \cm_{a,b,j}$ if $u\in \cm_{a,b,j}$. 
Let $E$ be the trivialization of $u$, i.e.,
$$\cases{E^{-1}E_x= a\l + u, &\cr E^{-1}E_t = b\l^j + v_1\l^{j-1} + \cdots +
v_j,&\cr}$$ where $v_i= Q_{b,i}(u)$.  Set
$$\hat E(x,t,\l)= E(r^{-1}x, r^{-j}t, r\l).$$ 
Then $$\cases{\hat E^{-1}(x,t,\l)\hat E_x(x,t,\l) = a\l + \hat u(x,t), &\cr 
\hat E^{-1}\hat E_t = b\l^j + \hat v_1(x,t)\l^{j-1} +\cdots + \hat v_j(x,t), &\cr}$$
where
$\hat v_i(x,t)= r^i v_i(r^{-1}x, r^{-j}t)$ for $1\leq i\leq j$ and $v_i=Q_{b,i}(u)$.
If $b$ is a polynomial in $a$, then $Q_{b,i}(u)$ is a polynomial differential operator.
Hence $Q_{b,i}(u)$ vanishes at both $\infty, -\infty$, and the $\hat v_i$'s are in the
Schwartz class.  By Proposition
\refgm{}, we conclude $\hat u\in \cm_{a,b,j}$. 
(In fact, this calculation works for either positive or negative $r$, negative $r$
reverses $\pm\infty$). 

\ms

Next we define a one dimension extension of the group $G_-^m$. 

\refpar[kp] Definition.
Let $R^\ast=\{r\in R\n r\not=0\}$ denote the multiplicative group, and
$R^\ast\sdp G_-^m$  the semi-direct product of $R^\ast$ and $G_-^m$ defined by the
homomorphism
$$\rho:R^\ast\to {\rm Aut\/}(G_-^m), \qquad \rho(r)(g)(\l)= g(r\l),$$
i.e., the multiplication in $R^\ast\sdp G_-^m$ is defined by
$$(r_1,g_1)\cdot(r_2,g_2)= (r_1r_2, g_1 (\rho(r_1)(g_2))).$$

\ms
\refclaim[kd] Theorem. Suppose $j\geq 1$ or $j=-1$. Then the action $\ast$ of $G_-^m$
(resp. $G_-^{m,\s}$) extends to an action of
$R^\ast\sdp G_-^m$ (resp. $R^\ast \sdp G_-^{m,\s}$) on the space $\cm_{j,a,b}$ (resp.
$\cm_{a,b,j}^\s$) by
$$r\ast u (x,t) =r^{-1}u(r^{-1}x, r^{-j} t).$$ 

\proof  It is easy to see that $(r_1r_2)\ast u = r_1\ast(r_2\ast u)$.  Since 
$$(r,I)\cdot(1,g)(r^{-1}, I)= (1, \rho(r)(g)),$$ the action  $\ast$ extends to an action of
$R^+\sdp G_-^m$ if
\refeq[bd]$$(\rho(r)(g))\ast u= r\ast (g\ast (r^{-1}\ast u)).$$
To see this, we let $E$ be the trivialization of $u$, and define 
$$\eqalign{(r\ast E)(x,t,\l) &= E(r^{-1}x,r^{-j}t, r\l),\cr 
(g\ast E)(x,t,\l) &= \tilde E(x,t,\l)\cr}$$ for $r\in R^\ast$ and $g\in G_-^m$, where $\tilde E$ is
obtained from the factorization
$$gE(x,t)=\tilde E(x,t) \tilde g(x,t)\in G_+\times G_-^m$$ as in Theorem \refas{}.
To prove equation \refbd{}, it suffices to prove 
\refeq[bd]$$(\rho(r)(g))\ast E = r\ast (g\ast (r^{-1}\ast E)).$$   Write  
$$\eqalign{&g(\l) (r^{-1}\ast E)(x,t,\l)= E_1(x,t,\l) g_1(x,t,\l), \cr 
&(\rho(r)(g))(\l) E(x,t,\l)= E_2(x,t,\l) g_2(x,t,\l),\cr}$$ such that $E_i(x,t)\in G_+$
and $g_i(x,t)\in G_-^m$ for all $(x,t)\in R^2$ and $i=1, 2$.  Note that the second
equation gives
$$g(r\l)E(x,t,\l)= E_2(x,t,\l) g_2(x,t,\l).$$
By definition, we have
$$g\ast (r^{-1}\ast E(x,t))= E_1(x,t), \qquad (\rho(r)(g))\ast E(x,t)= E_2(x,t).$$
Now a direct computation gives
$$\eqalign{ &(r\ast (g\ast (r^{-1}\ast E)))(x,t,\l)\cr &= (g\ast (r^{-1}\ast
E))(r^{-1}x, r^{-j}t, r\l)\cr &= E_1(r^{-1}x, r^{-j}t, r\l)\cr 
&=g(r\l)(r^{-1}\ast E)(r^{-1}x, r^{-j}t, r\l) g_1^{-1}(r^{-1}x,r^{-j}t, r\l)\cr 
&=g(r\l) E(x,t,\l) g_1^{-1}(r^{-1}x, r^{-j}t, r\l) \cr 
&=E_2(x,t,\l)g_2(x,t,\l) g_1^{-1}(r^{-1}x, r^{-j}t, r\l) = E_2(x,t,\l) g_3(x,t,\l).\cr}$$
But $E_2(x,t)\in G_+$ and $g_3(x,t)=g_2(x,t) g_1^{-1}(r^{-1}x, r^{-j}t) \in G_-^m$. So
$$r\ast (g\ast (r^{-1}\ast E))= E_2,$$ which is equal to $(\rho(r)(g))\ast E$.  This proves our
claim. 
\qed

Since $(r^{-1},1)(1, g_{z,\pi})(r, 1)= (1, g_{rz,\pi})$, we have

\refclaim[kh] Corollary. Suppose $j\geq 1$ or $j=-1$.  If $u\in \cm_{a,b,j}$ (resp.
$\cm_{a,b,j}^\s$), then 
$$r^{-1}\ast (g_{z,\pi}
\ast (r
\ast u))= g_{rz,\pi}\ast u.$$ 

\refpar[bo] Remark.  Corollary \refkh{} for the $-1$-flow in the $su(2)$-hierarchy twisted by
$\s(y)=-y^t$   is Proposition
\refko{}.

\bs

\newsection B\"acklund transformations for $n$-dimensional systems.\par

	  The integrable equations of evolution we have been describing up to 
this point have two independent variables. The flow of the first
variable, regarded as a spatial variable, is used to construct the initial
Cauchy data.
The second variable is considered to be the time variable, and the flow in
this variable is the evolution.  In this section, we turn our attention to a
family of geometric problems in $n$ spatial variables, which we shall call
$n$-{\it dimensional systems\/}. In the applications,
the $n$ variables are on an equal footing, and the flows in each variable is a
first flow. The flows commute, and hence the resulting geometric object is
always a flat connection on a region of $R^n$ with  special properties. 

These $n$-dimensional systems  have been discussed in a paper by the first author
([Te]). We give definitions and some of the basic examples.  The results on B\"acklund
transformations developed in previous sections apply easily to these
systems. 
	
\refpar[cg] Definition ([Te]). Let $U$ be a rank $n$, semi-simple Lie group,  $\ct$ a
maximal abelian subalgebra of the Lie algebra $\cu$,  $a_1, \cdots, a_n$
a basis of $\ct$, and $\ct^\perp$ the orthogonal complement of $\ct$ with respect to
$(y_1, y_2)=\tr(y_1y_2)$.  {\it The
$n$-dimensional system associated to
$U$\/} is the following first order system:
\refeq[eea]$$[a_i,v_{x_j}]-[a_j, v_{x_i}] = [[a_i,v], [a_j,v]], 
\qquad v:R^n\to \ct^\perp.$$

\refpar[ch] Definition ([Te]). Let $U/K$ be a rank $n$ symmetric space,
$\s:\cu\to \cu$ the corresponding involution, $\cu=\ck+\cp$ the Cartan decomposition, $\ca$ a maximal
abelian subalgebra in $\cp$, $a_1, \cdots, a_n$ a basis of $\ca$, and $\ca^\perp$ the
orthogonal complement of $\ca$ in $\cu$.  {\it The
$n$-dimensional system associated to $U/K$ \/} is the first order system:
\refeq[ef]$$[a_i,v_{x_j}]-[a_j, v_{x_i}] = [[a_i,v], [a_j,v]], \qquad v:R^n\to
\cp\cap\ca^\perp.$$

\refclaim[ci] Proposition. The following conditions are equivalent: 
\item {(i)} $v$ is a solution of equation \refeea{} (or \refef{}) 
\item {(ii)} $[{\partial \over \partial x_i} +(a_i\l + [a_i,v]),
{\partial \over \partial x_j} +(a_j\l + [a_j,v])]=0$ for all $ i\not=j$.
\ei

\refpar[cr] Example.  Let $U/K=U(n)/O(n)$, and  $\cu=\ck+\cp$ the Cartan
decomposition corresponding to the involution
$\s(y)=-y^t$. Then $i\cp$ is the set of all real symmetric $n\times n$
matrices, and the space $\ca$ of all diagonal matrices in $\cp$ is a maximal
abelian subalgebra in $\cp$.  Let $e_{ii}$ denote the diagonal matrix such that all
entries are zero except the $ii$-th entry is equal to $1$.  Then $ie_{11},
\cdots, ie_{nn}$ form a basis of $\ca$.  The space $i(\cp\cap \ca^\perp)$ is
the space of all real symmetric $n\times n$ matrices whose diagonal
entries are zero. The $n$-dimensional system \refef{} associated to
$U(n)/O(n)$ for
$v=iF$ can be written as the system for 
$$F=(f_{ij}):R^n\to gl(n, R), \quad f_{ij}=f_{ji},
\quad f_{ii}=0\quad {\rm if \/}\,\, 1\leq i\leq n$$
\refeq[em]$$\cases{(f_{ij})_{x_i} + (f_{ij})_{x_j} + \sum_k f_{ik}f_{kj}=0, 
& if $i\not=j$,\cr 
(f_{ij})_{x_k}= f_{ik}f_{kj}, & if $i, j, k$ are distinct.\cr}$$
By Proposition \refci{},  $F$ is a solution of system \refem{} if and only if 
$$\left[{\p\over \p x_i} + e_{ii} \l + [e_{ii}, F], \,\, {\p\over \p x_j} + e_{jj} \l + [e_{jj},
F]\right] = 0$$ for all $i\not=j$. 
The $n$-dimensional system \refem{} is the equation for the Levi-Civita connection of an
Egoroff metric being flat. Here a metric $ds^2$ on $R^n$ is called  an {\it Egoroff metric\/}
if  it is of the form
$$ds^2=\sum_{i=1}^n \phi_{x_i}dx_i^2$$ for some smooth function 
$\phi:R^n\to R$. Set 
$$f_{ij}=\cases{{\phi_{x_ix_j}\over 2\sqrt{\phi_{x_i}\phi_{x_j}} },& if $i\not=j$,\cr
0, & if $i=j$.\cr}$$
It is easy to see that the Levi-Civita connection $1$-form for $ds^2$ is 
$$w_{ij}= -f_{ij}(dx_i-dx_j).$$ The metric
$ds^2$ is flat (i.e., $dw=w\wedge w$) if and only if $F=(f_{ij})$ is a solution of
system
\refem{}.

Since the $n$-dimensional system is the system consisting of $n$
commuting first flows, B\"acklund theory developed in section previous sections
for the first flow generalizes easily to that of the $n$-dimensional
systems (with minor changes).  

Let $\cm$ denote the space of solutions of the
$n$-dimension system
\refeea{} associated to
$u(n)$.  Given $v\in \cm$, the trivialization $E$ of $v$  is the solution of
$$\cases{E^{-1}E_{x_j}=a_j\l + [a_j,v],& $1\leq j\leq n$\cr E(0,\l)=I.&\cr}$$
Then $E(x,\bar\l)^*E(x,\l)=I$, i.e., $E(x)\in G_-^m$, where $E(x)(\l)=E(x,\l)$. So the action
of
$G_-^m$ leaves
$\cm$ invariant and the action of simple elements give B\"acklund transformations.

Let $\s(y)=(y^t)^{-1}$ be the involution of $U(n)$, and
$\cm^\s$ denote the space of solutions of the $n$-dimensional system \refem{}
associated to $U(n)/O(n)$. Let $G_-^{m,\s}$ denote the subgroup of $g\in G_-^m$ such
that $\s(g(-\l))=g(\l)$.  Since the trivialization $E$ of
$v\in \cm^\s$ satisfies the reality condition $$\s(E(x,t,-\l))=
E(x,t,\l),$$ we have  $E(x,t)\in G_-^{m,\s}$. So the action of $G_-^{m,\s}$ leaves
$\cm^\s$ invariant. Hence we obtain a B\"acklund theory for the system
\refem{}. 

In the two theorems below, we write down the analogous B\"acklund transformations and
Permutability formula for the $n$-dimensional system \refeea{}.  Given
$y\in gl(n)$, we will let 
$y_\ast$ denote
$y$ with the diagonal entries replaced by zeros.  

\refclaim[cj] Theorem. The group $R^\ast\times G_-^m$ acts on the space $\cm$ of
solutions of the $n$-dimensional system \refeea{} associated to $U(n)$, and  the
action $\ast$ is constructed in the same manner as on the spaces of solutions of the
first flow.  In fact, given $g_{z,\pi}\in G_-^m$ and $v\in \cm$,  the
following initial value problem is solvable for $\tilde \pi$ and has a unique solution:
$$\cases{(\tilde \pi)_{x_j} +[a_jz+ [a_j,v], \tilde \pi] = (\bar z-z)[\tilde \pi,
a_j]\tilde
\pi,&\cr  \tilde \pi^*=\tilde \pi, \quad \tilde \pi^2=\tilde \pi, \quad \tilde
\pi(0)=\pi.&\cr}$$  Moreover, 
\item {(i)} $g_{z,\pi}\ast  v = v - (z-\bar z) (\tilde \pi)_\ast$, 
\item {(ii)} the trivialization of $g_{z,\pi}\ast  v$ is $g_{z,\pi}E g_{z,\tilde \pi} ^{-1}$,
where $E$ is the trivialization of $v$
\item {(iii)} $\tilde \pi(x)$ is the projection onto the linear subspace $E(x,z)^*(V)$, where
$V$ is the image of the projection $\pi$,
\item {(iv)} $(r\ast v)(x)=r^{-1}v(r^{-1}x)$ for $r\in R^\ast$.
\item {(v)} if $U/K$ is the symmetric space defined by the involution $\s$, then the
group $R^\ast\sdp G_-^{m,\s}$ leaves the space
$\cm^\s$ of solutions of the $n$-dimensional system \refef{} associated to $U/K$
invariant. \ei

\refclaim[ck] Theorem.  Let $z_1, z_2\in C\setminus R$ such that
$z_1\not= z_2$ and $z_1\not=\bar z_2$, and
$\pi_1,\pi_2$ projections of $C^n$. Let $v_0\in \cm$,  and 
$v_i=g_{z_i,\pi_i}\ast v_0=v_0+ (z_i-\bar z_i)(\tilde
\pi_i)_\ast$ for
$i=1,2$ as given in Theorem \refcj{}.   Set   
$$\eqalign{\phi &=(z_2-z_1)I + (z_1-\bar z_1) \pi_1 - (z_2-\bar z_2) \pi_2,
\cr \tilde \phi &= (z_2-z_1)I + (z_1-\bar z_1) \tilde\pi_1 - (z_2-\bar z_2)
\tilde\pi_2,\cr
\tau_i &=\phi \pi_i \phi^{-1},\cr
\tilde \tau_i &=\tilde \phi \tilde \pi_i \tilde \phi^{-1}.\cr}$$
Then  $g_{z_2,\tau_2}g_{z_1,\pi_1}= g_{z_1,\tau_1}g_{z_2,\pi_2}$, and
\refeq[jl]$$\eqalign{v_3 &= (g_{z_2,\tau_2}g_{z_1,\pi_1})\ast v_0
= v_0+(z_1-\bar z_1)(\tilde\pi_1)_\ast  + (z_2-\bar
z_2)(\tilde\tau_2)_\ast
\cr &= (g_{z_1,\tau_1}g_{z_2,\pi_2})\ast v_0 = v_0+ (z_1-\bar
z_1)(\tilde\tau_1)_\ast + (z_2-\bar z_2)(\tilde
\pi_2)_\ast.\cr}$$

\bs

\newsection Time periodic solutions.\par

In this section, we use the action of $G_-^m$ to obtain many solutions of the $j$-th flow that
are periodic in time.  This is an algebraic calculation, which shows that when the poles are
properly placed, the solutions are periodic in time.  Multi-solitons will be time periodic if
the periods of the component solitons are rationally related.  We also show that the classical
breather solution of the sine-Gordon equation is obtained from the action of  a
suitable quadratic element at the vacuum.  

\refclaim[df] Theorem.  Let $j>1$ be an integer, $a=\diag(ia_1,
\ldots, ia_n)$, and
$b=\diag(ib_1, \ldots, ib_n)$.  If  $b_1, \ldots, b_n$ are rational numbers. Then the
$j$-th flow equation on $\cs(R,\cu_a^\perp)$ defined by $b$
has infinitely many $m$-soliton solutions that are periodic in $t$. 

\proof First assume $j\geq 3$.  Let $z=\rho e^{2\pi i\over j}$, $U$ a constant
$n\times k$ complex matrix of rank $k$, and $\pi$ the projection of $C^n$ onto the
linear subspace spanned by the $k$ columns of $U$. Then $z^j=\rho^j$ and formula
\refcb{} implies that 
$$g_{z,\pi}\ast 0 = (z-\bar z)\left[e^{-a\bar z x - b \rho^j t}U(U^*e^{a (z-\bar z) x}U)^{-1}
U^*e^{azx+b\rho^j t}, \,\, a\right].$$
Since $b_1, \cdots, b_n$ are rational
numbers, there exist $\tau>0$ and integers $m_r$ such that $$\tau= {2\pi m_r\over b_r\rho^j}
$$ for all $1\leq r\leq n$. So $g_{z,\pi}\ast 0$ is periodic in $t$ with period $\tau$. 
	
Let $\rho_1, \cdots, \rho_m$ be distinct rational numbers, and $z_k=\rho_k e^{2\pi
i\over j}$. Let $\tau_i$ denote the periods for $g_{z_i,\pi_i}$. Then it follows from the
details in the proofs of Theorems \refki{} and \refjpa{}  that the
$m$-soliton
$$(g_{z_1,\pi_1}\cdots g_{z_m,\pi_m})\ast 0$$ is an algebraic function of 
$g_{z_1, \pi_1}\ast0, \cdots, g_{z_m, \pi_m}\ast 0$.  Since $\rho_1, \cdots, \rho_m$
are rational numbers, there exist
$T>0$ and integers $k_i$ such that $T= k_i \tau_i$. So the
$m$-soliton  is periodic in time with period $T$.  

For $j=2$, let $z=is$ with $s\in R$.  Then the $1$-soliton
$$g_{z,\pi}\ast 0= 2is \left[e^{iasx + s^2 bt} v\left(v^*e^{2iasx }v\right)^{-1} v^*
e^{iasx-bs^2t}, \, a\right] $$
is periodic in $t$.  The existence of time-periodic
$m$-solitons of the second flow can be proved the same way as for the $j$-th flow.   
\qed

\refpar[eu] Example. Let $s\in R$,  and $c=(c_1, \cdots , c_{n-1})^t\in C^{n-1}$.
Recall that the second flow in the $u(n)$-hierarchy defined by $a=b=\diag(i, -i,
\cdots, -i)$ (Example \refad{} (3)) is the matrix non-linear Schr\"odinger equation for:
$$q_t={i\over 2}(q_{xx} + 2qq^*q), \quad q:R^2\to C^{n-1}.$$ Let $\pi$ be the Hermitian 
projection of
$C^n$  onto the complex
line spanned by $(1,c)^t$, where $c=(c_1, \cdots, c_{n-1})$.    The $1$-soliton solution
computed in Example
\refbl{} is
$$q=g_{is,\pi}\ast 0 = {4s \bar c e^{2is^2t}\over \n c\n^2 e^{2sx}+ e^{-2sx}}.$$ This
solution is periodic in $t$. 

\ms
The same algebra works for the $-1$ flow \refbt{}. Rewrite the $-1$ flow in terms
of $g$:
\refeq[cs]$$(g^{-1}g_x)_t =[a,g^{-1}ag].$$
Note that $(x,t)$ are characteristic coordinates. Let $$X=x-t,\quad T=x+t$$ be the
space-time coordinates.  Then \refcs{} in $(X,T)$-coordinate is
\refeq[if]$$(g^{-1}g_T)_T - (g^{-1}g_X)_X  + [g^{-1}g_X, \, g^{-1}g_T] = [a, \,
g^{-1}ag].$$ We will obtain solutions periodic in physical time (or
space). The trivialization of the vacuum solution for the $-1$-flow \refbt{} on
$\cs(R,\cu_a^\perp)$ defined by
$a=\diag(i,\cdots, i, -i,\cdots, -i)$ is $$E(\l,x,t)=\exp(a(\l x+ \l^{-1}t)).$$ By formula
\refcb{}, the
$1$-soliton 
$g_{e^{i\o},\pi}\ast 0$ for the $-1$-flow  is a function of 
$$\exp(i\cos\o (x+t) - \sin\o (x-t)) = \exp(i\cos\o X - \sin\o T).$$   This
proves

\refclaim[pa] Theorem. If $z=e^{i\o}$ and  $a=\diag(i,\cdots, i, -i,
\cdots, -i)$, then the $1$-soliton $g_{z,\pi}\ast 0$ for the
$-1$-flow \refif{} is periodic in time $T$ with period
${2\pi\over\cos\o}$. A multiple soliton generated by a rational loop with poles at
$z_1=e^{i\o_1}, \cdots, z_r=e^{i\o_r}$ will be periodic with period $\tau$ if there
exists integers $k_1, \cdots, k_r$ such that 
$$\tau={2\pi k_j\over \cos \o_j} \qquad \forall \,\, 1\leq j\leq r.$$

The multi-solitons above satisfy the sine-Gordon equation if the rational loop
satisfies
$(f(-\l)^t)^{-1}=f(\l)$, or equivalently $\overline{f(\bar \l)}= f(\l)$.    Now use
Theorem \refpa{} and Proposition \refbm{}, with 
$$z_1=e^{i\o_1}, \, z_2= -e^{-i\o_1}, \cdots, z_{2k-1}= e^{i\o_k}, \,
 z_{2k}=-e^{-i\o_k},$$ to get
a $2k$-soliton for the sine-Gordon equation that is periodic in time $T$.  To summarize, we
have

\refclaim[pb] Corollary. Multiple-breather solutions exists for the
sine-Gordon equation. 

\refpar[kv] Example.  If $\pi$ is a real symmetric projection (i.e., $\pi^2=\pi$, $\pi^\ast=\pi$
and $\bar\pi =\pi$), then
$$(g_{e^{i\o},\pi}g_{-e^{-i\o}, \pi})\ast 0= 4 \tan^{-1} \left( {\sin \o \sin ( (x+t)\cos\o)\over
\cos\o \cosh ((x-t)\sin \o)}\right).$$ This is the classical breather solution for the
sine-Gordon equation.  Theorem \refdg{} gives $m$-breather solutions
explicitly.  

\bs

\newsection The $u(k,n-k)$-hierarchies.\par

 We are able to obtain
global B\"acklund transformations for flows in the $su(n)$-hierarchy and the group
structure of these B\"acklund transformations because  the following three results:
\item {(i)} simple elements generate the rational group $G_-^m$, 
\item {(ii)} we can always do the Birkhoff factorization,
\item {(iii)} solutions to the ODE B\"acklund transformations in $x$-coordinate lie in the
Schwartz class.   

\ni Example \refdj{} shows that all three results fail to be true for  the
$sl(n,C)$- and $sl(n,R)$- hierarchies.  We will see that (i)-(iii) again fails for other
$\cu$-hierarchies when
$\cu$ is the Lie algebra of some non-compact group.  We give explicit examples for the
$u(1,1)$-hierarchy to explain this phenomenon.  However, our computation in fact
works for any real semi-simple Lie algebra. 

\ms
 Let  $J=\diag(\e_1, \cdots, \e_n)$ with $\e_i= 1$ for $1\leq i\leq k$ and $\e_j= -1$ if
$k<j\leq n$, and let $$\li v_1, v_2\ri_J = v_1^*Jv_2$$ denote the Hermitian bilinear form
on $C^n$ defined by $J$. Let $U(k,n-k)$ denote the group of linear maps of $C^n$ that
preserve $\li\,  , \,\ri_J$, and $u(k,n-k)$ its Lie algebra.
  Given a linear map $A:C^n\to C^n$, let
$A^{*J}$ denote the adjoint of
$A$, i.e.,
$A^{*J}$ is defined so that 
$$\li A(v_1), \, v_2\ri_J =  \li v_1,\, A^{*J}(v_2)\ri_J$$ for all $v_1, v_2\in C^n$. 
A direct computation shows that
$$A^{*J}=J^{-1}A^*J.$$ 
A projection $\pi$ of $C^n$ is called a $J$-{\it projection\/} if $\pi^{*J}=\pi$.   

It is easy to check that if $\pi$ is a $J$-projection then the simple element $h_{\bar z, z,
\pi}$ satisfies the $U(k,n-k)$-reality condition:
$$f(\l)^{-1} = J^{-1} f(\bar\l)^* J = f(\bar \l)^{*J}.$$
 Given a global solution $u$ of the
$j$-th flow in the $u(k,n-k)$-hierarchy, will the
$U(k,n-k)$-reality condition prevent the new solution having singularities?  Let
  $f= h_{z,\bar z,\pi}$. A direct computation gives
$$\li f(\bar z)^{-1}(V_1), \, f(z)^{-1}(V_2)\ri_J = \li f(z)^{*J}(V_1), \,
f(z)^{-1}(V_2)\ri_J =\li V_1, V_2\ri_J=0,$$
where $V_1$ and $V_2$ are image of $\pi$ and $I-\pi$ respectively.  
This implies that all vectors in $(f(\bar z)^{-1}(V_1))\cap  (f(z)^{-1}(V_2))$ are null
vectors with respect to $\li\, , \,\ri_J$.  Since there are non-zero null vectors with respect to
$< , >_J$, we can not always able to do the factorizations in the $U(k,n-k)$ case. Hence the
corresponding B\"acklund transformations may produce singular solutions.  In fact, the
following example shows this does happen.

\refpar[dn] Example.  Apply B\"acklund transformation to
the vacuum solution of the defocusing non-linear Schr\"odinger equation \refdk{} (Example
\refad{} (4)) to get solutions
$$\tilde u(x,t)= {-2i(z-\bar z)\over ((c+1)e^{-\xi(x,t)}-(c-1)e^{\xi(x,t)})} \pmatrix{0&
ibe^{-\eta(x,t)}\cr
-i\bar b e^{\eta(x,t)} &0 \cr},$$ where
$\xi(x,t)=i((z-\bar z)x+(z^j-\bar z^j)t)$ is real and  $\eta(x,t)=i((z+\bar z)x+(z^j+\bar
z^j)t)$ is pure imaginary ($b\in C$, $c$ is real, and $c^2-1=\n b\n^2$).  
Note that $\tilde u$ blows up at
$(x_0, t_0)$ when
$e^{2\xi(x_0,t_0)} = {c+1\over c-1}$. 

\ms

We have seen that the space of pure soliton solutions of the $j$-th flow in the
$su(n)$-hierarchy is the orbit of the group
$G_-^m$ through the vacuum (recall that $G_-^m$ is the group of all rational maps
$g:C\to GL(n,C)$ satisfying the
$SU(n)$-reality condition
$g(\bar\l)^*g(\l)= I$ and $g(\infty)=I$).   However, Example
 \refdn{} shows that the structure of the space of pure solitons
of the $j$-th flow in the
$\cu$-hierarchy is not clearly understood if $\cu$ is the Lie algebra of a non-compact Lie
group. 

\bs

\newsection The KdV hierarchy.\par

There is a formulation of KdV as  a restriction of the third flow.  The
odd flows in the
$sl(2,R)$-hierarchy for $q,r:R\to R$ leaves the submanifold defined by
$r=1$ invariant. The KdV equation is the third flow:
\refeq[eg]$$q_t={1\over  4}(q_{xxx} - 6qq_x),$$ and its Lax pair is
$$\left[{\p\over \p x} + a\l + u, \,\, {\p\over \p t} + a\l^3 +
u\l^2 + Q_2\l + Q_3\right]=0,$$ where 
$$\eqalign{ & a= \pmatrix{1&0\cr 0&-1\cr}, \quad u=
\pmatrix{0& q\cr 1& 0\cr}, \cr
& Q_2= \pmatrix{-{q\over 2}& -{q_x\over
2}\cr 0 & {q\over 2}\cr},\quad Q_3= \pmatrix{{q_x\over 4}& {q_{xx} -
2q^2\over 4}\cr -{q\over 2}& -{q_x\over 4}\cr}.\cr}$$
This Lax pair satisfies the $sl(2,R)$-reality condition.  But there is a second reality
condition that gives the restriction $r=1$.  To see this, let $$\phi(\l)=\pmatrix{1&\l\cr 0
&1\cr}.$$  Then the Lax pair of KdV satisfies:
\refeq[hc]$$\cases{\overline{A(\bar \l)}=A(\l),& \cr
\phi(\l)^{-1}A(\l)\phi(\l) =
\phi(-\l)^{-1} A(-\l) \phi(-\l). & \cr}$$
We will call this the KdV reality condition.  It is useful to realize that the second
condition is another way of saying $\phi^{-1}A\phi$ is even in $\l$.  So the
trivialization of a solution of the KdV equation normalized at $(0,0)$ satisfies the
same reality conditions \refhc{}, i.e., 
$$\cases{\overline{E(x,t,\bar\l)}=E(x,t,\l), &\cr
\phi(\l)^{-1}E(x,t,\l)\phi(\l) =  \phi(-\l)^{-1}E(x,t,-\l)\phi(-\l).&\cr}$$ 

In this section, we prove that B\"acklund transformations
of the KdV equation can be obtained in a similar way as before by
factoring the product of a degree one rational map and the trivialization
of a solution in the opposite order in the loop group.  Since the
factorization in the loop group of
$SL(2,R)$ can not always be carried out as we have seen in the previous
section, the same phenomenon is expected for the KdV equation.  

We will show that the Lax pair of all the odd flows in the
$sl(2,R)$-hierarchy with $r=1$ satisfies the reality conditions \refhc{}.  First, we obtain
the following lemma by a direct computation.

\refclaim[fa] Lemma.  $A(\l)= a\l+\pmatrix{\xi & q\cr r & \eta}$
satisfies the KdV-reality conditions \refhc{} if and only if $r=1$, $q\in R$ and 
$\xi=\eta\in R$, i.e., $A(\l)=a\l +\pmatrix{\xi&q\cr 1& \xi\cr}$. 

\refclaim[fb] Proposition.  The Lax pairs of the odd flows in the
$sl(2,R)$-hierarchy for $u=\pmatrix{0&q\cr 1&0\cr}$ satisfies the KdV-reality
conditions \refhc{}.

\proof  Note that the formal power series 
$$Q(\l) \sim a + Q_1\l^{-1} + Q_2\l^{-2} + \cdots$$ satisfies the equation
\refeq[gj]$$[d+a\l + u, Q(\l)] \sim 0.$$  Since $A(\l)=a\l+u$ satisfies the reality
conditions \refhc{}, we have
$$\eqalign{&\phi(\l)^{-1}\left[d_x+ A(\l), Q(\l)\right]\phi(\l) = 0\cr
&= \left[\phi(\l)^{-1}(d_x+A(\l))\phi(\l), \phi(\l)^{-1}Q(\l) \phi(\l)\right]\cr
&=\left[d_x + \phi(\l)^{-1}A(\l)\phi(\l),\phi(\l)^{-1}Q(\l) \phi(\l)\right]\cr
&= \left[d_x + \phi(-\l)^{-1}A(-\l)\phi(-\l),\phi(\l)^{-1}Q(\l)
\phi(\l)\right]\cr
&=\phi(-\l)^{-1} \left[d_x + A(-\l),\phi(-\l)\phi(\l)^{-1}Q(\l)
\phi(\l)\phi(-\l)^{-1}\right]\phi(-\l).\cr}$$
So $\left[d_x + A(-\l),\phi(-\l)\phi(\l)^{-1}Q(\l)
\phi(\l)\phi(-\l)^{-1}\right]= 0$.  It follows from a direct computation that the first two terms
of the asymptotic expansion of $$\phi(-\l)\phi(\l)^{-1}Q(\l)\phi(\l)\phi(-\l)^{-1}$$ is $(-a +
\l^{-1} u)$. So uniqueness of $Q$ (Proposition \refgm{}) implies that 
$$\phi(-\l)\phi(\l)^{-1}Q(\l)\phi(\l)\phi(-\l)^{-1}= -Q(-\l).$$ 
In particular, $\l^{2j+1}\phi(\l)^{-1}Q(\l)\phi(\l)$ is even, i.e., a power
series in $\l^2$.  This implies that 
$$a\l^{2j+1} + u\l^{2j} + Q_2\l^{2j-1} + \cdots +
Q_{2j+1}$$ satisfies the KdV reality condition \refhc{}.  \qed

There are no linear fractional transformations $g:C\to GL(2,C)$ 
satisfying the KdV- reality conditions \refhc{}.  But Corollary \refie{} tells us that 
rational loops with values in the center of $GL(n,C)$ act trivially on the space of solutions of
the $j$-th flow. So the group $G_-^{KdV}$ for constructing B\"acklund
transformations of the KdV equation is given as follows:

\refpar[pm] Definition. Let $G_-^{KdV}$ be the group of rational maps
$g:C\to GL(2,C)$ such that $g(\infty)=I$ and $g$ satisfies the KdV-reality
condition up to center elements. In other words,
$G_-^{KdV}$ is the group of rational maps $g:C\to GL(2,C)$ such that 
\item {(1)} $g(\infty)=I$, and
\item {(ii)} there exists some rational function 
$f :C\to C$ such that $fg$ satisfies the KdV-reality condition.  

\ms 
By Corollary \refie{}, if $a\l+u$ satisfies the KdV-reality
condition and $g\in G_-^{KdV}$ then both the trivialization $E$ of $a\l+u$ and $g\sharp E$
satisfy the KdV-reality condition. 

\ms
We need to find the simplest kind of elements in $G_-^{KdV}$.   
It follows from Lemma
\reffa{} that given any
$\xi,k\in R$
\refeq[eb]$$p_{\xi, k}(\l)=a \l +
\pmatrix{\xi&\xi^2-k^2\cr 1&\xi}= \pmatrix{\l+\xi & \xi^2-k^2\cr 1&
-\l+\xi\cr}$$
 satisfies the KdV-reality condition \refhc{}. (Recall that $a=\diag(1,-1)$). So
$$g(\l)={p_{\xi,k}(\l)\over
\l-k}\,\,\in G_-^{KdV}.$$ The inverse of $p_{\xi,k}$ is 
$$p_{\xi,k}^{-1}(\l) = {p_{-\xi, k}(\l)\over \l^2-k^2}.$$  We call
$\l_0$ a {\it zero\/} of $p_{\xi,k}$ if 
$\det(p_{\xi,k}(\l_0))=0$.  The proof of the next Proposition is a direct
computation. 

\refclaim[mn] Proposition. 
\item {(i)} $k, -k$ are the
only zeros of $p_{\xi, k}$, 
\item {(ii)} $p_{\xi,k}(\pm k)(v_\pm)=0$, where $v_+=\pmatrix{k-\xi\cr 1\cr}$ and 
$v_-=\pmatrix{-(k+\xi) \cr 1\cr}$,
\item {(iii)} Im$(p_{\xi, k}(k))$ is spanned by $av_-$ and Im$(p_{\xi,k}(-k))$ is
spanned by $av_+$. 

As a consequence we get

\refclaim[ng] Corollary.  Let $B=\pmatrix{k-\xi&-(k+\xi)\cr 1&1\cr}$, and
$Y=-kaBaB^{-1}$. Then $Y = \pmatrix{\xi&\xi^2-k^2\cr 1&\xi\cr}$ and
$a\l + Y= p_{\xi,k}(\l)$.

  We use a method similar to that of the $su(n)$
and $sl(n,R)$ hierarchies to construct the B\"acklund transformation
for the KdV equation corresponding to $p_{\xi,k}(\l)/(\l-k)$.  We give an
outline here. Let
$q$ be a solution of the KdV, and $E$ its trivialization at $(0,0)$.  First, we
take as an Ansatz, that there exists a map
$\tilde \xi(x,t)$ such that 
\refeq[et]$$\tilde E(x,t,\l)=p_{\xi,k}(\l)E(x,t,\l)p_{\tilde
\xi(x,t),k}^{-1}(\l)={p_{\xi,k}(\l)E(x,t,\l)p_{-\tilde
\xi(x,t),k}\over \l^2-k^2}$$ is holomorphic in
$\l\in C$. In other word, our Ansatz is that we can factor ${p_{\xi,k}\over
\l-k}\,E(x,t)$ as 
$${p_{\xi,k}\over \l-k}\,E(x,t)=\tilde E(x,t)\,{p_{\tilde \xi(x,t),k}\over
\l-k}\,\in G_+\times G_-^{KdV},$$ where $G_+$ is the group of holomorphic maps from
$C$ to $GL(2,C)$.  Since $E(x,t)$, $p_{\xi,k}$ and $p_{-\tilde \xi(x,t), k}$ satisfy the KdV-
reality condition \refhc{}, so is $\tilde E(x,t)$.  
 Since the residues of $\tilde E(x,t,\l)$ at $\l=k$ is zero and 
$$\tilde E(x,t,\l)= {p_{\xi, k}(\l)E(x,t,\l) p_{-\tilde \xi(x,t)}(\l)\over
\l^2-k^2},$$ we
get 
$$p_{\xi,k}(k)E(x,t,k)p_{-\tilde \xi(x,t), k}(k)=0.$$
By Proposition \refmn{}, we can choose $\tilde v_-(x,t)$ such that $E(x,t,k)a\tilde v_-(x,t)$
is proportional to $v_+$.  A similar calculation as for the $sl(n)$-hierarchy gives
the well-known Darboux and B\"acklund transformations for KdV:

\refclaim[ec] Theorem. Let $q$ be a solution of the
KdV equation,  and $E$ the trivialization of $q$ normalized at
$(x,t)=(0,0)$.  Given $\xi,k\in R$ with $k\not=0$,  set
$$\eqalign{&\pmatrix{f_1(x,t)\cr f_2(x,t)\cr} = E(x,t,k)^{-1}\pmatrix{k-\xi\cr
1\cr},\cr & \tilde \xi(x,t) = k- {f_1(x,t)\over f_2(x,t)},\cr &\tilde q(x,t)= -q + 2(
\tilde\xi^2(x,t)-k^2),\cr &
\tilde E(x,t,\l)= {p_{\xi,k}(\l) E(x,t,\l) p_{-\tilde
\xi(x,t),k}(\l)\over \l^2-k^2}.\cr}$$ If $f_2$
does not vanish in $\co\subset R^2$, then $\tilde q$ is a
solution of the KdV equation defined on $\co$ and $\tilde E(x,t,\l)$ is the trivialization of
$\tilde q$. 

\ms
\refpar[ia] Definition. Let $p_{\xi,k}\ast q$ denote the new solution
$\tilde q$ obtained in Theorem \refec{}. 

\ms

Using the same method as in previous sections, we compute $\tilde
E^{-1}\tilde E_x$ and $\tilde E^{-1} \tilde E_t$ to get the usual
ordinary differential equations for B\"acklund
transformations for KdV:

\refclaim[ei] Theorem.  Let
$k\in R$ be a constant. Then the following first order system for $A:R^2\to R$ is
compatible if and only if $q$ is a solution of KdV:
$$\cases{A_x=q -A^2+k^2,&\cr A_t=
{q_{xx}-2q^2\over 4} -{q_xA\over 2} + {q(A^2+k^2)\over 2} -
k^2(A^2-k^2),&\cr A(0,0)= \xi_0.&\cr}$$
Moreover, if $A$ is a solution of the above system, then the new solution is
 $$\tilde q= p_{\xi_0,k}\ast q = -q+
2(A^2-k^2).$$

\refpar[mm] Remark. Write $E(x,t,\l)=\pmatrix{y_1&y_2\cr z_1& z_2\cr}$. It
follows from Theorem \refec{} and a direct computation that 
$$\tilde \xi=  {f'\over f} \quad{\rm and\,} \quad\tilde q= -q + 2 {(f')^2\over f^2} -
2k^2 = f(f^{-1})'' - k^2,$$ where $f=y_1-(k-\tilde \xi)z_1$.  Since $y_1, z_1$ are
solutions of $y''=qy + k^2 y$, so is $f$.  
In other words, Theorem \refec{} in $x$-variable is the
classical Darboux transformation for the Hills operator ${d^2\over dx^2} -
q-k^2=0$ and the transformation maps a solution
$q$ of KdV to a new solution
$\tilde q$ is the Darboux transformation (cf. [AM], [De]).  Our result gives an
interpretation of Darboux transformations in terms of the Birkhoff factorization theorem.

\refpar[ed] Example.  We compute $p_{\xi,k}\ast 0$ for the KdV equation.  A direct
computation shows that the trivialization of the vacuum solution $q=0$ of KdV
normalized at $(0,0)$ is
$$E_0(x,t,\l)=\pmatrix{e^{(\l x+\l^3 t)}& 0\cr {\sinh(\l x+\l^3 t))\over
\l} & e^{-(\l x+\l^3 t)}\cr}.$$  Use Theorem \refec{}  and a direct
computation to get $\tilde q= p_{\xi,k}\ast 0 = 2(k^2-\tilde \xi^2)$, where 
$$\tilde \xi(x,t)= \cases{k\tanh(kx+k^3t +x_0),& if $b>0$,\cr
k\coth(kx+k^3t + x_0), &if $b<0$,\cr}$$  
$b=k^2-\xi^2$ and $x_0={1\over 2}\ln({\n\xi+k\n\over \n\xi-k\n})$.  So 
$$\tilde q(x,t)=\cases{-2k^2\sech^2(kx+k^3t +x_0), & if
$b>0$,\cr  -2k^2\csch^2(kx+k^3t +x_0), & if $b<0$.\cr}$$
Note that $\tilde q$ is a $1$-soliton solution if $b>0$,  and $\tilde q$ blows
up on the line $kx+k^3 +x_0=0$ if $b<0$. 

\ms

\refpar[mgs] Remark. The method we discussed above still 
works when we choose $k=0$. In fact, set the coefficients of $\l$ and $\l^2$
in the expansion of $$p_{\xi,0}(\l)E_0(x,t,\l)p_{-\tilde \xi(x,t),0}$$ in
$\l$ equal to zero to get $\tilde \xi(x,t)= {\xi\over 1+ \xi x}$ and a 
rational solution of KdV:
$$\tilde q(x,t)=2\tilde \xi^2={2\xi^2\over (1+\xi x)^2}.$$

\ms

Use a computation similar to that of Proposition \refgt{} to get:

\refclaim[ek] Proposition.  Given $a_1, a_2, k_1, k_2\in R$, if
$a_1-a_2\not=0$ then there exist uniquely $\xi_1, \xi_2$ such that
$$p_{\xi_2,k_2}p_{a_1,k_1}=p_{\xi_1, k_1}p_{a_2,k_2}.$$ Moreover, 
$\xi_1=-a_2+{k_1^2-k_2^2\over a_1-a_2}$ and  $\xi_2 =-a_1+
{k_1^2-k_2^2\over a_1-a_2}$.

As a consequence, we get the Permutability Formula for the KdV equation:

\refclaim[ej] Corollary. Suppose $q_0$ is a solution of KdV, and
$$q_i=p_{a_i,k_i}\ast q_0= -q_0 + 2(\xi_i^2- k_i^2)$$ for $i=1,2$. 
Set $\xi_{12}= -\xi_1 + {k_1^2-k_2^2\over \xi_1-\xi_2}$.
Then $$q_{12}= -q_1+ 2( \xi_{12}^2- k_2^2) = q_0 - 2(\xi_1^2-k_1^2)
+  2( \xi_{12}^2- k_2^2)$$ is again a
solution of KdV.

\refpar[ep] Example.  Assume $0<k_1^2<k_2^2<\cdots < k_n^2$.
Set $b_i=k_i^2-a_i^2$, and $$q_i=p_{a_i,k_i}\ast 0=
2(\xi_i^2-k_i^2)$$ for $1\leq i\leq n$.  Given a permutation $i$ of $1,
\cdots, n$, define
$\xi_{i_1i_2\cdots i_n}$ and
$q_{i_1i_2\cdots i_n}$ by induction:
$$\eqalign{\xi_{i_1i_2\cdots i_rjm}&= -\xi_{i_1i_2\cdots i_r j} +
{k_j^2-k_m^2\over\xi_{i_1i_2\cdots i_r j}- \xi_{i_1i_2\cdots i_r m}},
\cr q_{i_1i_2\cdots i_{r+1}} &= -q_{i_1i_2\cdots i_r} + 2(
\xi_{i_1i_2\cdots i_{r+1}}^2- k_{i_{r+1}}^2).\cr}$$  Use Corollary \refej{}
repeatedly, we conclude that $q_{12\cdots j}$ is a local solution of the
KdV  for all $1\leq j\leq n$. In general, these solutions may have
singularities.  We suspect that if  $b_1>0$ and 
$b_ib_{i+1}<0$ for $1\leq i\leq n-1$, then
$q_{12\cdots n}$ is smooth global solution. We can prove this for
$n\leq 5$, but do not have a proof for general $n$. It seems fairly clear that
these conditions are necessary for the solutions to be non-singular. 

\bs

\newsection The Kupershmidt-Wilson hierarchy.\par

The $sl(2,C)$-hierarchy \refdh{} leaves the submanifold $q=r$ invariant, and the third
flow is the complex modified KdV equation: 
\refeq[fc]$$q_t= {1\over 4}(q_{xxx}-6q^2q).$$
On this submanifold, the Lax pairs satisfy the following reality condition
$$\pmatrix{0&1\cr 1&0\cr}^{-1}A(-\l)\pmatrix{0&1\cr
1&0\cr}= A(\l).$$  Kuperschmidt and
Wilson [KW] found a natural 
$n\times n$ generalization of this restricted hierarchy
(KW-hierarchy).  First we describe the required reality condition. Let 
$e_{ij}$ denote
the matrix with zero on all entries except the $ij$-th entry is equal to $1$, and $\tau\in
GL(n)$ the matrix representing the cyclic permutation $(12 \cdots n)$, i.e., 
$$\tau=e_{21} +e_{32} +\cdots + e_{n, n-1} + e_{1n}.$$  Or
equivalently,
$\tau(e_i)=e_{i+1}$ for $1\leq i\leq n$ (here we use the convention that $e_i=e_j$ if
$i\equiv j$ mod $n$).  

The  reality condition for the $n\times n$ KW-hierarchy is  
\refeq[hd]$$\tau^{-1}A(\a^{-1}\l)\tau= A(\l), \quad {\rm where\,\,}  \a=e^{{2\pi
i\over n}},$$  which is called the {\it KW-reality condition\/}.
Since $\tau^n=I$, the order of the automorphism Ad$(\tau^{-1})$
on $gl(n,C)$ is $n$.  Let $\cg_k$ denote the eigenspace of Ad$(\tau^{-1})$
corresponding to eigenvalue $\a^k$ for $k=0, 1, \cdots, n-1$, i.e., $y\in \cg_k$ if
and only if $\tau^{-1}y\tau= \a^k y$. Or equivalently, $y=(y_{ij})\in \cg_k$ if and only if
$y_{i+1,j+1}=\a^k y_{ij}$ for all $1\leq i,j\leq n$.  Then 
$$gl(n,C)=\cg_0 + \cdots + \cg_{n-1}.$$
 For example, for $n=3$  we have
$$\eqalign{&\cg_0=\left\{\pmatrix{c_1 &c_2 &c_3\cr c_3& c_1 & c_2\cr 
c_2&c_3&c_1\cr}\bigg| c_i\in C\right\}, \quad \cg_1=\left\{\pmatrix{c_1 &c_2
&c_3\cr \a c_3&
\a c_1 & \a c_2\cr  \a^2 c_2&\a^2 c_3&\a^2 c_1\cr}\bigg| c_i\in C\right\}, \cr
&\cg_2= \left\{\pmatrix{c_1 &c_2 &c_3\cr\a^2  c_3&\a^2 c_1 &\a^2 c_2\cr 
\a c_2&\a c_3&\a c_1\cr}\bigg| c_i\in C\right\}.\cr}$$
 Because Ad$(\tau^{-1})$ is a Lie algebra
homomorphism, we have 
$$[\cg_i,\cg_j]\subset \cg_{i+j}.$$ Here 
$\cg_i=\cg_k$ if $i\equiv k$ mod $n$. 
A direct computation shows that $A(\l)=\sum_{k\leq n_0} u_k \l^k$
satisfies the KW-reality condition  \refhd{} if and only if $u_k\in \cg_k$ for all $k$. 

\refclaim[he] Proposition.  Let $\a=e^{2\pi i\over n}$, and $a=\diag(1, \a, \a^2, \cdots,
\a^{n-1})$.  Then the
$nk+1$-th flow in the $sl(n,C)$-hierarchy leaves $\cs(R, \cg_0\cap
sl(n)_a^\perp)$ invariant, and its Lax pair satisfies the KW-reality condition
\refhd{}.

\proof  Use a proof similar to that of Proposition \reffb{} to
conclude that $\tau^{-1}Q(\a^{-1}\l)\tau= \a Q(\l)$.  Hence
$$Q_{a,j}(u)\in \cg_{1-j}.$$ Since $a\in\cg_1$ and $[\cg_i, \cg_1]\subset
\cg_{i+1}$, we obtain $[Q_{kn+2}(u), a]\subset
\cg_{-kn}=\cg_0$.  \qed

\refpar[hf] Definition. The KW-equation is the restriction of the $(n+1)$-th flow in the
$sl(n,C)$-hierarchy to $\cs(R, \cg_0\cap sl(n)_a^\perp)$
\refeq[ib]$$u_t=(Q_{a,n+1}(u))_x + [u, Q_{a,n+1}(u)], \quad u:R^2\to  \cg_0\cap
sl(n)_a^\perp,$$  and the KW-hierarchy consists
of restricted $1$-st flow, $(n+1)$-th, $(2n+1)$-th, $\cdots$ flows in
the $sl(n,C)$-hierarchy.

\bs

 When $n=2$, $\cs(R,\cg_a^\perp\cap\cg_0)$ is the space of Schwartz class
maps from $R$ to $sl(2,C)$ of the form $\pmatrix{0&q\cr q&0\cr}$ and 
the third flow is the  complex modified KdV equation \reffc{}. 
For $n=3$, $\cs(R,\cg_0\cap sl(n)_a^\perp)$ is the space of Schwartz class
maps from $R$ to $sl(3,C)$ of the form $$\pmatrix{0&q_2&q_3\cr q_3&0 &q_2\cr
q_2&q_3&0\cr}.$$ The fourth flow is the KW-equation, which is of the form
$$(q_2)_t=P_2(q_2, q_3), \quad (q_3)_t=P_3(q_2, q_3),$$ where $P_2, P_3$ are
$4$-th order polynomial differential operators. The explicit formulas for $P_2$ and
$P_3$ are long, but they are not difficult to compute (use the method described in the proof
of Theorem \refiz{}). So we will not present them here. 

Since the Lax pair of the KW-equation satisfies the KW-reality condition \refhd{},
the trivialization of a solution of the KW-equation also satisfies the same reality
condition. Next we outline our strategy for finding B\"acklund transformations for
this equation.  We have seen that rational loops with values in the center of
$GL(n,C)$ do not play any effective role in the factorization (see Corollary
\refie{}). So the group $G_-^{KW}$ for constructing B\"acklund transformations is
defined as follows:

\refpar[pn] Definition. Let $G_-^{KW}$ denote the group of rational maps $f:S^2\to
GL(n,C)$ such that
\item {(i)} $f(\infty)=I$,
\item {(ii)} there exists a rational function $g$ such that
$gf$ satisfies the KW-reality condition.  

To construct B\"acklund transformations for the KW-equation, we 
start with a  degree one rational  map $g(\l)={a\l +Y\over \l-k}$ with $Y\in \cg_0$
as in the KdV case. Note that although $g$ does not satisfy the KW-reality
condition, $(\l-k)g(\l)=a\l+Y$ does. So $g\in G_-^{KW}$.  Let
$u$ be a local solution of the KW-equation that admits a reduced wave function $m$,
and $E$ the trivialization
 of $u$. Suppose at each
$(x,t)$ we can find  $\tilde Y(x,t)\in \cg_0$ such that 
$$\tilde E(x,t,\l)= {a\l+Y\over \l-k} E(x,t,\l)\left({a\l+\tilde Y(x,t)\over
\l-k}\right)^{-1}$$ is holomorphic in $\l\in C$. Then we can proceed  as in the proof
of Theorem
\refhma{} to conclude that
$\tilde E$ is the trivialization of some local solution of the $(n+1)$-th flow $\tilde u$ in the
$sl(n,C)$-hierarchy. But 
\refeq[hh]$$\tilde E(x,t,\l)=
(a\l+Y)E(x,t,\l)(a\l + \tilde Y(x,t))^{-1}.$$ Since all three terms in the right hand side
of \refhh{} satisfy the KW-reality condition, $\tilde E$ also satisfies \refhd{}.
Hence $\tilde E$ corresponds to a new solution
$\tilde u$ of the KW-equation.  However, in order to prove the expression \refhh{} is
holomorphic for all
$\l\in C$, we need to understand the relation between the zeros and kernels of $a\l+Y$
and the poles and residues of $(a\l+Y)^{-1}$.  We do this in the next few
Propositions. First we show that $a\l+Y$ is determined by a complex number $k$ and a vector
$v\in C^n$. 

\refclaim[ey] Proposition.  Let $Y\in \cg_0$, and 
 $f(\l)=a\l + Y$ (so $f$ satisfies the KW-reality condition \refhd{}). Then 
\item {(i)} there is a constant $k$ such that  det$(f(\l))= (-1)^{n+1}(\l^n-k^n)$,
\item  {(ii)} if $f(k)(v)=0$, then $f(\a^j k)(\tau^{-j}(v))=0$, 
\item {(iii)} $\tau^{-j}(v)$ are eigenvector of $a^{-1}Y$ with eigenvalues $-\a^j k$ for
$1\leq j\leq n-1$,
\item {(iv)}  if $v, \tau^{-1}(v),\cdots, \tau^{-(n-1)}(v)$ are linearly independent, then  
$Y=-kaBaB^{-1}$, where $B$ is the matrix whose $j$-th column is $\tau^{j-1}(v)$ for
$1\leq j\leq n$,
\item {(v)} $\det(Y)= (-k)^n$.

\proof Since $\tau^{-1}f(\a^{-1}\l)\tau=f(\l)$, $\deg(f(\a^{-1}\l)) =\deg(f(\l))$.
Hence $\det(f(\l))$ is a polynomial in $\l^n$. But the leading term of det$(f(\l))$ is 
$\a^{n(n-1)\over 2}\l^n$, which is equal to $(-1)^{n+1}\l^n$. This proves (i). The
rest of the Proposition follows from elementary linear algebra. \qed

\refpar[jv] Definition. Let $B$ denote the map
\refeq[gg]$$\eqalign{&B:C^n\to  gl(n, C)\quad {\rm defined \,
by\,}\cr &\,\, v\mapsto B(v)=(v,
\tau^{-1}(v),
\cdots,
\tau^{-(n-1)}(v)),\cr}$$ i.e., the $i$-th column of $B(v)$ is $\tau^{-(i-1)}(v)$ for
$1\leq i\leq n$. In other words, 
\refeq[ih]$$B(v) = B\pmatrix{v_1\cr v_2\cr
\cdot\cr \cdot \cr v_n\cr}=\pmatrix{v_1&v_2&\cdot&\cdot&v_n\cr v_2&v_3&\cdot&\cdot
&v_1\cr
\cdot& &\cdot &\cdot &\cdot\cr \cdot& &\cdot &\cdot &\cdot\cr v_n& v_1&\cdot
&\cdot &v_{n-1}\cr}.$$

\refpar[hj] Definition. Given $v\in C^n$ and $k\in C$, if $B(v)$ is non-singular, we
define $$p_{v,k}(\l)=a\l- kaB(v)aB(v)^{-1},$$ where $B(v)$ is the operator defined by
\refgg{} (or \refih{}). 

\ms
As a consequence of Proposition \refey{} we have

\refclaim[ji] Corollary.  Suppose $Y\in \cg_0$ and $a^{-1}Yv=-kv$ for some
non-zero vector $v$.  Then $f(\l)=a\l + Y = p_{v,k}(\l)$. Or equivalently, if
$f(\l)=a\l+ Y$ satisfies the KW-reality condition \refhd{} and $f(k)(v)=0$ then
$f=p_{v,k}$.

\refclaim[ez] Proposition.  Given $v\in C^n$ and $k\in C$,
 if $B(v)$ is non-singular, then
\item {(i)} $p_{v,k}(\l)$ satisfies the reality condition \refhd{},
\item {(ii)}  $p_{v,k}(\a^i k)\tau^{-i}v=0$ for $0\leq i\leq (n-1)$.

\proof Set $Y=-kaB(v)aB(v)^{-1}$. To prove (i), it suffices to prove
 $\tau^{-1}Y\tau=Y$. Note that
$a^{-1}Yv_j=-k\a^{j}v_j$, where $v_j=\tau^{-(j-1)}v$. Since $\tau v_j=v_{j-1}$
and $\tau^{-1}a\tau= \a a$, we get 
$$\eqalign{a^{-1}\tau^{-1}Y\tau v_j &=
a^{-1}\tau^{-1} a (a^{-1} Y)\tau v_j\cr &= a^{-1}\tau^{-1} a (a^{-1}Y)v_{j-1}\cr
&=-k\a^{j-1}(a^{-1}\tau^{-1}a)v_{j-1}= -k\a^{j}v_j.\cr}$$ 
This proves that $a^{-1}Y$ and $a^{-1}\tau^{-1}Y\tau$ have the same eigenvalues
and eigenvectors.  Hence $a^{-1}Y= a^{-1}\tau^{-1}Y\tau$, which implies
$Y=\tau^{-1}Y\tau$. This proves (i). 

By definition of $Y$, $(ak+Y)v=0$.  Since  $p_{v,k}(\l)=a\l+ Y$ satisfies \refhd{},
statement (ii) follows from  Proposition \refey{}. 
\qed

\refclaim[hk] Proposition.  Suppose $p(\l)=a\l+Y$ satisfies the KW-reality condition \refhd{}
and $k$ is a zero of $p(\l)$, i.e., det$(p(k))=0$.  Then
\refeq[hn]$$p(\l)^{-1} = {(\l+ \a a^{-1}Y) (\l + \a^2 a^{-1}Y) \cdots (\l + \a^{n-1}
a^{-1} Y) a^{-1}\over \l^n-k^n}.$$

\proof It follows from Proposition \refey{} that the eigenvalues of $a^{-1}Y$ are $-k, -\a
k,\cdots$,   $-\a^{n-1}k$. 
So $(a^{-1}Y)^n= (- k)^n$. But $$(\l+ z) (\l+\a z) \,\cdots \,(\l + \a^{n-1}z)= \l^n
+(-1)^{n+1}z^n.$$ Hence 
$$(\l+a^{-1}Y)\cdots (\l+ \a^{n-1}a^{-1}Y) = \l^n+ (-1)^{n+1}(a^{-1}Y)^n =
\l^n-k^n,$$  which finishes the proof. \qed

Next we  factor $p_{v, k}^{-1}$ as the product of simple elements. 
Note that 
$$\eqalign{&(\l+\a a^{-1}Y) \cdots (\l+\a^{n-1} a^{-1}Y)a^{-1} \cr 
&=(\l+\a a^{-1}Y) \cdots (a^{n-1}\l + \a^{n-1}a^{-1}Y a^{n-1})\cr
&=(\l+\a a^{-1}Y) \cdots  a^{n-2}(a\l+ \a^{n-1} a^{-(n-1)} Y a^{n-1})\cr
&=\cdots \cr 
&= (a\l+ \a a^{-1}Ya) (a\l + \a^2 a^{-2}Ya^2) \cdots (a\l+ \a^{n-1} a^{-(n-1)} Y
a^{n-1}).\cr}$$
Now suppose $k$ is a zero of $a\l +Y$ and $a^{-1}Yv = -kv$.  By Corollary \refji{},
$a\l+Y= p_{k,v}(\l)$.  Set $p_j(\l)=a\l + \a^j a^{-j}Ya^j$ for $1\leq j\leq n-1$.  Then 
$$\eqalign{p_j(\a^j k) (a^{-j} v) &= \a^j k a (a^{-j}v) + \a^j a^{-j} Y a^j(a^{-j}v)\cr
&=k\a^j a^{1-j}v + \a^j a^{-j}Yv\cr 
&= k\a^j a^{1-j}v + \a^j a^{1-j}(a^{-1}Y)v = 0.\cr}$$
 By Corollary \refji{} again, $p_j=p_{a^{-j}v, \a^j k}$.  As a consequence of this
computation and  Proposition
\refhk{}, we see that for $p(\l)= p_{v,k}(\l)$,  formula \refhn{} can be written as
\refeq[ho]$$p_{v,k}(\l)^{-1}={ p_{a^{-1}v,\a k}(\l) p_{a^{-2}v, \a^2k}(\l)
\cdot\,\cdot\,\cdot\,  p_{a^{-(n-1)}v,\a^{(n-1)}k}(\l)\over \l^n-k^n}.$$ 

\refclaim[mp] Proposition. Suppose the entries of $p,h:C\to GL(n,C)$ are polynomial such that
$p(\l)h(\l)=f(\l)I$ for some polynomial $f:C\to C$. If $f(k)=0$ for some $k\in C$, then
$$\Im(h(k))\subset \Ker(p(k)).$$

\proof Since $p(k)h(k)=0$, $\Im(h(k))\subset $ Ker$(p(k))$. \qed 

As a consequence of Propositions \refhk{} and \refmp{}, we have

\refclaim[hp] Corollary. Write  $p_{v,k}(\l)^{-1}={h(\l)\over
\l^n-k^n}$,
where $h(\l)$ is the degree $(n-1)$ polynomial  in the numerator of the formula
\refho{}. Then
 the image of $h(\a^i k)$ is the one dimensional space spanned by $\tau^{-i}v$ for
$0\leq i\leq (n-1)$.

Now we are ready to construct B\"acklund and Darboux transformations. 

\refclaim[ge] Theorem.  Let $v\in C^n$, $k\in C$ non-zero, $u:\co_1\to \cg_0\cap
sl(n)_a^\perp$ a solution of the KW-equation, and $E$ the trivialization of $u$ normalized at
$(0,0)$. Let 
$\tilde v(x,t)=E(x,t,k)^{-1}v$, and
$B$ the operator from $C^n$ to
$gl(n)$ defined by formula \refgg{} (or \refih{}).     If $B(v)$ is non-singular, then there
exists an open subset $\co$ of $\co_1$ such that $B(\tilde v(x,t))$ is non-singular for all
$(x,t)\in \co$.  Moreover,  
\item {(i)} $\tilde u=aua^{-1} + [\tilde Y, a]a^{-1}$ is again a solution of the
KW-equation defined on $\co$, where $$\tilde Y(x,t)=- kaB(\tilde v(x,t))
aB^{-1}\tilde v(x,t),$$
\item {(ii)} $\tilde E(x,t,\l)= p_{v,k}(\l)E(x,t,\l) p_{\tilde v(x,t),
k}(\l)^{-1}$ is the trivialization of $\tilde u$ normalized at $(0,0)$,
\item {(iii)} $Y$ is a solution of 
$$\cases{Y_x= Y u - (aua^{-1} + [Y, a]a^{-1}) Y,&\cr 
Y_t = YQ_{a,n+1}(u) -Q_{a,n+1}(aua^{-1} + [Y, a]a^{-1}) Y,\cr
\tau^{-1}Y\tau = Y,&\cr}$$ where
$Q_{a,n+1}$ is the polynomial differential operator defined in the
$sl(n,C)$-hierarchy.  

\proof  First we prove that $\tilde E$ is holomorphic for $\l\in C$. It follows from
formula \refho{} that $\tilde E$ is holomorphic for $\l\in C\setminus \{k, \a k, \cdots,
\a^{n-1}k\}$, and has possible simple poles at
$\l= \a^i k$ for $i=0, 1, \cdots, (n-1)$.   We claim that the residue of $\tilde E(x,t,\l)$
at $\l=\a^i k$ is zero. To see this we use formula \refho{} to write $$p_{\tilde
v(x,t),k}(\l) = {\tilde h(x,t,\l)\over
\l^n-k^n}.$$  For $0\leq i\leq (n-1)$, set $$f_i(\l)=(\l- k)
\cdots (\l-\a^{i-1}k) (\l-\a^{i+1}k) \cdots (\l-\a^{n-1}k),$$ i.e.,
$f_i(\l)={\l^n-k^n\over \l-\a^i k}$.  The residue  of $\tilde E(x,t,\l)$ at $\l=\a^i k$ is
equal to 
$${p_{v,k}(\a^ik)E(x,t,\a^i k)\tilde h(x,t, \a^i k)\over f_i(\a^i k)}.$$  But definition of
$\tilde v(x,t)$ implies $E(x,t,k)(\tilde v(x,t))= v$.  Since $E$ satisfies the reality
condition \refhd{}, we have 
$$E(x,t,\a^i k) = \tau^{-i}E(x,t,k)\tau^i.$$ By Corollary
\refhp{}, the image of $\tilde h(x,t,\a^i k)$ is the space spanned by $\tau^{-i}\tilde
v(x,t)$.  So the image of 
$$p_{v,k}(\a^ik)E(x,t,\a^i k)\tilde h(x,t, \a^i k)$$ is spanned by
$$p_{v,k}(\a^ik)\tau^{-i}E(x,t,k)\tau^i(\tau^{-i}\tilde v(x,t))= p_{v,k}(\a^i
k)\tau^{-i}v,$$ which is zero as follows from Proposition \refez{} (ii). This proves
that the residue of
$\tilde E$ is zero at $\a^ik$. Hence $\tilde E(x,t,\l)$ is holomorphic for $\l\in C$. 

 The rest of the theorem can be proved exactly the same as Theorem
\refhma{}. \qed

Relations among simple elements can be obtained by a direct computation as in Proposition
\refgt{}:
  
\refclaim[gh] Proposition.  Let $Y,Z\in \cg_0$ such that $(Y-Z)$ is non-degenerate. Set 
 $$\eqalign{\tilde Y &= a(Y-Z)a^{-1}Z(Y-Z)^{-1},  \cr
\tilde Z &= a(Y-Z)a^{-1}Y(Y-Z)^{-1}.\cr}$$ 
Then
\item {(i)} $\tilde Y, \tilde Z\in \cg_0$,
\item {(ii)} $(a\l+\tilde Y)(a\l + Y) = (a\l+ \tilde Z)(a\l+ Z)$.

The dimension of $\cg_0\cap sl(n)_a^\perp$ is $(n-1)$. So the KW-equation \refib{},
its B\"acklund transformations and permutability formula should be expressed in
terms of 
$(n-1)$ independent functions. Since an element in $\cg_0$ is determined by its first row, we
identify the space
$\cm_{1\times n}$ of
$1\times n$ complex matrices as $\cg_0$ via the linear isomorphism:
\refeq[jh]$$\eqalign{&\zeta:\cm_{1\times n}\to \cg_0\cr
& \, \, v=(v_1,\cdots, v_n)\mapsto \zeta(v), \,\, {\rm
where\, \,} (\zeta(v))_{ij}=v_{j-i+1}.\cr}$$
(Again $v_i=v_j$ if $i\equiv j$ mod $n$). 

Let $(0, q_2, \cdots, q_n)$ denote the first row of 
 $u\in  \cs(R, \cg_0\cap sl(n)_a^\perp)$, i.e., $u=\zeta((0, q_2, \cdots, q_n))$.
Let $(L_1(q), \cdots, L_n(q))$denote the first row of $Q_{a,n+1}(u)$. It follows from
Theorem \refiz{} that each
$L_j(q)$ is an order $n$ polynomial differential
operator in $q_2,\cdots, q_n$.
 Since $Q_{a,n+1}(u)\in \cg_0$, we can write 
$$Q_{a, n+1}(u)= \zeta(L_1(q), \cdots,
L_n(q)).$$   A direct computation implies that the KW-equation \refib{} written in terms of
$q$ is 
$${\p q_j\over \p t} = {\p \over \p x}L_j(q)+ \sum_{i=1}^n q_iL_{j-i+1}(q) -L_i(q)q_{j-i+1}
=  {\p\over \p x} L_j(q),$$ i.e., 
\refeq[hu]$${\p q_j\over \p t}= {\p\over \p x}(L_j(q)), \quad 2\leq j\leq n. $$

\ms
The Darboux transformation for the KW equation in Theorem \refge{} written in
terms of $q$ gives:

\refclaim[hr] Corollary.  Let $q= (0, q_2, \cdots, q_{n})$ be a solution of the $n\times n$
KW-equation \refhu{}, and $E$ the trivialization normalized at $(0,0)$. Given $v\in C^n$ and
$k\in C$, let $\tilde v(x,t)=E(x,t,k)^{-1}(v)$, and $B$ the operator defined by the formula
\refih{}. Let
$y(x,t)=(y_1(x,t), y_2(x,t), \cdots, y_{n}(x,t))$ denote the first row of the matrix 
$$-kaB(\tilde v(x,t))aB(\tilde v(x,t))^{-1},$$ and  $\tilde q=
(q-y)a^{-1} + y$, i.e., 
$$\tilde q_j=\a^{1-j}q_j +(1-\a^{1-j})y_j, \quad {\rm for\,}\, 2\leq j\leq n.$$ Then
$\tilde q=(0,\tilde q_2, \cdots, \tilde q_n)$ is again a solution of the KW-equation \refhu{}.

The ODE version of B\"acklund transformations for the KW-equation in Theorem \refge{}
written in terms of
$q$ gives:

\refclaim[hv] Corollary.  Suppose $q=(0, q_2, \cdots, q_n)$ is a solution of the KW-equation
\hfill\break
\refhu{}. Then the following systems for $y=(y_1,
\cdots, y_n)$ are compatible:
$$\cases{(y_j)_x= \sum_{i=1}^n (1-\a^{i-j})y_i(q_{j+1-i}-y_{j+1-i}),&
\cr (y_j)_t= \sum_{i=1}^n
L_{j+1-i}(q)y_i-L_i((q-y)a^{-1}+ y)y_{j+1-i}.&\cr}\eqno(BT^{KW}_q)$$ Moreover, if
$y(x,t)$ is a solution of $BT^{KW}_q$, then
\item {(i)} $\tilde q= (q-y)a^{-1} + y$ is again a solution of the KW-equation,
\item {(ii)} $\det(\zeta(y(x,t))$ is a constant.

As a consequence of Proposition \refgh{}, we have

\refclaim[hx] Corollary.  Suppose $q=(0,q_2, \cdots, q_n)$ is a solution of \refhu{}, $\xi,
\eta$ are solutions of $BT^{KW}_q$, and
$$q'=(q-\xi)a^{-1}+\xi, \quad q''=(q-\eta)a^{-1}+ \eta$$ are the
corresponding new solutions of the KW-equation.  Assume
$\det(\zeta(\xi-\eta))\not=0$, where $\zeta$ is the operator defined by formula
\refjh{}.  Set  
\refeq[hz]$$\eqalign{\tilde \xi &=
(\xi-\eta)a^{-1}\zeta(\eta)a^{-1}(\zeta(\xi-\eta))^{-1},\cr
\tilde \eta &= (\xi-\eta)a^{-1}\zeta(\xi)a^{-1}(\zeta(\xi-\eta))^{-1}.\cr}$$
Then 
 $\tilde \xi$ is a solution of $BT^{KW}_{q'}$,
$\tilde \eta$ is a solution of $BT^{KW}_{q''}$, and
\refeq[hy]$$\eqalign{\tilde q &= q'a^{-1}+\tilde \xi(I-a^{-1})=
(\xi+(q-\xi)a^{-1})a^{-1} +\tilde \xi (I-a^{-1})\cr  &= q''a^{-1}+ \tilde \eta
(I-a^{-1})= (\eta+(q-\eta)a^{-1})a^{-1} +
\tilde \eta (I-a^{-1})\cr}$$ is again a solution of the KW-equation.

\refclaim[mh] Corollary.  Let $q, q', q''$ and $\tilde q$ be as in Corollary \refhx{}. Then $\tilde
q$ is an algebraic function of $q, q', q''$. 

\proof
Use formulas \refhz{} and \refhy{} to write $\tilde q$ in terms of
$q, \xi$ and $\eta$.  Corollary \refhr{} implies
$$\xi_j={q'_j-\a^{1-j}q_j\over 1-\a^{1-j}}, \quad \eta_j={q''_j-\a^{1-j}q_j\over
1-\a^{1-j}}, \quad 2\leq  j\leq n.$$ So $\xi_2, \cdots, \xi_n, \eta_2, \cdots, \eta_n$ are
algebraic functions of $q, q', q''$.  Corollary \refhv{} (ii) implies that $\det \zeta(\xi)=c_1$
and $\det\zeta(\eta)=c_2$ are constant. 
But $\det\zeta(\xi)$ (resp. $\det\zeta(\eta)$) is a degree $n$ polynomial in
$\xi_1$ (resp. $\eta_1$). This implies that $\xi_1$ (resp. $\eta_1$) can be written as an
algebraic function of $\xi_2, \cdots, \xi_n$ ( resp. $\eta_2, \cdots, \eta_n$). 
Hence $\tilde q$ is  an algebraic function of $q, q'$ and
$q''$.  \qed

\bs

\newsection The Gel'fand-Dikii Hierarchy.\par

The Gel\'fand-Dikii ($GD_n$-) hierarchy is a hierarchy of flows on the  space
$\cp_n$  of $n$-th order scalar differential operators
$$L= D^n -(p_1D^{n-2} + p_2 D^{n-3} +\cdots + p_{n-1}),$$  where $D={d\over
dx}$ and $p_i\in \cs(R, C)$. Flows in this hierarchy  are given by
$${\p L\over \p t}= [L^{j/n}_+, L],$$ where $L^{j/n}_+$ is the differential 
operator part of the pseudo-differential operator $L^{k/n}$.  The spectral
problem 
$Ly_1=\l^n y_1$ is equivalent to the spectral problem of the following first
order system for $y=(y_1,
\cdots, y_n)$: 
\refeq[mq]$$\eqalign{{d\over dx}(y_1, \cdots, y_n) &= (y_1, \cdots,
y_n)(e_{1n}\l^n + b + v),\cr &=(y_1, \cdots, y_n)
\pmatrix{0&\cdot
&\cdot  &\cdot & p_{n-1}+\l^n\cr 1&0 &\cdot
 &\cdot & p_{n-2} \cr 0&1 &0  &\cdot & p_{n-3} \cr 0&& &
\cr 0&\cdot&\cdot &1 &0\cr}.\cr}$$ Here
\refeq[ka]$$\eqalign{ b&=e_{21} +\cdots + e_{n n-1}=\pmatrix{0&\cdot
&\cdot  &\cdot & 0\cr 1&0 &\cdot
 &\cdot & 0\cr 0&1 &0  &\cdot & 0 \cr 0&& &
\cr 0&\cdot&\cdot &1 &0\cr}\cr
 v &=\pmatrix{0&0&\cdots& p_{n-1}\cr 0&0&\cdots& p_{n-2}\cr
&&\cdots&\cr 0&0&\cdots &p_1\cr 0&0&\cdots &0\cr},\cr}$$  and $e_{ij}$
is the $n\times n$ matrix whose $ij$-th entry is $1$ and all other
entries are $0$.  The
$GD_n$-hierarchy gives rise to a hierarchy on the space $\cm_n$ of all
$v:R\to sl(n,C)$ of the form \refmq{} with $p_i\in \cs(R,C)$ for 
$1\leq i\leq n-1$. 
 But unlike all the hierarchies we have discussed in previous sections, when
$n\geq 3$,
$\cm_n$ is not determined by a reality condition and the  corresponding flows
on $\cm_n$
 are not the restriction of the flows in the
$sl(n,C)$-hierarchy to $\cm_n$.  Drinfeld and Sokolov [DS1, 2] gave a 
description of the symplectic structures  and the flows on $\cm_n$ using a
symplectic quotient. 

When $n=2$, the $GD_2$-hierarchy is the complex KdV hierarchy.  In section
12, we saw that this hierarchy is obtained by restricting the
$sl(2,C)$-hierarchy to the submanifold that is defined by the reality
condition:
\refeq[kb]$$\phi(\l)^{-1}A(\l)\phi(\l)=\phi^{-1}(-\l)A(-\l)\phi(-\l), \quad
\phi(\l)=\pmatrix{1&\l\cr 0&1\cr}.$$  In other words, $A(\l)$ lies in the Lie
subalgebra of fixed points of the involution $\o$ defined by
$$\o= \Ad(\phi)\circ \tau \circ\Ad(\phi)^{-1},$$ where $\tau(A)(\l)= A(-\l)$. 
In fact, we showed that the odd flows in the
$sl(2,C)$-hierarchy leaves  invariant the space $\cs_\phi$ of all $u:R\to
sl(2)_a^\perp$ such that 
$A(\l)=a\l+ u$ satisfies the reality condition \refkb{}. The $GD_2$- hierarchy
is the $sl(2,C)$-hierarchy restricted to $S_\phi$.   
The main purpose of this
section is to generalize this construction to the $GD_n$-hierarchy.
However, we need to use a different
$gl(n)$-valued first order linear operator ${d\over dx} + A(\l,x)$ than the
one given by the formula \refmq{}.  Here $A(\cdot, x)$ is fixed by certain 
order $n$ Lie algebra homomorphism $\s_n$:
$$\s_n=\Ad(\phi_n)\circ \tau_n\circ\Ad(\phi_n)^{-1}.$$ 
 To motivate the choice of $\phi_n$,
we  first explain the relation between  $$\phi_2(\l)=\phi(\l)=\pmatrix{1&\l\cr
0&1\cr}$$ and the phase space of the
$GD_2$-hierarchy.  The vacuum
$L_0=D^2-\l^2$ corresponds to
$${d\over dx}+ A_0(\l), \quad {\rm where \,\,} A_0(\l)=\pmatrix{0&\l^2\cr
1&0\cr}.$$  $A_0(\l)$ can be diagonalized by $V(\l)=\pmatrix{1&\l\cr 1&-\l\cr}$:
$$V(\l)\pmatrix{0&\l^2\cr 1&0\cr} V(\l)^{-1}
= \pmatrix{1&0\cr 0&-1\cr}\l.$$ 
Now factor
$V(\l)=\phi_-(\l)\phi_+(\l)$ so that $\phi_-(\l)$ is lower-triangular and $\phi_+(\l)$ is
upper-triangular with $1$ on the diagonal: 
$$ V(\l)= \phi_-(\l)\phi_+(\l)= \pmatrix{1&0\cr
1&-2\l\cr}\pmatrix{1&\l\cr 0&1\cr}.$$ Note that $\phi_+$ is exactly the $\phi$ used to
define the reality condition for the KdV-hierarchy.  This gives us a hint for the choice of the
$GD_n$-reality condition.  We proceed as follows: The vacuum for the
$GD_n$-hierarchy is the order $n$ operator
$L_0=D^n-\l^n$, and the corresponding first order system is 
${d\over dx} + A_0(\l)$, where $A_0(\l)=b+e_{1n}\l^n$.  It is
easy to check that
$$V(\l)=\pmatrix{1&\l & \cdots &\l^{n-1}\cr
1&\a\l&\cdots&(\a\l)^{n-1}\cr 1&\a^2\l&\cdots&(\a^2\l)^{n-1}\cr &&\cdots&\cr
1&\a^{n-1}\l&\cdots &(\a^{n-1}\l)^{n-1}\cr} =\left((\a^{i-1}\l)^{j-1}\right)$$
diagonalizes $A_0(\l)$, where $\a= e^{2\pi i\over n}$. In fact,
$$V(\l)A_0(\l)V(\l)^{-1} = a\l, \quad a = \diag(1,\a,
\cdots, \a^{n-1}).$$   

\refclaim[pc] Proposition.  $V(\l)$ can be factored uniquely as
$$V(\l)=\phi_n^-(\l)\phi_n(\l),$$ where
$\phi_n^-$ is  lower-triangular and $\phi_n$ is upper-triangular  
with $1$'s on the diagonal.

\proof 
It is an elementary result in linear algebra
that the factorization of $V=\phi_n^-\phi_n$ can be
carried out using the Gaussian elimination  if the all the principal $k\times k$
minors $\D_k$ of $V=(v_{ij})$ are non-zero.  But
$$\D_k=\det((v_{ij})_{1\leq i, j\leq k})= \l^k \prod_{0\leq i< j\leq
k-1}(\a^j-\a^i),$$ which is not zero for $\l\not=0$. \qed

 For example, we use Gaussian elimination to factor $V$ and get 
$$\eqalign{\phi_3(\l)&= \pmatrix{1&\l&\l^2\cr 0&1&(1+\a)\l\cr
0&0&1\cr},  \quad \a=e^{2\pi i \over 3}\cr 
\phi_4(\l) &=\pmatrix{1&\l&\l^2&\l^3\cr 0&1&(1+\a)\l&(1+\a+\a^2)\l^2\cr
0&0&1&(1+\a+\a^2)\l\cr 0&0&0&1\cr}, \quad \a=e^{2\pi i\over 4} =i.\cr}$$
However,  it is difficult to write down an explicit formula for $\phi_n$ using
Gaussian elimination. In order to do this, we need first prove some properties of
$\phi_n$. 
  Let  $\cg(k)$ denote the subspace spanned by
$$\{e_{i,i+k}\n 1\leq i, i+k\leq n\}.$$ Then 
\refeq[pg]$$\cg(k) \cg(m)\subset \cg(k+m).$$
 Let $\cn_+=\sum_{i=1}^{n-1} \cg(i)$ denote the
subalgebra of strictly upper triangular matrices, $\cb_+=\sum_{i=0}^n
\cg(i)$ the subalgebra of upper triangular matrices, and $N_+, B_+$ the Lie
group associated to $\cn_+, \cb_+$ respectively.  

\refclaim[pd] Proposition.  Let $\phi_n(\l)$ be the polynomial obtained in
Proposition \refpc{}.  Then
\item {(i)} $\phi_n= I +\sum_{i=1}^{n-1}f_i\l^i$ for some constant $f_i\in \cn_+$,
\item {(ii)} $\phi_n(\l)(e_{1n}\l^n+ b)\phi_n(\l)^{-1} =a \l+b$.

\proof  The Gaussian elimination  proves  (i).  

Since $V(\l)(e_{1n}\l^n + b)V(\l)^{-1}= a\l$ and 
$V(\l)=\phi_-(\l)\phi_n(\l)$, 
$$\phi_n(\l)(e_{1n}\l^n+ b)\phi_n(\l)^{-1} = 
\phi_-(\l)^{-1}a \l \phi_-(\l).$$ 
Note that the left hand side lies in $\sum_{k\geq -1}\cg(k)$ and
the right hand side lies in $\sum_{k\leq 0}\cg(k)$. Moreover, the $\cg_0$-component of the
right hand side is $a\l$ and the 
$\cg_{-1}$-component of the left hand side is  $b$. This proves (ii).
\qed

 In the following,  we  use  Proposition
\refpd{} (ii) to get an explicit formula for $\phi_n$.   We need a Lemma,
which is proved by a direct computation and \refpg{}. 

\refclaim[nt] Lemma.  Let $b=e_{21} + e_{32} +\cdots + e_{n n-1}$. Then
\item {(1)} $[b,\cg(i)]\subset \cg(i-1)$ and $\ad(b)$ is injective on $\cn_+$,
\item {(2)} $[b,x]\in \cg(i-1)$ if and only if $x\in \cg(i)$,
\item {(3)} if  $[b,x]=\sum_{k=1}^{n-i+1} c_k e_{k,k+i-1} \in \cg(i-1)$ for
$i\geq 1$, then $\sum_{k=1}^{n-i+1} c_k=0$ and $x=-\sum_{k=1}^{n-k}
(\sum_{j=1}^k c_j)e_{k, k+i}$.\ei

\refclaim[ku] Proposition.  Suppose
$\phi_n(\l)= I + f_1\l +
\cdots + f_{n-1}\l^n$ is a $N_+$-valued map of degree
$n-1$ in $\l$.  Then 
\refeq[na]$$\phi_n(\l)(e_{1n}\l^n + b) = (a\l+b) \phi_n(\l)$$ if and only if
\refeq[nr]$$f_i= (1+ \a+ \cdots +\a^{i-1})^{-1}\L^i,\qquad {\rm where\,}$$  
\refeq[ns]$$\L= \sum_{i=0}^{n-1} (1+\a+\cdots +\a^{i-1}) e_{i,i+1}.$$

The verification that this formula gives a solution is quite tedious. However, it
helps to know there is a unique solution $\phi_n$ for equation \refna{}.  It is also
helpful to note that, after doing the computation for $n=3,4$ by the Gaussian
elimination, that $f_i\in \cg(i)$ and $f_if_j=f_jf_i$.  Hence we expect
$f_i=c_if_1^i$ for some constant $c_i$.  

\ms

\proof Compare coefficients of $\l^i$ in equation \refna{}
$$(I+f_1\l + \cdots + f_{n-1}\l^{n-1}) (e_{1n}\l^n + b) = (a\l + b) 
(I+f_1\l + \cdots + f_{n-1}\l^{n-1})$$
 for $1\leq i\leq n$ to get 
\refeq[mz]$$\cases{f_1b=bf_1 + a, &\cr  f_i b=bf_i + f_{i-1},& if $2\leq i\leq
n-1$,\cr e_{1n}= af_{n-1}.&\cr}$$ 
Since $a=\sum_{i=1}^n \a^{i-1}e_{ii}$ and $\sum_{i=0}^{n-1} \a^i=0$ and the
first equation of system \refmz{} is $[b,f_1]= -a$, Lemma
\refnt{} (3) implies that  
$$f_1 =\L =\sum_{i=1}^{n-1}(1+\a + \cdots +
\a^{i-1}) e_{i,i+1}\in \cg(1).$$  In particular $f_1\in \cg(1)$.  The second to
the $(n-1)$-th equation of \refmz{} is 
$[f_i, b]=f_{i-1}$.  So Lemma \refnt{} (2) and induction imply that $f_i\in
\cg(i)$ for $1\leq i\leq n-1$. 

Since $\ad(b)$ is injective on $\cn_+$, system \refmz{} has at most one
solution. Proposition \refpd{} shows that the overdetermined system \refmz{} has
a solution.  So system \refmz{} has a unique solution. 

It remains to prove the formula for $f_j$. 
We need the following simple equalities:
\refeq[nu]$$\L b - b\L = a,$$
\refeq[nv]$$\L a= \a a\L.$$
(They can be proved by a direct computation.)
Next we claim that 
\refeq[nw]$$\L^k b- b\L^k = (1+\a + \cdots + \a^{k-1})a\L^{k-1}$$
is true for $1\leq k\leq n-1$. 
This equality implies that $f_k=(1+\a + \cdots + \a^{k-1})^{-1} \L^k$ solves
system \refmz{}.   We use induction to prove equality \refnw{}. When $k=1$,
equality \refnw{} is \refnu{}.  Now suppose equality \refnw{} is true for $k$.
Then 
$$\eqalign{&\L^{k+1}b - b\L^{k+1} = \L(\L^k b) - b\L^{k+1}\cr
&\, =\L\left(b\L^k + \left(\sum_{i=0}^{k-1}\a^i\right)a\L^{k-1}\right)-
b\L^{k+1}\cr &\, = (b\L + a)\L^k + \left(\sum_{i=0}^{k-1}\a^i\right)\L
a\L^{k-1} - b\L^{k+1}\cr &\,= b\L^{k+1} + a\L^k +\left(
\sum_{i=0}^{k-1}\a^i\right)
\a \,a \L \L^{k-1} - b\L^{k+1}\cr &\, = (1+\a + \cdots + \a^k) a\L^k.\cr}$$
This completes the proof. \qed

Since $f_i = c_i \L^i$ with $c_i=(\sum_{k=0}^{i-1}\a^k)^{-1}$, we have

 \refclaim[nk] Corollary. Let $f_i\in
\cg(i)$ be as in Proposition \refku{}. Then  $f_if_j=f_jf_i$ and $f_{1, i+1}=1$
for all $1\leq i, j\leq n-1$.

Next consider the $GD_n$-{\it reality condition\/}:
\refeq[kc]$$\phi_n(\l)^{-1}A(\l)\phi_n(\l)= \phi_n(\a \l)^{-1}A(\a \l)
\phi_n(\a\l),\quad
\a=e^{2\pi i\over n}.$$ The following statements are easily seen to be
equivalent:
\item {(i)} $A$ satisfies the $GD_n$-reality condition,
\item {(ii)} all entries of $\phi_n(\l)^{-1}A(\l)\phi_n(\l)$ are polynomial in
$\l^n$,
\item {(iii)} $A$ is a fixed point of the order $n$ automorphism
$$\Ad(\phi_n)\circ \tau_n \circ \Ad(\phi_n)^{-1},$$ where $\tau_n(A)(\l)=
A(\a\l)$. 

\ms

The first step in the construction of B\"acklund transformations is to
determine the condition for
$a\l+Y$ to satisfy the $GD_n$-reality condition. We need  two Lemmas:

\refclaim[mo] Lemma.  $\phi_n(\l)^{-1}= I+g_1 \l + \cdots + g_{n-1}\l^{n-1}$,
where
$g_i\in \cg(i)$ are constant.

\proof Note that 
\refeq[mr]$$(I+g_1\l + \cdots + g_{n-1}\l^{n-1})(I+f_1\l + \cdots +
f_{n-1}\l^{n-1})=I$$ holds if and only if the coefficients of $\l^j$ are zero for
$1\leq j\leq 2(n-1)$.  Since 
$\cg(i)\cg(j)\subset \cg(i+j)$ and $\cg(m)=0$ if $m\geq n$,  the coefficient
of $\l^j$ in
\refmr{} is zero for all $j\geq n$. The coefficient of $\l$ in \refmr{} is zero
implies that 
$g_1=-f_1\in \cg(1)$.  The coefficient of $\l^j$ is zero implies that 
$$-g_j = f_j + \sum_{i=1}^{j-1} g_i f_{j-i}.$$ By induction on $j$, we conclude
that $g_j\in
\cg(j)$. \qed

\refclaim[nq] Lemma.  Let $\L=\sum_{i=1}^{n-1}(1+\a +\cdots + \a^{i-1})
e_{i, i+1}\in \cg(1)$ be as in Proposition \refku{}. Then the centralizer 
$$gl(n)_\L=\{Z\in gl(n)\n Z\L= \L Z\} = \sum_{i=0}^{n-1} C\L^i.$$ 

\proof Let $c_i=1+\a + \cdots + \a^{i-1}$. Write $Z=\sum_{1\leq
i,j\leq n} z_{ij} e_{ij}$.  Then
$Z\L=
\L Z$ if and only if $z_{i, j-1}c_{j-1} = c_i z_{i+1, j}$.  Hence 
$$z_{i, i+j} = \left(\prod_{k=1}^{i-1} {c_{k+j}\over c_k}\right) z_{1, 1+j}.$$ So
the dimension of $gl(n)_\L$ is $n$.  But $I, \L, \cdots, \L^{n-1}\in gl(n)_\L$ and
are linearly independent.  Hence $gl(n)_\L= \sum_{i=0}^{n-1}C\L^i$.  \qed

\refclaim[kd] Proposition.  $A(\l)=a\l+Y$ satisfies  the $GD_n$-reality
condition if and only if $$Y=b+ \sum_{i=0}^{n-1} y_i f_i,$$ where $f_i$'s are
defined in Proposition
\refku{}, $f_0=I$ and $b= e_{21} +e_{32} + \cdots + e_{n n-1}$. 

\proof  Suppose $a\l+Y$ satisfies the $GD_n$-reality condition. By
Proposition \refmo{},
$\phi_n(\l)^{-1}(a\l+ Y)\phi_n(\l)$ is a polynomial in $\l$ with degree $\leq
2n-1$ with constant term $Y$. But the $GD_n$-reality condition implies that it is a
polynomial in
$\l^n$.  So   
$$\phi_n(\l)^{-1}(a\l+ Y)\phi_n(\l) = C_0 \l^n + Y$$ for some $C_0\in gl(n)$. Write
$Y=b+Z$.  It follows from Proposition \refku{} that we have
$$\eqalign{&\phi_n(\l)^{-1}(a\l+Y) \phi_n(\l) = \phi_n(\l)^{-1}(a\l + b+
Z)\phi_n(\l)\cr  &=e_{1n}\l^n + b + \phi_n(\l)^{-1}Z\phi_n(\l).\cr}$$  Hence
$$\phi_n(\l)^{-1}Z\phi_n(\l)= C\l^n + Z,$$ where $C= C_0-e_{1n}$.  So
$Z\phi_n(\l)= \phi_n(\l)(C\l^n + Z)$, i.e.,
$$Z(I+f_1\l + \cdots + f_{n-1}\l^{n-1}) = (I+f_1\l + \cdots +
f_{n-1}\l^{n-1})(C\l^n + Z).$$ Because the left hand side has degree $n-1$ in
$\l$, the coefficient of $\l^n$ of the right hand side is zero.  This implies
that $C=0$. So we have
$$\phi_n(\l)^{-1}Z\phi_n(\l)=Z,$$ i.e.,
$$Z(I + f_1\l + \cdots + f_{n-1}\l^{n-1}) = (I+f_1\l + \cdots +
f_{n-1}\l^{n-1})Z.$$ 
Compare coefficient of $\l^i$ in the above equation to get
$f_jZ=Zf_j$ for all $1\leq j\leq n-1$.  But recall that $f_i = c_i \L^i$ for
some non-zero constant $c_i$. Hence $Z\L=\L Z$.  By Lemma \refnq{}, 
$Z=\sum_{i=0}^{n-1} y_i f_i$ for some constant $y_0, y_1, \cdots, y_{n-1}$.
\qed

For example, for $n=3,4$, $a\l+Y_n$ satisfies the $GD_n$-reality condition
if and only if 
$$\eqalign{Y_3&=\pmatrix{y_0&y_1&y_2\cr 1&y_0&(1+\a)y_1\cr
0&1&y_0\cr}, \quad
\a=e^{2\pi i\over 3}\cr Y_4 &=\pmatrix{y_0&y_1&y_2&y_3\cr
1&y_0&(1+\a)y_1&(1+\a+\a^2) y_2\cr 0&1&y_0&(1+\a+\a^2)y_1\cr
0&0&1&y_0\cr}, \quad \a=e^{2\pi i\over 4}.\cr}$$

\refpar[pe] Definition.  For $y=(y_0, \cdots, y_{n-1})$, we set
$$\eqalign{&Y_y=b+\sum_{i=0^n}y_if_i,\cr
&f_y(\l)= a\l + Y_y,\cr},$$ where 
$f_0=I$ and $f_i$'s are given in Proposition \refku{}.

\refpar[nn] Definition.  Let $\cs_{\phi_n}$ denote the space of Schwartz maps
$u:R\to sl(n)_a^\perp$ such that $ A(\l)=a\l+b+u$ satisfies the $GD_n$-reality
condition, where $b=e_{21}+\cdots + e_{n n-1}$.  In other words, $\cs_{\phi_n}$ is
the space of
$u=Y_q$ for some smooth Schwartz map $q=(0, q_1, \cdot, q_{n-1})$. 

\refclaim[ks] Proposition. Given $u:R\to sl(n)_a^\perp$, if $a\l+ u$ satisfies
the
$GD_n$-reality condition, then the Lax pair of the $(nj+1)$-th flow in the
$sl(n,C)$-hierarchy satisfies the $GD_n$-reality condition. 

\proof  The proof is similar to that of Proposition \reffb{} except we
replace $-\l$ by $\a\l$, where $\a=e^{{2\pi i\over n}}$.  So we have 
$$\phi_n(\a\l)\phi_n(\l)^{-1} Q(\l) \phi_n(\l)\phi_n(\a\l)^{-1} = \a Q(\a\l),$$
i.e., $\phi_n(\l)^{-1}Q(\l) \phi_n(\l)= \a \phi_n(\a\l)^{-1}Q(\a\l)\phi_n(\a\l)$. 
This implies that $$\l^{nj+1}\phi_n(\l)^{-1}Q(\l)\phi_n(\l)$$ is a power
series in $\l^n$. In particular, 
$$\eqalign{&\phi_n(\l)^{-1}(\l^{nj+1}Q(\l))_+\phi_n(\l)\cr &\,\,=
\phi_n(\l)(a\l^{nj+1} + u\l^{nj} + Q_2\l^{nj-1} +\cdots +
Q_{nj+1})\phi_n(\l)\cr}$$ is a polynomial in $\l^n$. 
\qed

As a consequence we have

\refclaim[ke] Corollary.   The $jn+1$-th flow in the $sl(n,C)$-hierarchy
leaves the submanifold $\cs_{\phi_n}$ invariant for all $j\geq 0$. 

The first order $n\times n$ system 
$${d\over dx}(y_1, \cdots, y_n)=(y_1, \cdots, y_n)(a\l+Y_q)$$ is equivalent to a
unique  order $n$ differential operator of $y_1$.  Hence $\cs_{\phi_n}$ is
isomorphic to the phase space  $\cp_n$ of the
$GD_n$-hierarchy, and the restriction of the $sl(n,C)$-hierarchy to
$\cs_{\phi_n}$ corresponds to the $GD_n$-hierarchy.  It will still be called the
$GD_n$-hierarchy.  

\refpar[ms] Example.  $\cs_{\phi_3}$ is the space of 
$u= q_1 f_1 +q_2 f_2$.  A direct computation shows that the third order
scalar differential operator corresponding to the first order system 
${dy\over dx} - y(a\l+ u)=0$ is 
$$D^3 -((1-\a^2)q_1 D + ((q_1)_x+ q_2)) = \l^3.$$ Hence 
$$p_1= (1-\a^2)q_1, \quad p_2= (q_1)_x + q_2$$  defines a linear
isomorphism from $\cs_{\phi_3}$ to  $\cp_3$, and the inverse is given by 
$$q_1= {p_1\over 1-\a^2}, \quad q_2= p_2 - {(p_1)_x\over 1-\a^2}.$$

\ms
The general outline for constructing B\"acklund and Darboux transformations is the
same as in our previous examples. But the computations are quite involved.  We
give an admittedly rather brief description of how the construction goes. 

\refpar[po] Definition. 
Let $G_-^{GD}$ denote the group of rational maps $f:S^2\to GL(n,C)$ such that
\item {(i)} $f(\infty)=I$,
\item {(ii)} there exists a rational map $g:C\to C$ such that $gf$ satisfies the
$GD_n$-reality condition. 

\ms
$G_-^{GD}$ acts on the space of local solutions of the $(n+1)$-th flow in the
$GD_n$-hierarchy, and the action of a linear fractional map in $G^{GD}$ gives a
B\"acklund transformation.  
 To construct B\"acklund transformations for the
$GD_n$-equation, we  start with a  degree one rational  map $\o_y(\l)={a\l
+Y_y\over \l-k}$, where
$-k$ is an $n$-th root of $\det(Y_y)$ (i.e., 
$\det(Y_y)=(-k)^n$). Note that although $\o_y$ does not satisfy the
$GD_n$-reality condition,
$(\l-k)\o_y(\l)=a\l+Y_y$ does.  So $\o_y\in G_-^{GD}$.

\refclaim[ph] Proposition.   $f_y(\l)=a\l + Y_y$ has the following properties:
\item {(i)} $\det(f_y(\l))= (-1)^n(\l^n-k^n)$, where
$\det(Y_y)= (-k)^n$.  
\item {(ii)} $f_y(\l)^{-1}= {h(\l)\over \l^n-k^n}$ for some $gl(n)$-valued 
polynomial of degree $(n-1)$.  

\proof  Since $f_y(\l)$ satisfies the $GD_n$-reality condition,
$$\det(f_y(\l))=\det(f_y(\a
\l))$$ for $\a=e^{2\pi i/n}$. So $\det(f_y(\l))$ is a polynomial in $\l^n$. But
$\det(f_y(\l))$ is of degree $n$ in $\l$ whose leading term is
$$\prod_{i=0}^{n-1}\a^i = \a^{{n(n-1)\over 2}}= (-1)^{n+1}$$ and the constant term
is $\det(Y_y)$.  So
$\det(f_y(\l))= (-1)^{n+1}(\l^n-k^n)$ for some $k\in C$.  This proves (i).  Statement (ii) follows
from the Cramer's rule. 
\qed

Let $u$ be a local solution of $n+1$-th flow in the $GD_n$-hierarchy that admits a
reduced wave function
$m$, and $E$ the trivialization
 of $u$. Suppose at each
$(x,t)$ we can find  $\tilde y(x,t)$ such that 
\item {(i)} $\det(Y_{\tilde y(x,t)})= (-k)^n$, and 
\item {(ii)} $\tilde E(x,t,\l)= {a\l+Y_y\over \l-k} E(x,t,\l)\left({a\l+ Y_{\tilde
y(x,t)}\over
\l-k}\right)^{-1}$ is holomorphic in $\l\in C$. 
 
\ni Then we can proceed the same way  as in the proof
of Theorem \refhma{} to conclude that
$\tilde E$ is the trivialization of some local solution of the $(n+1)$-th flow $\tilde u$ in the
$sl(n,C)$-hierarchy. First notice that the denominators in $\tilde E$ can be
canceled. So we get
\refeq[pf]$$\tilde E(x,t,\l)=
(a\l+Y_y)E(x,t,\l)(a\l + Y_{\tilde y(x,t)})^{-1}.$$ Since all three terms in the right
hand side of \refhh{} satisfy the $GD_n$-reality condition, $\tilde E$ also satisfies
the $GD_n$-reality condition. Hence $\tilde E$ corresponds to a new solution
$\tilde u$ of the $(n+1)$-th flow in the $GD_n$-hierarchy.  Proposition \refmp{}
implies that the image of
$(ak+Y_y)^{-1}$ is the kernel of $(ak+Y_y)$.  But formula
\refpf{} is holomorphic in $\l\in C$ implies that the residue of the right hand side
at $\l=k$ is zero.  Hence 
$${\rm Ker\/}(ak + Y_y) = E(x,t,k)(\Im(ak+Y_{\tilde y(x,t)})^{-1}) 
=E(x,t,k)(\Ker(ak+Y_{\tilde y(x,t)}).$$ Therefore, we need to find the relation
between the zeros and kernels of 
$a\l+Y_y$ and $y$. We do this in the following few Propositions.  

\refpar[kf] Definition.  Let $C:C^n\to gl(n)$ denote the map defined by
 $C(v)=$ the matrix whose first column is $v$ and whose $i+1$-th column is
$\phi_n(\a^ik)\phi_n(k)^{-1}(v)$ for  $1\leq i\leq n-1$.

Let $\ell_n:C^n\to C$ denote the projection onto the $n$-th coordinate. 

\refclaim[kg] Proposition.  Suppose $\det(f_y(\l))= (-1)^{n+1}(\l^n-k^n)$ and
$0\not=v\in C^n$ such that 
$f_y(k)(v)=0$.  Then
\item {(i)} $f(\a^jk)v_j=0$, where $v_j=\phi_n(\a^jk)\phi_n(k)^{-1}v$
for $1\leq j\leq n-1$,
\item {(ii)} if $k\not=0$, then $C(v)$ is non-singular and
$Y_y=-kaC(v)aC(v)^{-1}$. 

\proof 
Since $\phi_n(k)^{-1}f(k)\phi_n(k) = \phi(\a k)^{-1}f(\a k) \phi_n(\a k)$, (i)
follows. 

Use (i), we get $(\a^j ka + Y_y)v_j=0$, so $Y_yv_j= -\a^jkav_j$.  Write this in
terms of matrix to get $Y_yC(c)= -kaC(v)a$,  which proves (ii). 
\qed

Let $k\in C$, and $v\in C^n$ a non-zero vector.  Set  
\refeq[pk]$$h_{k,v}(\l)=
a\l-kaC(v)aC(v)^{-1}.$$ The above Proposition says that if $\det(Y_y)=(-k)^n$ and
$(ak+Y_y)v=0$, then
$$f_y(\l)= h_{k,v}(\l).$$   Now given any $k\in
C$ and $v\in C^n$, does $h_{k,v}$ satisfies the $GD_n$-reality condition? We will
answer this next.  

Let $\cy_n$ denote the set of all $y=(y_0,y_1, \cdots, y_{n-1})\in C^n$ such
that $$\det (Y_y)=\det(b+\sum_{i=0}^{n-1} y_if_i)\not=0,$$ and
$$\eqalign{\D&=\{k\in C\n 0\leq {\rm arg\/}(k) <{2\pi \over n}, k\not=0\}, \cr
V_n&=\{(r_1,\cdots, r_{n-1}, 1)^t\n r_i\in C\}.\cr}$$ 

\refclaim[nx] Proposition.  Let $K_n:\cy_n\to \D\times V_n$ be the map defined by  
$K_n(y)= (k, v)$, where $k\in \D$ such that
$\det(Y_y)=(-k)^n$ and  $(ak+Y_y)v=0$. Then
\item {(i)} $K_n$ is bijective,
\item {(ii)} both $K_n$ and $K_n^{-1}$ are algebraic maps. 

\proof  Proposition \refkg{} implies that $K_n$ is one to one.  

Let $k\in C$ be non-zero, and $v=(r_1, \cdots, r_{n-1}, 1)^t$.  
To prove $K_n$ is onto is
equivalent to prove the following linear system  has a non-zero solution $y$:
\refeq[ny]$$(ak +b+ \sum_{i=0}^{n-1} y_if_i) v=0.$$
Write
$f_i=\sum_{k=1}^{n-i} c_{k,k+i} e_{k,k+i}$.  We claim that 
$c_{k,k+i}\not=0$ for all $1\leq k\leq n-i$. 
To see this, we recall that
$f_i= s_i^{-1}
\L^i$ and
$\L=\sum_{i=1}^{n-1} s_i e_{i, i+1}$, where $s_i=\sum_{k=0}^{i-1}\a^k$.
Since  $s_1, \cdots, s_{n-1}$ are non-zero,  our claim is proved. 
System \refny{} in matrix form is
\refeq[nya]$$\pmatrix{k+y_0&c_{11}y_1&\cdot &\cdot &\cdot &\cdot
&c_{1n}y_{n-1}\cr
1& \a k+y_0&&&&&\cr
0&1&&&&&\cr
&&&&&&\cr
&&&&&&\cr
&&&&1&\a^{n-2}k + y_0&c_{n-1,n}y_1\cr
&&&&&1&\a^{n-1}k + y_0\cr}\pmatrix{r_1\cr r_2 \cr \cdot\cr\cr\cr
r_{n-1}\cr1\cr}=
\pmatrix{0\cr 0\cr \cr\cr\cr\cr 0\cr}.$$ 
The last equation in \refnya{} implies $y_0=-(r_{n-1} + \a^{n-1}k)$.
Substitute $y_0$ to the $(n-1)$-th equation to get
$$y_1= -{1\over c_{n-1, n}}(r_{n-2} +(\a^{n-2}k -r_{n-1} -\a^{n-1}k)
r_{n-1}).$$ The $(n-j)$-th equation gives a recursive formula for $y_j$ in terms of
$y_1, \cdots, y_{j-1}$. This is due to the fact that the matrix is the sum of
$b=e_{21}+e_{32} + \cdots + e_{n n-1}$ and a upper triangular matrix. So $y$ is
solved explicitly in terms of $k,v$. 
 \qed

Using the same kind of arguments as for the KdV- and KW-hierarchies
we get  local Darboux and B\"acklund transformations.  We give an outline of the
results.  Let $q=(0, q_1, \cdots, q_{n-1})$, and $$u=Y_q= b+ \sum_{i=1}^{n-1}
q_i f_i$$  a solution of the
$n+1$-th flow in the
$GD_n$-hierarchy on $\cs_{\phi_n}$, and $E$ the trivialization of $u$
normalized at
$(0,0)$.  Let $(k,v)\in \D\times V_n$.  Set 
$$\tilde v(x,t)=E(x,t,k)^{-1}(v).$$  Since $\tilde v(0,0)=v$, 
 there exists an open neighborhood $\co$ of
$(0,0)$ in $R^2$ such that $\ell_n(\tilde v(x,t))\not=0$ for all $(x,t)\in \co$.  Set
 $$\tilde y(x,t)=K_n^{-1}\left(k, {\tilde v(x,t)\over \ell_n(\tilde v(x,t))}\right).$$
Then

 \item {(i)} $\tilde u = aua^{-1} +
[Y_{\tilde y(x,t)}, a]a^{-1}$ is again a solution of the $n+1$-th flow in the
$GD_n$-hierarchy.  In other words, $\tilde u= Y_{\tilde q}$, where $\tilde q=(0,
\tilde q_1, \cdots, \tilde q_{n-1})$ and 
$$\tilde q_i=\a^{-i} q_i + (1-\a^{-i}) y_i, \qquad 1\leq i\leq n-1.$$
Let
$$\o_{k,v}= {h_{k,v}\over \l-k}, \quad \tilde u= \o_{k,v}\ast u, \quad \tilde q=
\o_{k,v}\ast q,$$ where $h_{k,v}$ is defined by formula \refpk{} ($h_{k,v}=f_y$ if
$y=K_n^{-1}(k,v)$).  This gives the Darboux transformation for the $(n+1)$-th flow
in the
$GD_n$-hierarchy.

\item {(ii)} $Y_{\tilde y(x,t)}$ is a solution of 
$$\cases{Y_x= Y u - (aua^{-1} + [Y, a]a^{-1}) Y,&\cr  Y_t = YQ_{a,n+1}(u)
-Q_{a,n+1}(aua^{-1} + [Y, a]a^{-1}) Y,\cr}$$ where
$Q_{a,n+1}$ is the polynomial differential operator defined in the
$sl(n,C)$-hierarchy.   This gives the B\"acklund transformation. 

\refpar[pj] Example.  The trivialization $E_0$ of the vacuum $u=0$ is 
$$E_0(x,t,\l)= e^{a(\l x+ b) + a\l^{n+1} t},$$ where $b=e_{21} +\cdots + e_{n, n-1}$. 
 Given $(k,v)\in
\D\times V_n$,  set $$\eqalign{\tilde v(x,t)&= e^{-a(k x+ b) - ak^{n+1} t}v,\cr
\tilde y(x,t) &=K_n^{-1}\left(k,{\tilde v(x,t)\over \ell_n(\tilde v(x,t))}\right),\cr
q_j&= (1-\a^{-j})\tilde y_j.\cr}$$ 
Then $u=Y_q$ is a solution of the $j$-th flow in the $GD_n$-hierarchy.  

We obtain the following relation among simple elements:

\refclaim[kn] Proposition.   Let $(k_1,v_1),(k_2,v_2)\in \D\times V_n$ such that
$k_1^n\not= k_2^n$. Then
$$h_{k_2,\xi_2}h_{k_1,v_1}=h_{k_1,\xi_1}h_{k_2,v_2}$$ if and only if 
$\xi_1$ is parallel to $ h_{k_2,v_2}(k_1)(v_1)$ and $\xi_2$ is parallel to
$h_{v_1,k_1}(k_2)(v_2)$.

\proof Write $h_{k_i, v_i}= a\l + Y_i$ and $h_{k_i, \xi_i}= a\l + Z_i$. If 
$$(a\l+Z_2)(a\l + Y_1) = (a\l+ Z_1)(a\l+ Y_2),$$ then set $\l=k_2$ to get 
$$(ak_2+Z_2)(ak_1+Y_1)(v_2) = (ak_2+Z_1)(ak_2+ Y_2)(v_2).$$  But the
right hand side is zero by definition of $Y_{v,k}$.  So $(ak_2 +Y_1)(v_2)$ lies
in the kernel of $(ak_2+Z_2)$, which is $C\xi_2$.  This implies that  $\xi_2$ 
is parallel to
$(ak_2 +Y_1)(v_2)$.  Similarly, $\xi_1$ is parallel to $(ak_1+Y_2)(v_1)$.
\qed

We give the Permutability formula next.  Let $u_0=Y_{q^{(0)}}$
be a solution of the $(n+1)$-th flow in the
$GD_n$-hierarchy, and $E(x,t,\l)$ the trivialization of $u_0$ normalized at $(0,0)$.
Let $(k_1,v_1), (k_2, v_2)\in \D\times V_n$ such that $k_1^n\not=k_2^n$, and 
$$u_i= \o_{k_i, v_i}\ast u_0, \quad q^{(i)}= \o_{k_i,v_i}\ast q^{(0)}.$$ 
Set $$\eqalign{\tilde v_i(x,t)&= E(x,t,k_i)^{-1}(v_i), \quad i=1,2,\cr
\xi_2(x,t)&= h_{k_1,\tilde v_1(x,t)}(k_2)(\tilde v_2(x,t)),\cr
y^{(i)}(x,t)&= K_n^{-1}\left(k_i, {\tilde v_i(x,t)\over \ell_n(\tilde
v_i(x,t))}\right)\quad i=1,2,\cr
 q_j^{(i)}&= \a^{-j}q_j^{(0)} + (1-\a^{-j})y_j^{(i)}, \quad i=1, 2, \,{\rm and\,}\, 1\leq
j\leq n-1\cr
\tilde y^{(2)}(x,t)&=K_n^{-1}\left(k_2, {\xi_2(x,t)\over
\ell_n(\xi_2(x,t))}\right).\cr}$$ Set $$q_j^{(3)}=\a^{-j} q_j^{(1)} + (1-\a^{-j}) \tilde
y_j^{(2)}.$$ Then 
\item {(i)} $u_3= Y_{q^{(3)}}$ is a solution of the $GD_n$-hierarchy,
\item {(ii)} $u_3= (\o_{k_2,\xi_2(0,0)}\o_{k_1, v_1})\ast u_0$.  
\item {(iii)} since $y^{(i)}$ can be written as an algebraic function of $q^{(0)}$ and
$q^{(i)}$,  $q^{(3)}$ can be written as an algebraic
function of
$q^{(0)}, q^{(1)}, q^{(2)}$; this is the  permutability formula for the
$GD_n$-hierarchy. 

\vfil\eject

\Bibliography

\b //AC//Ablowitz, M.J.,Clarkson, P.A.//Solitons,
non-linear evolution equations and inverse scattering//\break Cambridge Univ.
Press////1991//////

\a //AKNS//Ablowitz, M.J., Kaup, D.J., Newell, A.C. and Segur, H.//The inverse
scattering transform - Fourier analysis for nonlinear
problems//Stud. Appl. Math.//53//1974//249-315////

\a //AM//Adler, M., Moser, J.//On a class of polynomials connected with the
Korteweg-deVries equation//Comm. Math. Phys.//61//1978//1-30////

\b //Ba//B\"acklund, A.V.//Concerning surfaces with constant negative
curvature//New Era Printing Co., Lancaster, PA////////original 1883, translation 1905//

\a //BC1//Beals, R., Coifman, R.R.//Scattering and inverse scattering for
first order systems//Commun. Pure Appl. Math.//37//1984//39-90////

\a //BC2//Beals, R., Coifman, R.R.//Inverse scattering and evolution
equations//Commun. Pure Appl. Math.//38//1985//29-42////

\b //BDT//Beals, R., Deift, P., Tomei, C.//Direct and inverse
scattering on the line//American Math. Soc.////1988//////

\a //BP//Bobenko, A., Pinkall, P.//Discrete surfaces with constant negative Gaussian
curvature and the Hirota equation//J. Differential Geometry//43//1996//527-611////

\b //Ch//Cherednik, I.V.//Basic methods of soliton theory, Adv. ser.  Math.
Phys. v. 25//World Scientific////1996//////

\b //Da//G. Darboux//Lecon sur la th\'eorie
g\'en\'erale des surfaces//Chelsea////1972////3rd edition//

\a //De//Deift, P.//Applications of a commutation formula//Duke. Math.
J.//45//1978//267-310////

\a //DS1//Drinfel'd, V.G., and Sokolov, V.V.//Equations of Korteweg-de Vries type
and simple Lie algebras//Dokl. Akad. Nauk SSSR//258//1981//11-16//(Trans. as
Soviet Math. Dokl. 23, 457-462//

\a //DS2//Drinfel'd, V.G., and Sokolov, V.V.//Lie algebras and equations of Korteweg-de Vries
type//Itogi nauki//24//1984//81-180////

\b //Ei//Eisenhart, L.P.//A treatise on the
differential geometry of curves and surfaces//Ginn////1909//////

\b //FT//Faddeev, L.D., Takhtajan, L.A.//Hamiltonian Methods in the theory
of Solitons//Springer-Verlag////1987//////

\a //FK//Fordy, A.P., Kulish, P.P.//Nonlinear Schr\"odinger equations and simple
Lie algebra//Commun. Math. Phys.//89//1983//427-443////

\a //GZ//Gu, G.H., Zhou, Z.X.//On Darboux transformations for soliton
equations in high dimensional spacetime//Letters in Math.
Phys.//32//1994//1-10////

\a //KW//Kuperschmidt, B.A., Wilson, G.//Modifying Lax equations and the second
Hamiltonian structure//Invent. Math.//62//1981//403-436////

\a //P//Pohlmeyer, K.//Integrable Hamiltonian systems and interactions through
quadratic constraints//Comm. Pure Appl. Math.//46//1976//207-221////

\b //PS//Pressley, A. and Segal, G. B.//Loop Groups//Oxford Science Publ., 
Clarendon Press, Oxford////1986//////

\a //Sa//Sattinger, D.H.//Hamiltonian hierarchies on semi-simple Lie
algebras//Stud. Appl. Math.//72//1984//65-86////

\a //SZ1//Sattinger, D.H.,Zurkowski, V.D.//Gauge theory of B\"acklund
transformations.I//Dynamics of infinite dimensional
systems, Nato Sci. Inst. Ser. F. Comput.
Systems Sci.//37//1987//273-300//Springer-Verlag//

\a //SZ2//Sattinger, D.H.,Zurkowski, V.D.//Gauge theory of B\"acklund
transformations.II//Physica//26D//1987//225-250////

\a //SS//Shatah, J., Strauss, W.//Breathers as homoclinic geometric wave
maps//Physics, D//99//1996//113-133////

\a //Te//Terng, C.L.//Soliton equations and differential
geometry//J. Differential Geometry//45//1997//407-445////

\p //TU1//Terng, C.L., Uhlenbeck, K.//Poisson actions and scattering theory
for integrable systems//////////preprint dg-ga 9707004//

\p //TU2//Terng, C.L., Uhlenbeck, K.//Homoclinic wave maps into compact
symmetric spaces//////////in preparation//

\a //U1//Uhlenbeck, K.//Harmonic maps into Lie group (classical solutions of the
Chiral model)//J. Differential Geometry//30//1989//1-50////

\a //U2//Uhlenbeck, K.//On the connection between harmonic maps and the self-dual
Yang-Mills and the sine-Gordon equations//Geometry \& Physics//2//1993//////

\a //ZMa1//Zakharov, V.E., Manakov, S.V.//On resonant interaction of wave packets
in non-linear media//JETP Letters//18//1973//243-247////

\a //ZMa2//Zakharov, V.E., Manakov, S.V.//The theory of resonant interaction of wave
packets in non-linear media//Sov. Phys. JETP//42//1974//842-850////

\a //ZS1//Zakharov, V.E., Shabat, A.B.//Exact theory of two-dimensional
self-focusing and one-dimensional of waves in nonlinear
media//Sov. Phys. JETP//34//1972//62-69////

\a //ZS2//Zakharov, V.E., Shabat, A.B.//Integration of non-linear equations of
mathematical physics by the inverse scattering method, II//Funct. Anal.
Appl.//13//1979//166-174////

\def\uhlenaddr{{\hsize=2in \vsize=1in
\hbox{\vbox{
\leftline{Karen Uhlenbeck}
\leftline{Department of Mathematics}
\leftline{The University of Texas at Austin}
\leftline{RLM8.100 Austin, Texas 78712}
\leftline{email:uhlen@math.utexas.edu}
}}}}

\def\cltaddr{{\hsize=2in \vsize=1in
\hbox{\vbox{
\leftline{Chuu-lian Terng}
\leftline{Department of Mathematics}
\leftline{Northeastern University}
\leftline{Boston, MA  02115}
\leftline{email: terng@neu.edu}
}}}}

\bigskip\bigskip
\leftline{\cltaddr \hfill \uhlenaddr}

\end